\documentclass[11pt]{ip-journal}
\usepackage{graphicx}
\usepackage{amsmath,amsthm}
\usepackage{amssymb}
\usepackage{caption}
\usepackage[title]{appendix}
\usepackage[all]{xy}
\usepackage{mathabx}
\usepackage{slashed}
\usepackage{enumerate}
\theoremstyle{plain}
\usepackage{bm}
\newtheorem{theo}{Theorem}
\newtheorem{utheo}{Upcoming Theorem}
\newtheorem{proposition}[theo]{Proposition}
\newtheorem{define}[theo]{Definition}
\newtheorem{lemma}[theo]{Lemma}
\newtheorem{corollary}[theo]{Corollary}
\DeclareRobustCommand{\gobblefour}[5]{}
\newcommand*{\SkipTocEntry}{\addtocontents{toc}{\gobblefour}}
\usepackage[colorlinks=true,urlcolor=blue]{hyperref}
\usepackage{xcolor}
\hypersetup{
    colorlinks,
    linkcolor={red!20!black},
    citecolor={blue!50!black},
    urlcolor={blue!80!black}
}
\usepackage{url}

\newcommand{\IR}{\mathbb{R}}
\newcommand{\tr}{\mathrm{tr}}
\newcommand{\IZ}{\mathbb{Z}}
\newcommand{\p}{\partial}
\newcommand{\e}{\mathbf{e}}
\newcommand{\D}{\mathpzc{D}}
\newcommand{\Fe}{\mathpzc{e}}
\newcommand{\FE}{\mathpzc{E}}
\renewcommand{\S}{\mathcal{S}}
\newcommand{\T}{\mathbb{T}}
\newcommand{\IC}{\mathbb{C}}
\newcommand{\IH}{\mathbb{H}}
\newcommand{\ev}{\mathrm{ev}}

\newcommand{\TT}{\mathcal{T}}

\renewcommand{\t}{\,t\!\!\!t}
\newcommand{\dds}{\frac{d}{ds}}
\newcommand{\ka}{\kappa}

\newcommand{\ii}{\mathrm{i}}
\newcommand{\tslash}{\slashed{t}}
\newcommand{\Mbar}{\overline{\mathcal{M}}}
\newcommand{\rr}{\bm{r}}
\newcommand{\Cbar}{\overline{C}}
\newcommand{\ebar}{\overline{\Fe}}
\newcommand{\Ff}{\mathpzc{f}}
\newcommand{\W}{\mathcal{W}}
\newcommand{\ch}{\mathrm{ch}\,}
\newcommand{\Ebar}{\overline{\mathcal{E}}}
\newcommand{\F}{\mathpzc{F}}

\DeclareMathAlphabet{\mathpzc}{OT1}{pzc}{m}{it}

\usepackage{tikz}
\tikzset{>=latex}
\usetikzlibrary{calc,decorations.markings,positioning}

\usetikzlibrary{arrows,positioning}
\usepackage{caption}
\usepackage{subcaption}

\begin{document}

\title[Instantons on multi-Taub-NUT Spaces II]{Instantons on multi-Taub-NUT Spaces II:\\
Bow Construction}

\author[S.A. Cherkis,   A. Larra\'in-Hubach, and  M. Stern]{Sergey A. Cherkis$^*$,  Andr\'es Larra\'in-Hubach$^{**}$,\\ and Mark Stern$^{\star}$}
\address{$^*$ School of Mathematics, Institute for Advanced Study, Princeton, NJ 08540, USA; 
Department of Mathematics, University of Arizona, Tucson, AZ 85721-0089, USA.}
\address{$^{**}$ Department of Mathematics, University of Dayton, Dayton, OH 45469, USA.}
\address{$\star$ Department of Mathematics, Duke University, Durham, NC 27708-0320, USA.}
\address{\small\tt cherkis@math.arizona.edu, alarrainhubach1@udayton.edu, }
\address{\small\tt  stern@math.duke.edu}

\begin{abstract}
Unitary anti-self-dual connections on Asymptotically Locally Flat (ALF) 
hyperk\"ahler spaces  are constructed in terms of  data organized in a 
bow.  Bows generalize quivers, and the relevant bow gives rise to the 
underlying ALF space as the  moduli space of its particular 
representation -- the small representation.  Any other representation of 
that bow gives rise to anti-self-dual connections on that ALF space.

We prove that each resulting connection has finite action, i.e. it is an instanton.  Moreover, we derive the asymptotic form of such a  connection and compute its topological class.
\end{abstract}

\maketitle

\newpage
\tableofcontents
\newpage

\section{Introduction}
\subsection{Perspective}
The complete construction of instantons on $\mathbb{R}^4$ discovered by Atiyah, Hitchin, Drinfeld, and Manin, and known as the ADHM construction  \cite{Atiyah:1978ri},  
has been generalized in several important ways: by Nahm to calorons, which are instantons on $\mathbb{R}^3\times S^1,$  and to monopoles on $\mathbb{R}^3$ \cite{Nahm:1983sv, Nahm:1979yw, NahmADHM}, and by Kronheimer and Nakajima to instantons on orbifolds and their deformations called ALE spaces \cite{KN}.  In each of these cases  a gauge equivalence class of a Hermitian connection with finite action functional and  anti-self-dual curvature is in bijective correspondence with certain algebraic or ODE solutions data,  called ADHM-Nahm data \cite{Hitchin:1983ay,KN}.  Correspondingly, in all of these cases there is a map, which we call the {\em Down transform}, mapping any instanton to its ADHM-Nahm data and another transform, which we call the {\em Up transform}, mapping any ADHM-Nahm data to an instanton.  The rationale for these names is that the ADHM,  as well as the Kronheimer-Nakajima transform, maps instantons in four dimensions to  some quiver data which is zero-dimensional.  The Nahm transform, in turn, relates calorons (instantons on four-dimensional $\mathbb{R}^3\times S^1$) or monopoles (in three dimensions) to the Nahm data on a circle or an interval.  Thus the Down transform maps a higher-dimensional object to a lower-dimensional one.  

The Up transform acts similarly, but in the opposite direction. In all of the known constructions, the Down and Up transforms are inverse to each other \cite{Hitchin:1983ay, Corrigan:1983sv, KN}.

There is a natural hyperk\"ahler structure on the instanton moduli space, as well as on the space of ADHM data.  The Down and Up transforms are not only bijections, but  are also isometries of these hyperk\"ahler manifolds.  For monopoles this was proved in \cite{Nakajima:1990zx}; while for instantons on ALE manifolds it was proved in \cite{KN}. 
(There is also a relation between instantons on a four-torus and those of its dual four-torus proved in \cite{Braam:1988qk}. This is one of the handful of cases, together with monopole walls \cite{Cherkis:2012qs,Cherkis:2014vfa,Cross:2015hla} and Hitchin systems on a plane \cite{Szabo07,Szabo17}, in which the dimension of the underlying space is the same on both sides of the transform, making our Up and Down nomenclature less felicitous.) 

A further generalization of the ADHM construction to instantons on Asymptotically Locally Flat (ALF) spaces was formulated in \cite{Cherkis:2010bn}.  In this case, the relevant instanton data is organized into a {\em bow}, which is essentially one-dimensional.  The construction of \cite{Cherkis:2010bn} presents the up transform, which is used in \cite{Cherkis:2009jm} to find explicit generic $SU(2)$ instanton connections on the Taub-NUT space.  The metric on the bow moduli space corresponding to one $SU(2)$ instanton on the Taub-NUT space is computed in \cite{Cherkis:2008ip}, where it is conjectured to be isometric to the instanton moduli space. 

\subsection{Main Results}
A generalization of the ADHM-Nahm construction, presented in \cite{Cherkis:2010bn}, transforms data associated to bows to Hermitian bundles with  compatible anti-self-dual  connections on multi-Taub-NUT spaces.  This {\em Up transform} maps (a gauge equivalence class of) a bow solution, which is a solution of a one-dimensional ODE, to (a gauge equivalence class of) an instanton on the multi-center Taub-NUT space, which is a solution of a four-dimensional PDE.  It was also proved in \cite{Cherkis:2010bn} that the Up transform produces an anti-self-dual connection. 

In \cite{First} we established core analytic results for 
instantons on the multi-Taub-NUT spaces, including the asymptotic form 
of the connection and the index of the Dirac operator. 
The goal of this  paper is to establish core analytic results 
for the moduli spaces of bow data and to formulate the Up 
transform and to establish the asymptotic form of any connection 
obtained through it.  

This sets the stage for a forthcoming 
paper \cite{Third}, where we formulate the Down transform and 
prove that the bow construction is in fact complete, as the 
two transforms are inverse to each other, and  the transforms 
are hyperk\"ahler isometries of respective moduli spaces.

The bows we consider in this paper can be specified by $k$ points 
on a circle:  $\{p_\sigma\}_{\sigma =1}^k\subset S^1$. 
A bow representation requires (among other data) 
an additional finite collection of points 
$\Lambda=\{\lambda_j\}_{j=1}^{|\Lambda|}\subset S^1\setminus \{p_\sigma\}_{\sigma =1}^k.$  
Let $\mathcal{E}$ denote the Hermitian bundle  with connection $A$, produced via the Up transform. In Theorem \ref{Thm:Index}, we compute its rank: 
\begin{theo}\label{rankth} 
	The rank of $\mathcal{E}$ is $|\Lambda|$. 
\end{theo}
We establish sharp asymptotics for the connection in Theorem \ref{Thm:Asymp}. The multi-Taub-NUT spaces are asymptotically  circle bundles over $\IR^3$. Let $\vec{t}$ denote coordinates on the $\IR^3$ base and $\tau$ a local coordinate for the circle fiber. 
\begin{theo}\label{MainTh} The curvature of $A$
 is $L^2$. Moreover, 
there is a local frame for $\mathcal{E}$ in which 
\begin{equation}\label{goal-1}
A(\frac{\p}{\p t^j})=  O(|\vec{t}\,|^{-2}),
\end{equation}
and 
\begin{equation}\label{goal}
A(\frac{\p}{\p\tau} )= -\ii\, {\rm diag}\frac{{\lambda_{i}+\frac{\hat{m}_{i}}{2|\vec{t}\,|}}}{V}+O(|\vec{t}\,|^{-2}),
\end{equation}
where $\hat{m}_{i}\in \IZ$ is given in terms of Bow data in Theorem \ref{Thm:Asymp}.
\end{theo}

In the third paper \cite{Third} in this series we formulate the Down transform and prove the completeness of the bow construction:
\begin{utheo}
The Up transform is bijective. 
\end{utheo}
\noindent and
\begin{utheo}
The Up transform is an isometry from the hyperk\"ahler moduli space of a bow representation  to the moduli space of the corresponding instanton on the multi-Taub-NUT space.
\end{utheo}
To prove these last two theorems we formulate the {\em Down transform} that maps an instanton on a multi-Taub-NUT space to a solution of a bow representation.  In this map, the multi-Taub-NUT space determines the {\em bow}, the topological type of the instanton itself determines the bow {\em representation}, and the instanton determines a {\em solution} of this bow representation (up to gauge transformations). 

\subsection*{Outline}
After reviewing definitions of the bow and its representation in Sec.~\ref{Sec:Bow}, we introduce its associated data and moduli space in Sec.~\ref{Sec:Data}. This moduli space is the hyperk\"ahler reduction of the bow representation  data by the gauge group of this bow representation.  
One ingredient of the bow data is the triplet of Nahm matrices with poles at some of the $\lambda$-points.  
In Sec.~\ref{ggmmnm} we establish the subleading order behavior of the Nahm matrices near  these  poles. 
In Sec.~\ref{Sec:SmallRep} we focus on the small representation $\mathfrak{s}$, with data $\mathrm{Dat}(\mathfrak{s}),$ gauge group $\mathcal{G}^0(\mathfrak{s})$, and hyperk\"ahler moment map $\mu_\mathfrak{s}$.  We show that its moduli space is the multi-Taub-NUT space:  TN$_k^\nu=\mu^{-1}_{\mathfrak{s}}(\ii\nu)/\mathcal{G}^0(\mathfrak{s}).$

The main objective of this paper is to show that the datum from any bow representation $\mathfrak{R},$ that satisfies the moment map conditions $\mu_\mathfrak{R}=\ii\nu$ and a nondegeneracy (WAF) condition given in  Sec.~\ref{Sec:Up}, can be used to construct a bundle with an anti-self-dual  connection with $L^2$-curvature, on the multi-Taub-NUT space TN$_k^\nu$. This is done in several steps.

First, in Sec.~\ref{BDDef}, we use both large and small  bow representations to construct  a family of ordinary differential  Dirac type operators $\D_{t,b}$, parametrized by the small bow representation level set, $\mu^{-1}_{\mathfrak{s}}(\ii\nu)\ni(t,b).$  These operators act on appropriate Sobolev spaces and  form a locally continuous family of Fredholm operators.

The moment map equations, coming from the gauge group action on bow data, imply that the $L^2$-kernel of $\D_{t,b}$ is trivial, since $\D_{t,b}^\dagger\D_{t,b}$ is strictly positive (as explained at the end of Sec.~\ref{BDDef}). Therefore, $\mathrm{Ind}\,\D^\dagger_{t,b}=\dim \mathrm{Ker}\, \D^\dagger_{t,b}$ and 
 $\mathrm{Ker}\, \D^\dagger_{t,b}$ forms a vector bundle over $\mu^{-1}_{\mathfrak{s}}(\ii\nu).$  Moreover, this bundle inherits the natural action of the gauge group $\mathcal{G}^0(\mathfrak{s})$ of the small representation $\mathfrak{s},$ and thus descends to a bundle $\mathcal{E}$ over the quotient  TN$_k^\nu=\mu^{-1}_{\mathfrak{s}}(\ii\nu)/\mathcal{G}^0(\mathfrak{s}).$

The induced connection $d_A$ on the bundle $\mathcal{E}$ over TN$_k^\nu$ is anti-self-dual, as proved in \cite[Sec.7]{Cherkis:2009jm} and in  \cite[Sec.3.2]{Cherkis:2010bn}. For completeness, we  include the proof in Appendix~\ref{ASDproof}.

In Theorem~\ref{Thm:Index} of Sec.~\ref{BowIndex}  we prove that the index of the operator family $\D^\dagger_{t,b},$ (which equals  the rank of its index bundle over TN$_k^\nu$)  equals the total number $|\Lambda|$ of $\lambda$-points in the large bow representation $\mathfrak{R}.$ This  implies Theorem~\ref{rankth}.

In Sec.~\ref{Sec:AsymInst} we analyze the asymptotic structure of the connection $d_A$  as $|(\vec{t},b)|\to \infty$. This depends on the decay rate of the Green's function $ (\D^\dagger_{t,b} \D_{t,b})^{-1}$. We establish quadratic decay for this Green's function in Lemma~\ref{Greendecay} in Sec.~\ref{Sec:GreenFn}. Using these results, we construct  in Sec.~\ref{Sec:AsymCon} a special frame for the index bundle and compute the asymptotic form of the connection matrix one-form in this frame, proving Theorem~\ref{Thm:Asymp}, which in turn implies our main  Theorem~\ref{MainTh}.

\SkipTocEntry\section*{Acknowledgements}
The work of SCh is supported by the Charles Simonyi Endowment at the Institute for Advanced Study. 
The work of ALH was supported, in part, by a SEED grant from the University of Dayton. The work of MS is supported by the Simons Foundation Grant 3553857.

\section{Bow Representations}\label{Sec:Bow}
The relevant bow for the $k$-centered Taub-NUT space ${\rm TN}_k$ is the affine $A_{k-1}$-bow, consisting of $k$ oriented intervals $I_\sigma=[p_{\sigma-1}+,p_\sigma-],$ $\sigma=1,\ldots,k$, with $k$ oriented edges $e_\sigma$ connecting the ends of consecutive intervals,  with $e_\sigma$  beginning at $p_\sigma-,$ the right endpoint of $I_\sigma,$ and ending at $p_\sigma+,$ the left endpoint of $I_{\sigma+1}.$  We understand the intervals to be cyclically ordered with $I_{\sigma+k}=I_\sigma.$ We denote the length of the $I_\sigma$ interval by $l_\sigma,$ and the total length by $l:=\sum_{\sigma=1}^k l_\sigma.$

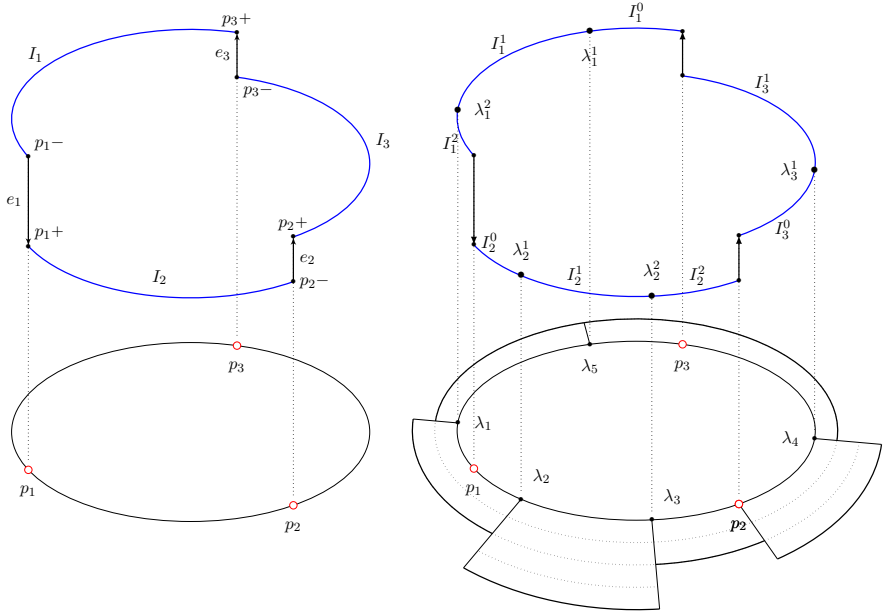
\begin{figure}[htb]
\centering
\resizebox{\textwidth}{!}{
%
%
%
%
%
%
%
%
%
%
\begin{tikzpicture}[>=latex']
    \draw (0,0) ellipse (4 and 2);
           \node(p1) at ($(-155:4 and 2)$) [label=below:$p_1$]{\tikz\draw[red,fill=white] (0,0) circle (.5ex);};
           \node(p2) at ($(-55:4 and 2)$) [label=below:$p_2$] {\tikz\draw[red,fill=white] (0,0) circle (.5ex);};
           \node(p3) at ($(75:4 and 2)$)  [label=below:$p_3$] {\tikz\draw[red,fill=white] (0,0) circle (.5ex);};
    \draw[blue,thick] (p1)+(0,5)  arc (-155:-55:4 and 2);
     \draw[blue,thick] (p2)+(0,6) arc (-55:75:4 and 2);
     \draw[blue,thick] (p3)+(0,7) arc (75:205:4 and 2);
\draw[dotted] (p1) -- ++(0,7);
\draw[dotted] (p2) -- ++(0,6);
\draw[dotted] (p3) -- ++(0,7);
 \node(Int1) at ($(0,7)+(150:4 and 2)$) [label=above:$I_1$] {};
 \node(Int2) at ($(0,5)+(-100:4 and 2)$) [label=above:$I_2$] {};
 \node(Int3) at ($(0,6)+(15:4 and 2)$) [label=right:$I_3$] {};
\draw[->,thick] (p1)+(0,7) node {\tiny$\bullet$} node[anchor=south west] {$p_1-$} -- ++(0,6) node[left] {$e_1$} -- ++(0,-1) node {\tiny$\bullet$} node[anchor=south west] {$p_1+$};
\draw[->,thick] (p2)+(0,5) node {\tiny$\bullet$} node[anchor= west] {$p_2-$} -- ++(0,5.5) node[right]{$e_2$} -- ++ (0,0.5) node {\tiny$\bullet$} node[anchor=south] {$p_2+$};
\draw[->,thick] (p3)+(0,6) node {\tiny$\bullet$} node[anchor=north west] {$p_3-$} -- ++(0,6.5) node[left]{$e_3$} -- ++(0,0.5) node {\tiny$\bullet$} node[anchor=south] {$p_3+$};
    \draw[opacity=0] ($(-130:6 and 4)$)  arc (-130:-85:6 and 4);
\end{tikzpicture}%
%
%
%
%
%
%
%
\begin{tikzpicture}
    \draw (0,0) ellipse (4 and 2);
           \node(p1) at ($(-155:4 and 2)$) [label=below:$p_1$] {\tikz\draw[red,fill=white] (0,0) circle (.5ex);}; 
           \node(p2) at ($(-55:4 and 2)$) [label=below:$p_2$] {\tikz\draw[red,fill=white] (0,0) circle (.5ex);};
           \node(p3) at ($(75:4 and 2)$)  [label=below:$p_3$] {\tikz\draw[red,fill=white] (0,0) circle (.5ex);};
     \draw[blue,thick] (p1)+(0,5)  arc (-155:-55:4 and 2);
     \draw[blue,thick] (p2)+(0,6) arc (-55:75:4 and 2);
     \draw[blue,thick] (p3)+(0,7) arc (75:205:4 and 2);
\draw[dotted] (p1) -- ++(0,7);
\draw[dotted] (p2) -- ++(0,6);
\draw[dotted] (p3) -- ++(0,7);
\draw[->,thick] (p1)+(0,7) node {\tiny$\bullet$}  -- ++(0,5) node {\tiny$\bullet$};
\draw[->,thick] (p2)+(0,5) node {\tiny$\bullet$}  -- ++(0,6) node {\tiny$\bullet$} ;
\draw[->,thick] (p3)+(0,6) node {\tiny$\bullet$}  -- ++(0,7) node {\tiny$\bullet$} ;
\coordinate (lam1) at ($(-185:4 and 2)$);
\coordinate (lam2) at ($(-130:4 and 2)$);
\coordinate (lam3) at ($(-85:4 and 2)$);
\coordinate (lam4) at ($(-5:4 and 2)$);
\coordinate (lam5) at ($(105:4 and 2)$);
           \node(l1) at (lam1) [label={[label distance=1pt]0:$\lambda_1$}] {\tiny$\bullet$};
           \node(l2) at (lam2) [label={[label distance=1]60:$\lambda_2$}] {\tiny$\bullet$};
           \node(l3) at (lam3)  [label={[label distance=1]60:$\lambda_3$}] {\tiny$\bullet$};
           \node(l4) at (lam4) [label=left:$\lambda_4$] {\tiny$\bullet$};
           \node(l5) at (lam5)  [label=below:$\lambda_5$] {\tiny$\bullet$};
\node at ($(90:4 and 2)+(0,7)$) [above] {$I_1^0$};
\node at ($(140:4 and 2)+(0,7)$) [above] {$I_1^1$};
\node at ($(-163:4 and 2)+(0,7)$) [left] {$I_1^2$};
\node at ($(-145:4 and 2)+(0,5)$) [above] {$I_2^0$};
\node at ($(-110:4 and 2)+(0,5)$) [above] {$I_2^1$};
\node at ($(-70:4 and 2)+(0,5)$) [above] {$I_2^2$};
\node at ($(-35:4 and 2)+(0,6)$) [below] {$I_3^0$};
\node at ($(45:4 and 2)+(0,6)$) [above] {$I_3^1$};
\draw[dotted] (l1) -- ++(0,7) node(lup1)[anchor=center,label=right:$\lambda_{1}^2$]{$\bullet$};
\draw[dotted] (l2) -- ++(0,5) node(lup1)[anchor=center,label=above:$\lambda_{2}^1$]{$\bullet$};
\draw[dotted] (l3) -- ++(0,5) node(lup1)[anchor=center,label=above:$\lambda_{2}^2$]{$\bullet$};
\draw[dotted] (l4) -- ++(0,6) node(lup1)[anchor=center,label=left:$\lambda_{3}^1$]{$\bullet$};
\draw[dotted] (l5) -- ++(0,7) node(lup1)[anchor=center,label=below:$\lambda_{1}^1$]{$\bullet$};
    \draw[thick] ($(-185:5 and 3)$) node(o1){} arc (-185:-130:5 and 3);
    \draw[thick] ($(-130:6 and 4)$) node(o2){}  arc (-130:-85:6 and 4) node(o3){};
    \draw[thick] ($(-85:5 and 3)$)  arc (-85:-55:5 and 3);
    \draw[thick] ($(-55:5.5 and 3.5)$) node(o35){}  arc (-55:-5:5.5 and 3.5) node(o4){};
    \draw[thick] ($(-5:4.5 and 2.5)$)  arc (-5:105:4.5 and 2.5);
    \draw[thick] ($(105:4.5 and 2.5)$) node(o5){} arc (105:175:4.5 and 2.5);
   \draw (lam1) -- ($(-185:5 and 3)$);
   \draw (lam2)--(o2.center);
   \draw (lam3)--(o3.center);
   \draw (p2.center)--(o35.center);
   \draw (lam4)--(o4.center);
   \draw (lam5)--(o5.center);
\node(p2) at ($(-55:4 and 2)$) [label=below:$p_2$] {\tikz\draw[red,fill=white] (0,0) circle (.5ex);};
       \draw[very thin,dotted] (0,0) ellipse (4.5 and 2.5);
       \draw[very thin,dotted] (o1.center) arc (-185:-5:5 and 3);
              \draw[very thin,dotted] (-85:5.5 and 3.5) arc (-85:-130:5.5 and 3.5);
\end{tikzpicture}%
}
\caption{On the left is a bow as encoded in its circle diagram. On the right is a bow representation with its ranks encoded in the $R$ function on the circle.}\label{fig:Rep}
\end{figure}

Another convenient way to specify an affine $A_{k-1}$ bow, illustrated in Figure~\ref{fig:Rep}, is to mark $k$ points $p_\sigma,$ forming the set 
$$P:=\{p_\sigma\}_{\sigma=1}^k,$$
 on a circle $S^1$ of perimeter $l.$  We denote the coordinate along this circle by $s;$ so, we identify $s\sim s+l.$  These marked points cut the circle into $k$ disjoint open intervals,  with $I_\sigma=[p_{\sigma-1}+,p_\sigma-]$ above being the closure of the resulting $\sigma^\mathrm{th}$ interval $(p_{\sigma-1},p_\sigma)$.  The projection $\pi: \sqcup_\sigma I_\sigma \rightarrow S^1$ is the identity on the complement of the $p$-points and maps both $p_{\sigma}+$ and $p_{\sigma}-$ to $p_\sigma.$

 \begin{define}           
A {\em representation} of a bow, denoted by $\mathfrak{R}$, is specified by 
\begin{enumerate}
\item
A finite set of points  $\Lambda=\{\lambda_i\}_{i=1}^n$ on the circle, disjoint from the  $\{p_\sigma\}_{\sigma=1}^k$, and
\item
A  locally constant rank function $R(s)$, with values in non-negative integers, defined on the circle complement of the sets $\Lambda$ and $P.$  
\end{enumerate}
\end{define}
Elements of $P$ and $\Lambda$ are called $p-$points and $\lambda-$points respectively. 
Since the set of $\lambda$-points is disjoint from the set of $p$-points, for each $\lambda_i$ its preimage $\pi^{-1}(\lambda_i)$ is a single point belonging to one of the intervals $I_{\sigma}$.  
 The order on $I_\sigma$ descends to an order on $\Lambda_\sigma:=\Lambda \cap \pi(I_\sigma)$; thus, we label the elements of $\Lambda_\sigma$ as $\lambda_{\sigma}^\alpha$, in increasing order, with $\alpha=1,2,\ldots,n_\sigma$, where $n_\sigma=\left|\Lambda_\sigma\right|$ is the total number of $\lambda$-points on the interval $I_\sigma.$  Each interval $I_\sigma$ is divided by the $\lambda$-points into subintervals $I_\sigma^\alpha$ with $\alpha=0,1,\ldots,n_\sigma.$ The value of $R(s)$ over  $I_\sigma^\alpha$ is  denoted by $R_\sigma^\alpha$.

We classify the $\lambda$-points by the behavior of $R(s)$ near $\lambda$.  Any $\lambda\in\Lambda$ has
\begin{itemize}
\item
{\em Continuous Rank} if $R(\lambda-)=R(\lambda+)$, and we denote the set of  such points  $\Lambda^0,$
\item
{\em Increasing Rank}  if $R(\lambda-)<R(\lambda+)$, and we denote the set of  such points  $\Lambda^+,$
\item
{\em  Decreasing Rank}  if $R(\lambda-)>R(\lambda+)$, and we denote the set of  such points $\Lambda^-.$ 
\end{itemize}
Let $\Lambda_\sigma^\ast = \Lambda_\sigma\cap \Lambda^\ast,$ for $\ast = 0,+,-$. 

\begin{figure}[ht]
\centering
\resizebox{0.8\textwidth}{!}{
%
%
%
%
\begin{tikzpicture}
    \draw (0,0) ellipse (4 and 2);
           \node(p1) at ($(-155:4 and 2)$) [label=below:$p_1$] {\tikz\draw[red,fill=white] (0,0) circle (.5ex);}; 
           \node(p2) at ($(-55:4 and 2)$) [label=below:$p_2$] {\tikz\draw[red,fill=white] (0,0) circle (.5ex);};
           \node(p3) at ($(75:4 and 2)$)  [label=below:$p_3$] {\tikz\draw[red,fill=white] (0,0) circle (.5ex);};
     \draw[blue,thick] (p1)+(0,5)  arc (-155:-55:4 and 2);
     \draw[blue,thick] (p2)+(0,6) arc (-55:75:4 and 2);
     \draw[blue,thick] (p3)+(0,7) arc (75:205:4 and 2);
\draw[dotted] (p1) -- ++(0,7);
\draw[dotted] (p2) -- ++(0,6);
\draw[dotted] (p3) -- ++(0,7);
\coordinate (lam1) at ($(-185:4 and 2)$);
\coordinate (lam2) at ($(-130:4 and 2)$);
\coordinate (lam3) at ($(-85:4 and 2)$);
\coordinate (lam4) at ($(-5:4 and 2)$);
\coordinate (lam5) at ($(105:4 and 2)$);
           \node(l1) at (lam1) [label={[label distance=1pt]0:$\lambda_1$}] {\tiny$\bullet$};
           \node(l2) at (lam2) [label={[label distance=1]60:$\lambda_2$}] {\tiny$\bullet$};
           \node(l3) at (lam3)  [label={[label distance=1]60:$\lambda_3$}] {\tiny$\bullet$};
           \node(l4) at (lam4) [label=left:$\lambda_4$] {\tiny$\bullet$};
           \node(l5) at (lam5)  [label=below:$\lambda_5$] {\tiny$\bullet$};
\draw[dotted] (l1) -- ++(0,7) node(lup1)[anchor=center,label=right:$\lambda_{1}^2$]{$\bullet$};
\draw[dotted] (l2) -- ++(0,5) node(lup1)[anchor=center,label=above:$\lambda_{2}^1$]{$\bullet$};
\draw[dotted] (l3) -- ++(0,5) node(lup1)[anchor=center,label=above:$\lambda_{2}^2$]{$\bullet$};
\draw[dotted] (l4) -- ++(0,6) node(lup1)[anchor=center,label=left:$\lambda_{3}^1$]{$\bullet$};
\draw[dotted] (l5) -- ++(0,7) node(lup1)[anchor=center,label=below:$\lambda_{1}^1$]{$\bullet$};
    \draw[thick] ($(-185:5 and 3)$) node(o1){} arc (-185:-130:5 and 3) node(o2in){};
    \draw[thick] ($(-130:6 and 4)$) node(o2){}  arc (-130:-85:6 and 4) node(o3){};
    \draw[thick] ($(-85:5 and 3)$) node(o3in){} arc (-85:-55:5 and 3) node(o35in){};
    \draw[thick] ($(-55:5.5 and 3.5)$) node(o35){}  arc (-55:-5:5.5 and 3.5) node(o4){};
    \draw[thick] ($(-5:4.5 and 2.5)$) node(o4in){} arc (-5:105:4.5 and 2.5) ;
    \draw[thick] ($(105:4.5 and 2.5)$) node(o5){} arc (105:175:4.5 and 2.5) node(o1in){};
   \draw (lam1) -- ($(-185:5 and 3)$);
   \draw (lam2)--(o2.center);
   \draw (lam3)--(o3.center);
   \draw (p2.center)--(o35.center);
   \draw (lam4)--(o4.center);
   \draw (lam5)--(o5.center);
\node(p2) at ($(-55:4 and 2)$) [label=below:$p_2$] {\tikz\draw[red,fill=white] (0,0) circle (.5ex);};
%
%
%
%
%
\fill [draw=none, fill=white, fill opacity=0.7]
   (-5.5,-4.5) -- (-5.5,3) -- (5.5,3) -- (5.5,-4.5) -- cycle;
\draw[dotted] (p1) -- ++(0,7);
\draw[dotted] (p2) -- ++(0,6);
\draw[dotted] (p3) -- ++(0,7);
\draw[->,thick] (p1)+(0,7) node(p1m) {\tiny$\bullet$}  to[out=-105,in=105] ++(0,5) node(p1p) {\tiny$\bullet$};
\draw[->,thick] (p2)+(0,5) node(p2m) {\tiny$\bullet$}  to[out=65,in=-65] ++(0,6) node(p2p) {\tiny$\bullet$} ;
\draw[->,thick] (p3)+(0,6) node(p3m) {\tiny$\bullet$}  to[out=60,in=-60] ++(0,7) node(p3p) {\tiny$\bullet$} ;
\node at ($(p1p)+(-0.5,1.1)$) {\small$B_{21}$};
\node at ($(p2m)+(-0.45,0.6)$) {\small$B_{23}$};
\node at ($(p3m)+(-0.42,0.4)$) {\small$B_{31}$};
\draw[->,thick] (p1)+(0,5)   to[out=75,in=-75] ++(0,7);
\draw[->,thick] (p2)+(0,6)    to[out=-115,in=115] ++(0,5);
\draw[->,thick] (p3)+(0,7)   to[out=-120,in=120] ++(0,6);
\node at ($(p1p)+(0.5,1.1)$) {\small$B_{12}$};
\node at ($(p2m)+(0.33,0.17)$) {\small$B_{32}$};
\node at ($(p3m)+(0.45,0.75)$) {\small$B_{13}$};
\draw[dotted] (l1) -- ++(0,7)  node {\tiny$\bullet$};
\draw[dotted] (l2) -- ++(0,5)  node {\tiny$\bullet$};
\draw[dotted] (l3) -- ++(0,5)  node {\tiny$\bullet$};
\draw[dotted] (l4) -- ++(0,6)  node {\tiny$\bullet$};
\draw[dotted] (l5) -- ++(0,7) node(l5up) {\tiny$\bullet$};
\draw (p1m.center)--++($(-155:1)$);
\draw (p1p.center)--++($(-155:1)$)  node(q1p){};
\draw (p2m.center)--++($(-55:1)$) node(q2m){};
\draw (p2p.center)--++($(-55:1.5)$)  node(q2p){};
\draw (p3m.center)--++($(75:0.5)$) node(q3m){};
\draw (p3p.center)--++($(75:0.5)$) node(q3p){};
   \draw (o1in.center)+(0,7) -- ($(o1.center)+(0,7)$) node(o1up){};
   \draw (o2in.center)+(0,5)--($(o2.center)+(0,5)$) node(o2up){};
   \draw (o3in.center)+(0,5)--($(o3.center)+(0,5)$) node(o3up){};
  \draw (o4in.center)+(0,6)--($(o4.center)+(0,6)$) node(o4up){};
    \draw[thick] (o1up)  arc (-185:-155:5 and 3);
    \draw[thick] (q1p)  arc (-155:-130:5 and 3);
    \draw[thick] (o2up) node(o2){}  arc (-130:-85:6 and 4);
    \draw[thick] (q2m)  arc (-55:-85:5 and 3);
    \draw[thick] (q2p)  arc (-55:-5:5.5 and 3.5);
    \draw[thick] ($(-5:4.5 and 2.5)+(0,6)$) arc (-5:75:4.5 and 2.5);
    \draw[thick] (q3p)  arc (75:175:4.5 and 2.5);
       \draw[very thin,dotted] (0,0) ellipse (4.5 and 2.5);
       \draw[very thin,dotted] (o1.center) arc (-185:-5:5 and 3);
                    \draw[very thin,dotted] (-85:5.5 and 3.5) arc (-85:-130:5.5 and 3.5);
\draw[->,thick] (l5up.center)    to[out=75,in=-75] ++(0,2) node(w5){};
\draw[->,thick] (w5.center)   to[out=-105,in=105] (l5up.center);
\node at ($(l5up)+(-0.4,1.3)$) {\small$I_{\lambda_5}$};
\node at ($(l5up)+(0.45,1.3)$) {\small$J_{\lambda_5}$};
\node[label={[label distance=-5]{$W_{\lambda_5}$}}] at (w5)  {$\bullet$};
\node at ($(90:4 and 2)+(0,6.9)$) [above] {$I_1^0$};
\node at ($(140:4 and 2)+(0,6.97)$) [above] {$I_1^1$};
\node at ($(-163:4 and 2)+(0,7)$) [left] {$I_1^2$};
\node at ($(-145:4 and 2)+(-0.1,4.35)$) [above] {$I_2^0$};
\node at ($(-110:4 and 2)+(0,4.35)$) [above] {$I_2^1$};
\node at ($(-70:4 and 2)+(0,4.35)$) [above] {$I_2^2$};
\node at ($(-35:4 and 2)+(0,6)$) [below] {$I_3^0$};
\node at ($(45:4 and 2)+(0,5.9)$) [above] {$I_3^1$};
\end{tikzpicture}
}
\caption[An affine $A_{k-1}$ bow representation above its circle diagram. Note the appearance of the auxiliary space at the continuous rank $\lambda$-point $\lambda_1^1=\pi^{-1}(\lambda_5).$]{An affine $A_{k-1}$ bow representation above its circle diagram. Note the appearance of the auxiliary space at the continuous rank $\lambda$-point $\lambda_1^1=\pi^{-1}(\lambda_5).$
$Q_{\lambda_5}=\left(\protect
\begin{smallmatrix} 
J_{\lambda_5}^\dagger \\ I_{\lambda_5}
\end{smallmatrix} 
\right)$ and 
$B_\sigma=\left(\protect
\begin{smallmatrix}
B_{\sigma+1,\sigma}^\dagger\\ B_{\sigma,\sigma+1}
\end{smallmatrix}\right).$
}
\end{figure}
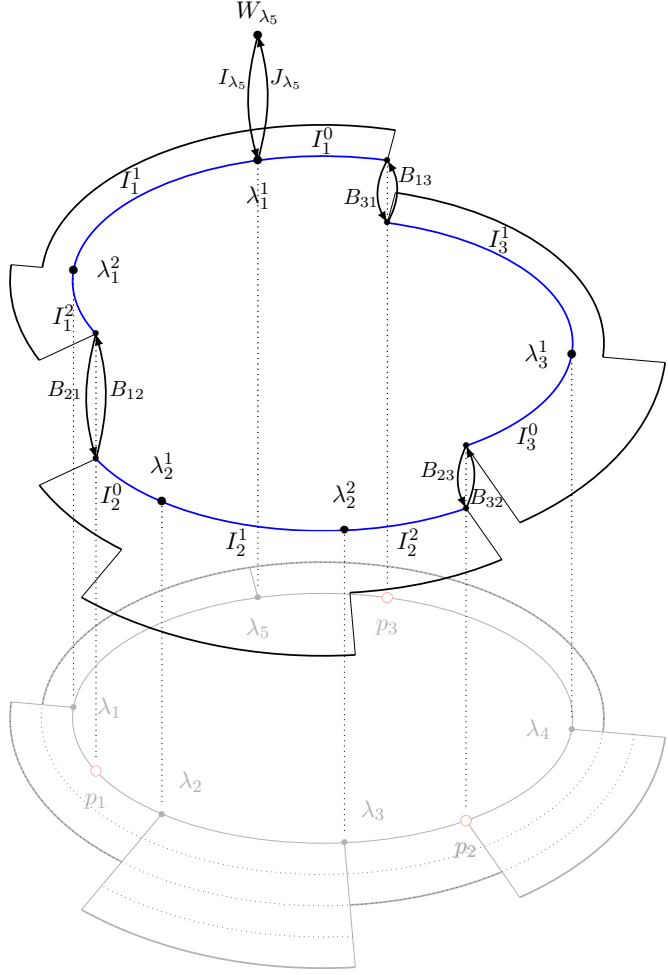

We define   a  Hermitian {vector bundle} $\FE$ over the bow to be a collection of Hermitian vector bundles $\FE_\sigma^\alpha\rightarrow I_\sigma^\alpha$, with $\mathrm{Rank}(\FE_\sigma^\alpha)=R_\sigma^\alpha$.  In the above notation, $R_\sigma^\alpha<R_\sigma^{\alpha+1}$ if $\lambda_\sigma^{\alpha+1}\in \Lambda^+$, $R_\sigma^\alpha>R_\sigma^{\alpha+1}$ if $\lambda_\sigma^{\alpha+1}\in \Lambda^-$, and $R_\sigma^\alpha=R_\sigma^{\alpha+1}$ if $\lambda_\sigma^{\alpha+1}\in \Lambda^0.$  

We  fix isometric injective maps $\kappa_\lambda$  of fibers at $\lambda-$points :
\begin{itemize}
\item At $\lambda_\sigma^{\alpha+1}\in \Lambda_\sigma^-$ we fix an injection $\ka_{\lambda_\sigma^{\alpha+1}}:\FE_{\sigma}^{\alpha+1}|_{\lambda_\sigma^{\alpha+1}}\to \FE_{\sigma}^{\alpha}|_{\lambda_\sigma^{\alpha+1}}.$ 
\item At $\lambda_\sigma^{\alpha+1}\in \Lambda_\sigma^+$  we fix an injection 
$\ka_{\lambda_\sigma^{\alpha+1}}:\FE_{\sigma}^{\alpha}|_{\lambda_\sigma^{\alpha+1}}\to \FE_{\sigma}^{\alpha+1}|_{\lambda_\sigma^{\alpha+1}}.$ 
\item At $\lambda_\sigma^{\alpha+1}\in \Lambda_\sigma^0$ we fix  an isomorphism 
$\ka_{\lambda_\sigma^{\alpha+1}}:\FE_{\sigma}^{\alpha}|_{\lambda_\sigma^{\alpha+1}}\to \FE_{\sigma}^{\alpha+1}|_{\lambda_\sigma^{\alpha+1}}.$ 
\end{itemize}
These provide an identification of the fiber of the lower rank bundle on one side of the $\lambda$-point with a subspace of the fiber of the higher rank bundle on  its other side. Henceforth, we will treat $\ka_{\lambda_\sigma^{\alpha+1}}$ as an inclusion and will frequently omit it  from our notation. 
\begin{define}
 We call $\mathrm{Im}\, \ka_\lambda$  
 the fiber of {\em continuous components}, denoted by $\FE_\lambda^{\rm cont}.$   We call the orthogonal complement of $\FE_\lambda^{\rm cont}$ in the larger fiber, the fiber of {\em terminating components}, denoted $\FE_\lambda^{\rm term}.$
\end{define}
 For $\lambda\in\Lambda^0$, $\FE_\lambda^{\rm term}=\{0\}$, and the fibers on the two sides are identified via $\kappa_\lambda$ and denoted $\FE_\lambda.$ 

The last ingredient in defining a bow representation is a set of Hermitian spaces $\{W_\lambda\}_{\lambda\in\Lambda^0}.$    
If we assign a multiplicity greater than one to some of the $\lambda_i$, then the dimension of $W_{\lambda_i}$ equals  the multiplicity of $\lambda_i.$  In this paper, we will only consider the case where each $\lambda$ has multiplicity one, and thus each $W_\lambda=\mathbb{C}.$

\section{Bow Data}\label{Sec:Data}

In this section, we associate to each bow representation $\mathfrak{R}$ an
 affine hyperk\"ahler space Dat$(\mathfrak{R})$, called the {\em data} of the bow representation $\mathfrak{R},$ as well as a hyperk\"ahler space $\mathcal{M}_\mathfrak{R}$ called the {\em moduli space} of $\mathfrak{R}.$

\subsection{Charge Conjugation}\label{Sec:ChargeConj}

 Let $\S$ be a one-dimensional right quaternionic vector space. Then $\S$ itself can be identified with $\IH$ and admits a left quaternionic action also. Choosing an embedding of $\IC$ in the quaternions, acting on the right (e.g. by sending $\ii$ to $I$),  allows us to view $\S$ as a Hermitian complex vector space equipped with an antilinear map $v\to vJ$, with $J^2 = -1$. Because the left action of $\IH$ commutes with the right action, the left action is complex linear with respect to this  (right acting) complex structure. This antilinear map and the Hermitian structure define a complex linear isomorphism   $\epsilon:\S\to \S^*$ by 
\begin{align}\epsilon(v)= \langle v J,\cdot\rangle,\end{align}
where our convention is that the inner product is antilinear in the first argument and linear in the second. 

The choice of $I$ and $J$ above amounts to a specific choice of basis in $\S,$ with the identification $\S\simeq\mathbb{C}^2$ given by $v=v_1+J v_2.$ Then the Hermitian structure is $\langle u,v\rangle:=\frac12(\bar{u}v- I\bar{u}v I),$  the $J$ action is $vJ=-\bar{v}_2+J\bar{v}_1$ amounts to an antilinear map $\left(\begin{smallmatrix} 
v_1   \\  
v_2
\end{smallmatrix} \right)\mapsto 
\left(\begin{smallmatrix} 
-\bar{v}_2\\
\bar{v}_1
\end{smallmatrix} \right)$, giving $\epsilon:\left(\begin{smallmatrix} 
v_1   \\  
v_2
\end{smallmatrix} \right)\mapsto (-v_2,v_1).$
 
Given a second complex vector space $B$, the isomorphism $\epsilon$ extends in a natural way to a complex linear isomorphism 
$$\epsilon :\S\otimes_{\IC}B\to \S^*\otimes_{\IC}B.$$
If $A$ is a third Hermitian complex vector space, we have 
\begin{align}\label{chargec}
\mathrm{Hom}\,(A,\S\mathop{\otimes}_{\IC}B)= \S\mathop{\otimes}_{\IC}B\mathop{\otimes}_{\IC}A^*\stackrel{\epsilon }{\rightarrow}
\S^*\mathop{\otimes}_{\IC}B\mathop{\otimes}_{\IC}A^*=\mathrm{Hom}\,(A\mathop{\otimes}_{\IC}\S,B).
\end{align}
Given $\beta\in \mathrm{Hom}\,(A,\S\otimes_{\IC}B)$, we let $\epsilon(\beta)\in \mathrm{Hom}\,(A\otimes_{\IC}\S,B)$ denote the image of $\beta$ under the maps in \eqref{chargec}. Then its Hermitian adjoint is $\epsilon(\beta)^*\in \mathrm{Hom}\,(B,\S\otimes_{\IC}A)$, 
and we denote $\epsilon(\beta)^*=:\beta^c$. We call the resulting antilinear map $\beta\to \beta^c$ {\em charge conjugation.} 

Since $\mathrm{Hom}\,(B,A\otimes_{\IC}\S)=A\otimes_{\IC}\S\otimes B^*= \mathrm{Hom}\,(A^*,\S\otimes_{\IC}B^*)$ one can view $\epsilon(\beta)^*$ as an element of either of these Hom spaces.  The two are related by transposition, and to distinguish the two we have the charge conjugate of $\beta$ in the former space
$$\beta^c=\epsilon(\beta)^*\in\mathrm{Hom}\,(B,A\otimes_{\IC}\S),$$ 
and its transpose (denoted by $^r$) in the latter space
$$^r\!\beta^c=^r\!\!(\epsilon(\beta)^*)\in \mathrm{Hom}\,(A^*,\S\otimes_{\IC}B^*).$$
Explicitly, for $\beta=\left(\begin{smallmatrix} 
\beta_1  \\  
\beta_2
\end{smallmatrix} \right)\in \mathrm{Hom}\,(A,\S\otimes_{\IC}B), \epsilon(\beta)=(-\beta_2,\beta_1),$ and
\begin{align*}
\beta^c&=\epsilon(\beta)^*=\left(\begin{smallmatrix} 
-\beta_2^*  \\  
\beta_1^*
\end{smallmatrix} \right)\in \mathrm{Hom}\,(B,\S\otimes_{\IC}A),\\
^r\!\beta^c&=\left(\begin{smallmatrix} 
-(\beta_2^*)^t  \\  
(\beta_1^*)^t
\end{smallmatrix} \right)\in \mathrm{Hom}\,(A^*,\S\otimes_{\IC}B^*).
\end{align*}

\subsection{The Hyperk\"ahler Space of Bow Data }

Let $\mathfrak{R}$ be a bow representation, then, as in \cite[Sec.2.3.2 and Sec.2.4]{Cherkis:2010bn}, ${\rm Dat}(\mathfrak{R}):=\mathcal{F}\times\mathcal{B}\times\mathcal{N}$  is a direct product of three spaces:
\begin{itemize}
\item[-]
The {\em Fundamental} linear space $\mathcal{F}:=\mathop{\oplus}\limits_{\lambda\in\Lambda^0} \mathrm{Hom}\,(W_\lambda, \S\otimes \FE_\lambda),$ with elements  denoted by $(Q_\lambda)_{\lambda\in\Lambda^0},$
\item[-]
The {\em Bifundamental} linear space $\mathcal{B}:=\mathop{\oplus}\limits_{\sigma=1}^k \mathrm{Hom}\,(\FE_{p_{\sigma}+}, \S\otimes \FE_{p_{\sigma}-}),$ with elements denoted by $(B_{\sigma})_{\sigma=1}^k$. Alternatively, we may use charge conjugation to identify  $\mathcal{B}$ antilinearly with another space $\mathop{\oplus}\limits_{\sigma=1}^k \mathrm{Hom}\,(\FE_{p_{\sigma}-}, \S\otimes \FE_{p_{\sigma}+}),$ with elements $(B_\sigma^c)_{\sigma=1}^k.$ 
 
\item[-]
Fix a smooth Hermitian connection $\nabla^0$ on $ \FE$. The {\em Nahm} affine space, $\mathcal{N}$, is the space of quadruplets $(\nabla, T_1, T_2, T_3)$ of a Hermitian connection $\nabla$ and three Hermitian sections $\{T_j\}_{j=1}^3$  of $\mathrm{End}(\FE),$ which satisfy the following conditions.
\begin{itemize}
\item[i)] $\nabla - \nabla^0\in L^2(\cup_\sigma I_\sigma \setminus\Lambda,\FE)$. 
\item[ii)] For each $\sigma,$ ${T_j}\in L^2_{loc}(I_\sigma \setminus(\Lambda_\sigma^+\cup\Lambda_\sigma^-))$.
\item[iii)]
If $\lambda\in\Lambda^+$, i.e. the rank change $m_\lambda:=R(\lambda+)- R(\lambda-)$ across $\lambda$ is positive,  then for some $\epsilon >0,$ in a $\nabla$ covariant constant frame,
\begin{align}\label{Pole1}
T_j(\lambda+s) + 
\begin{pmatrix}
\frac{\ii{\bm{\rho}}_j}{2s}  & 0\\
0 & 0
\end{pmatrix}\in L^2([0,\epsilon)). 
\end{align}
\item[iv)]
 If $\lambda\in\Lambda^-$, i.e. the rank change $m_\lambda:=R(\lambda+)- R(\lambda-)$ across $\lambda$ is negative,
then for some $\epsilon >0,$ in a $\nabla$ covariant constant frame,
\begin{align}\label{Pole2}
{T}_j(\lambda+s)  +   
\begin{pmatrix}
\frac{\ii{\bm{\rho}}_j}{2s}  & 0\\
0 &0
\end{pmatrix}\in L^2((-\epsilon,0]).
\end{align}
\end{itemize}
\end{itemize}
Here,  ${\bm{\rho}}_1, {\bm{\rho}}_2,$ and ${\bm{\rho}}_3$ form an $|m_\lambda|$-dimensional irreducible representation of $su(2),$ i.e. $[{\bm{\rho}}_i, {\bm{\rho}}_j]=2{\bm{\rho}}_k$ for any cyclic permutation $(i,j,k)$ of $(1,2,3).$   The block matrix decomposition  corresponds to the (covariant constant extension) of the decomposition of $\FE_{\lambda}$  into terminating and continuing summands denoted $\left(\begin{array}{c}\FE^{\text{term}}\\\FE^{\text{cont}}\end{array}\right).$   

Each of these spaces admits a Hermitian structure that is compatible with the induced left action of the quaternions, making them (flat) hyperk\"ahler Hilbert manifolds as follows.  
The fundamental and bifundamental spaces inherit the left action of quaternions from $\S.$ For the Nahm space tangent vector, set $\delta T_0:=-\ii\delta\nabla, $ $\delta\tilde \T:=\e_1\otimes\delta T_1+\e_2\otimes\delta T_2+\e_3\otimes\delta T_3,$ and $\delta\T:= \delta T_0 + \delta\tilde \T,$ where  $\{\e_1,\e_2,\e_3\}$ is a standard basis of unit imaginary quaternions. Note, that since the $T_j$ are Hermitian, the Hermitian conjugate is $\delta\T^\dagger=\delta T_0-\delta\tilde\T.$   Then the tangent space to the Nahm space is $\oplus_\sigma L^2(I_\sigma,\mathbb{H}\otimes\mathrm{End}\,\FE).$  Its quaternionic structure is given by the quaternionic unit action on $ \delta\T$  by multiplication on the left.

The hyperk\"ahler metric on Dat$\,\mathfrak{R}$ is 
\begin{align}
g=\sum_{\lambda\in\Lambda^0} \tr_{\FE_\lambda} \delta Q_\lambda^\dagger \delta Q_\lambda
     +\sum_\sigma \tr_{\FE_{p_\sigma+}} \delta B_\sigma^\dagger \delta B_\sigma
     +\frac12\int \tr_{\S\otimes\FE}  \delta\T^\dagger \delta\T ds.
\end{align} 
The three symplectic forms $\omega_j:=g(\e_j\cdot,\cdot)$ can be assembled into an imaginary quaternion-valued 2-form
$\omega:=\e_1\otimes\omega_1+\e_2\otimes\omega_2+\e_3\otimes\omega_3.$ Then 
\begin{align}
\omega=\sum_{\lambda\in\Lambda^0}  \tr_{W_\lambda}\delta Q\wedge\delta Q^\dagger
              &+ \sum_\sigma\tr_{\FE_{p_\sigma-}}\delta B_\sigma\wedge\delta B_\sigma^\dagger+\frac12\int \tr_{\FE} \delta\T \wedge \delta\T^\dagger ds.
\end{align}

\subsection{Gauge Group, Moment Maps, and Nahm Data}\label{ggmmnm}

We define the Sobolev space $H^1$ for various vector bundles over the bow to consist of sections $s$ which satisfy 
$\|s\|^2_{H^1}:= \|s\|^2_{L^2}+\|\nabla^0 s\|^2_{L^2}<\infty.$ As our bundles are only smooth over $\bigsqcup_\sigma I_\sigma \setminus (\Lambda\setminus \Lambda^0)$, we require the sections $s$ lie in the maximal domain of $\nabla^0 $ on this incomplete manifold.  We  take the group 
$\mathcal{G}(\mathfrak{R})$ of allowed gauge transformations of $\mathfrak{R}$ to consist of unitary $H^1$-sections of $\text{End}(\FE)$ on $\bigsqcup_\sigma I_\sigma \setminus (\Lambda\setminus \Lambda^0)$ which satisfy the transmission condition at each $\lambda\in \Lambda^+$ : 
\begin{align}\label{transmission}\lim_{s\to \lambda^+} g(s) = \left(\begin{array}{cc}I_{\text{term}}&0\\0&\lim\limits_{s\to \lambda^-}g(s)\end{array}\right),
\end{align}
and the analogous condition at $\lambda \in \Lambda^-$.   Here $I_{\text{term}}$ denotes the identity on the terminating summand. In particular, 
 these gauge transformations preserve the subspace of continuing components.

By Sobolev's multiplication theorem in dimension one (see  \cite[appendix]{DK1}), the multiplication map $H^k\times H^j\rightarrow H^j$ is bounded for $k>\frac{1}{2}$ and $k\geq j;$ thus, $\mathcal{G}(\mathfrak{R})$ is indeed a group.

The gauge group $\mathcal{G}(\mathfrak{R})$ acts naturally on sections of $\FE$ and  admits the following  triholomorphic action on each of the spaces 
$\mathcal{F},\,\,\mathcal{ B}, $ and $\mathcal{N}$.   

An element  $g\in\mathcal{G}(\mathfrak{R})$  
\begin{itemize}
\item acts on   $Q_\lambda\in \mathcal{F}$ as
\begin{align*}
g\cdot Q_\lambda=(1_\S\otimes g^{-1}(\lambda))Q_\lambda.
\end{align*}

\item acts on  $B_\sigma\in\mathcal{B}$ as 
\begin{align*}
g\cdot B_{\sigma}= (1_\S\otimes g^{-1}(p_{\sigma}-)) B_\sigma  g(p_{\sigma}+).
\end{align*}

\item acts on the Nahm data $(\nabla,T_1,T_2,T_3)\in\mathcal N$ as 
\begin{align*}
g\cdot \nabla&=g^{-1}\nabla\, g, & 
g\cdot{T}_j&=g^{-1}{T}_j\, g,\,\,\, j=1,2,3.
\end{align*}
\end{itemize}
Let $\tilde{\T}:=\e_1\otimes T_1+\e_2\otimes T_2+\e_3\otimes T_3.$ 
We denote the moment map of this action by $\mu_\mathfrak{R}: {\rm Dat}(\mathfrak{R})\rightarrow \mathpzc{g}_\mathfrak{R}^*\otimes\mathbb{R}^3=\mathpzc{g}_\mathfrak{R}^*\otimes\mathrm{Im}\,\mathbb{H},$ where $\mathpzc{g}_\mathfrak{R}^*$ is the dual of the Lie algebra $\mathpzc{g}_\mathfrak{R}$ of $\mathcal{G}(\mathfrak{R}).$ As $\mathpzc{g}_\mathfrak{R}$ is a closed subspace of $H^1$, $\mathpzc{g}_\mathfrak{R}^*\subset H^{-1}$ and $\mu_\mathfrak{R}$ is distribution valued.   
Any  element of the gauge algebra $\mathpzc{g}_\mathfrak{R}$ is a skew-hermitian function $F$ on the bow.  Let $F_\lambda, F_\sigma^+$ and $F_\sigma^-$ denote its respective values at $s=\lambda, p_\sigma+$ and $p_\sigma-.$ It generates a vector field 
$$X_F=
\begin{pmatrix}
\delta Q_\lambda\\ \delta B_\sigma\\ \delta T_0\\ \delta\tilde{\T}
\end{pmatrix}
=
\begin{pmatrix}
- F_\lambda Q_\lambda\\  
B_\sigma F_\sigma^+  - F_\sigma^- B_\sigma\\
-\ii [\nabla, F]\\
[\tilde{\T},F]
\end{pmatrix}.
$$
The defining equation for the moment map $\bm{\mu}_\mathfrak{R}$ is $d\bm{\mu}_\mathfrak{R}[F]=\omega(X_F,\cdot).$
The explicit  expression for the moment map at a point $(Q,B,\nabla,\tilde{\T})\in {\rm Dat}(\mathfrak{R})$  evaluated on  $F\in \mathpzc{g}_\mathfrak{R}$ is
\begin{multline}
\bm{\mu}_\mathfrak{R}(Q,B,\nabla,\tilde{\T})[F]=\int \tr_{\FE}\, \mu_\mathfrak{R}(Q,B,\nabla,\tilde{\T})(s) F(s) ds\\=
-\sum_{\lambda\in\Lambda^0} \tr_{\FE_\lambda} F(\lambda) Q_\lambda Q_\lambda^\dagger
+ \sum_\sigma \tr_{\FE_{p_\sigma-}}   (B_\sigma F_\sigma^+ B_\sigma^\dagger 
                                                          -  F_\sigma^- B_\sigma B_\sigma^\dagger)\\
-\int \tr_{\FE_s} (F [\ii\nabla,\tilde{\T}]+\frac12(\tilde{\T} F\tilde{\T}-F\tilde{\T}\tilde{\T}))ds
+\sum_\sigma\tr_{\FE} \ii\tilde{\T} F|^{p_\sigma-}_{p_\sigma+}.
\end{multline}
For any $\mathbb{C}$-linear space $V$ and for $M\in\mathrm{End}(\S\otimes V)$ define the imaginary part of $M$ by
\begin{align}\label{imdef}\mathrm{Im}\,M:=\frac12 \e_j\otimes\tr_{\S} M\e_j^\dagger=M-\frac12 1_\S\otimes\tr_\S M,
\end{align} so that for $M=1\otimes M_0+\e_j\otimes M_j$ one  has  $\mathrm{Im}\,M=\e_j\otimes M_j$.  
Using 
\begin{align}\tr_{\FE_{p_\sigma-}} \mathrm{Im}\, B_\sigma F B_\sigma^\dagger=-\tr_{\FE_{p_\sigma+}} \mathrm{Im}\,F B_\sigma^c(B_\sigma^c)^\dagger,
\end{align} 
\begin{multline}\label{Moment}
\mu_\mathfrak{R}(Q,B,\nabla,\tilde{\T})
=\mathrm{Im}\Bigg(
                              -\sum_{\lambda\in\Lambda^0}\delta(s-\lambda)Q_\lambda Q_\lambda^\dagger-[\ii\nabla,\tilde{\T}]+\tilde{\T}\tilde{\T}\\
                       +\sum_\sigma\Big(
       \delta(s-p_\sigma-)(\ii\tilde{\T}(p_\sigma-)-B_\sigma B_\sigma^\dagger  )\\
      -\delta(s-p_\sigma+)(\ii\tilde{\T}(p_{\sigma}+)+B_\sigma^c(B_\sigma^c)^\dagger)
                                                \Big)
                             \Bigg).
\end{multline}
For any $\ii\nu\in \mathpzc{g}_\mathfrak{R}^*\otimes\mathbb{R}^3=\mathpzc{g}_\mathfrak{R}^*\otimes\mathrm{Im}\,\mathbb{H}$ that is invariant under the coadjoint action of $\mathcal{G}(\mathfrak{R})$,  
the level set of $\ii\nu$,  $\bm{\mu}_\mathfrak{R}^{-1}(\ii\nu),$ is preserved  by $\mathcal{G}(\mathfrak{R})$. Hence, the quotient  $\bm{\mu}_\mathfrak{R}^{-1}(\ii\nu)/\mathcal{G}(\mathfrak{R})$  is well defined; we call it the {\em moduli space}  {\em of the bow representation} $\mathfrak{R}$ at level $\ii\nu $ and denote it $\mathcal{M}^\nu_\mathfrak{R}$. 
If the gauge group action has no fixed points in $\bm{\mu}_\mathfrak{R}^{-1}(\ii\nu)$, $\mathcal{M}_\mathfrak{R}^\nu$ is a smooth hyperk\"ahler manifold.  It is the hyperk\"ahler quotient of  $\mathrm{Dat}(\mathfrak{R})$ by  $\mathcal{G}(\mathfrak{R})$, as defined in \cite{Hitchin:1986ea}.

For reasons discussed in \cite[Sec.~2.5]{Cherkis:2010bn}, we consider only the moment map level $\ii\nu$ with $\nu$ in the form
\begin{align}\label{Level}
\nu=\sum_\sigma \nu_\sigma\otimes 1_{\FE} \left(\delta_{p_{\sigma -}}-\delta_{p_{\sigma+}}\right),
\end{align}
with each $\nu_\sigma$ a pure imaginary quaternion.  
Note: for brevity, we use $\delta_\lambda$ and $\delta_{p_\sigma\pm}$ to denote the delta functions $\delta_\lambda:=\delta(s-\lambda)$ and $\delta_{p_{\sigma\pm}}:=\delta(s-p_{\sigma\pm}).$  We also assume $\nu_\sigma\not =\nu_{\sigma'}$ for $\sigma\not = \sigma'$, which will ensure that the resulting instanton base space is nonsingular.

Imposing the moment map condition $\mu_\mathfrak{R}=\ii\nu$ with \eqref{Moment} and \eqref{Level} implies\newline 
1) the Nahm equations on the $T_j$:
\begin{align}\label{nahmeq0}
[\ii\nabla_{\frac{d}{ds}},T_1]=[T_2,T_3], \text{ and cyclic permutations,}
\end{align} 
2) the boundary conditions at $\lambda\in\Lambda_0:$
\begin{align}\label{Tdisc}
\tilde{\T}(\lambda+)-\tilde{\T}(\lambda-)=\mathrm{Im}\,\ii Q_\lambda Q_\lambda^\dagger,\text{  and}
\end{align}
3) the boundary conditions at the ends of the intervals
\begin{align}\label{pEnd}
\tilde{\T}(p_\sigma-)-\nu_\sigma&=-\mathrm{Im}\, \ii B_\sigma B_\sigma^\dagger,&
\tilde{\T}(p_\sigma+)-\nu_\sigma&=\mathrm{Im}\, \ii B_\sigma^c \left(B_\sigma^c\right)^\dagger.&
\end{align}
Set $\rho_j:= \begin{pmatrix}
\bm{\rho_j} & 0\\
0 & 0
\end{pmatrix}.$ 
In the remainder of this section we prove that in the vicinity of each $\lambda$-point, the function 
$\hat{b}_j(s):=T_j(\lambda+s) + \frac{{\ii}\rho_j}{2s}$ is bounded. 
\begin{lemma}\label{firstbound}
 Let 
$\{T_j\}_{j=1}^3$ satisfy the Nahm equations and conditions \eqref{Pole1} and \eqref{Pole2}. Then \newline
$T_j(\lambda+s) +  \frac{{\ii}\rho_j}{2s}  \in s^{-1/2}L^\infty([0,\epsilon])$, for $\lambda\in \Lambda^+$ and \newline
$T_j(\lambda+s)  +
\frac{{\ii}\rho_j}{2s} \in |s|^{-1/2}L^\infty([-\epsilon,0]),$ for $\lambda\in \Lambda^-.$
\end{lemma}
\begin{proof}
Suppose $\lambda\in \Lambda^+$ and $\{T_j\}_{j=1}^3$ satisfies the Nahm equations and \eqref{Pole1}.  Write, as above,  
$T_j = \hat{b}_j-\frac{{\ii}\rho_j}{2s}.$ Then 
\begin{align}\label{l1est}
\frac{d}{ds}(s^{\frac{2}{3}}\ii\hat{b}_1) = \frac{2}{3}s^{-\frac{1}{3}}\ii\hat{b}_1+s^{\frac{2}{3}}[\hat{b}_2,\hat{b}_3]-s^{-\frac{1}{3}}[\hat{b}_2,\frac{{\ii}\rho_3}{2}]-s^{-\frac{1}{3}}[\frac{{\ii}\rho_2}{2},\hat{b}_3].
\end{align}
The right hand side of \eqref{l1est} is $L^1$. Therefore $s^{\frac{2}{3}}\hat{b}_1$ is bounded and has a well defined limit $\beta$ at $s = 0$.  The assumption $\hat{b}_1\in L^2$ implies $\beta = 0.$ Hence integrating \eqref{l1est} from $0$ to $s$ and applying Cauchy-Schwartz, we deduce $s^{\frac{1}{2}}\hat{b}_1$ is bounded. By symmetry $s^{\frac{1}{2}}\hat{b}_j$ is bounded for all $j$. 
\end{proof}
For  convenience, we include a trivial ODE lemma which we will use repeatedly to bound functions.
\begin{lemma}\label{minor}
Let  
$A\in \mathrm{End}(\IR^N)$ and $x:(0,a]\to \IC^N$ such that  
$|s^{\frac{1}{2}}x|$ is bounded and 
\begin{align}\label{Eq:minoreq}
\frac{d}{ds} x +\frac{1}{s}Ax +E(s)=0,
\end{align}
with $|sE(s)|$ continuous and bounded. If $A$ is diagonalizable and all eigenvalues of $A$ lie in $(-\infty,0)\cup(\frac{1}{2},\infty)$, then  $|x(s)|$ is bounded. 
\end{lemma}
\begin{proof}Diagonalizing $A$, it suffices to consider the scalar equation 
\begin{align}\frac{d}{ds}y+\frac{\lambda}{s}y +E_\lambda(s)=0,
\end{align}
with $|sE_\lambda(s)|$ continuous and bounded. For $\lambda<0$, we have for $s>0$, 
\begin{align}y(s)=s^{-\lambda}a^{\lambda}y(a)+s^{-\lambda}\int_s^at^\lambda E_\lambda(t)dt.
\end{align}
The integral, and therefore $y(s)$ is uniformly bounded. For $\lambda>\frac{1}{2}$, we have 
\begin{align}y(s)= -s^{-\lambda}\int_0^st^\lambda E_\lambda(t)dt,
\end{align}
which is also uniformly bounded. 
\end{proof}

\begin{lemma}\label{secondbound} Let 
$\{T_j\}_{j=1}^3$ satisfy the Nahm equations and conditions \eqref{Pole1} and \eqref{Pole2}. Then \newline
$T_j(\lambda+s) +\frac{   {\ii}\rho_j}{2s}  \in L^\infty([0,\epsilon])$, for $\lambda\in \Lambda^+$ and\newline 
$T_j(\lambda+s)  + \frac{  {\ii} \rho_j}{2s} \in L^\infty([-\epsilon,0]),$ for $\lambda\in \Lambda^-.$
\end{lemma}
\begin{proof}

Since the $\{T_j\}_{j=1}^3$ satisfy the Nahm equations,
 the Jacobi identity implies 
\begin{align}\label{jacobi}\sum_j[T_j,\dot T_j] = 0.
\end{align}

We further decompose each $T_j$  as 
\begin{align}\ii T_j(\lambda+s) = \left(\frac{1}{2s}+f_j(s)\right)\rho_j+\gamma_j(s),
\end{align}
where $\gamma_j\perp\rho_j$ and, by Lemma \ref{firstbound}, $|\gamma_j|+|f_j|\leq Cs^{-1/2}$.  
In this notation, equation \eqref{jacobi} becomes 
\begin{align}\label{jacobi3}0 &= \sum_j[(\frac{1}{2s}+f_j)\rho_j+\gamma_j,(-\frac{1}{2s^2}+\frac{d}{ds}f_j)\rho_j+\frac{d}{ds} \gamma_j]\nonumber\\
&= \frac{1}{s^2}\frac{d}{ds}\sum_j [\frac{\rho_j}{2},   s  \gamma_j]+\sum_j[\gamma_j+f_j\rho_j,\frac{d}{ds} \gamma_j+\frac{d}{ds}f_j\rho_j],
\end{align}
and the Nahm equations 
\eqref{nahmeq0} become the pair of equations (and their cyclic permutations):
\begin{align}\label{diag} 
\left(\frac{d}{ds} f_1+\frac{f_2+f_3}{s}+2f_2f_3\right)|\rho_1|^2
+\langle [ \gamma_2, \gamma_3],\rho_1\rangle = 0,
\end{align}
and
\begin{multline}\label{fineb1}
\frac{d}{ds} \gamma_1+(\frac{1}{2s}+f_2)[\rho_2,\gamma_3]+ (\frac{1}{2s}+f_3)[ \gamma_2,\rho_3 ]\\
+ [ \gamma_2, \gamma_3]- \langle [ \gamma_2, \gamma_3],\rho_1\rangle\frac{\rho_1}{|\rho_1|^2} = 0.
\end{multline}
 
First, we prove that each $f_j$ is bounded. Equation \eqref{diag} and the $O(s^{-1/2})$ bound on the $f_j$ and $\gamma_j$ from Lemma~\ref{firstbound} imply $|\frac{d}{ds}f_j|=O(s^{-3/2})$. Similarly, \eqref{fineb1} implies $|\frac{d}{ds}\gamma_j|=O(s^{-3/2})$. 
Hence \eqref{jacobi3} yields 
\begin{multline}\label{jacobi4}   
\sum_j[\rho_j, \gamma_j(s)]=\\
 \frac{2}{s}\int_0^su^2\sum_j \left[\frac{d}{du} \gamma_j(u)+\frac{d}{du} f_j(u)\rho_j, \gamma_j(u)+f_j(u)\rho_j\right]du = O(1).
\end{multline}

From \eqref{diag} we see that 
$
\frac{d}{ds}[s^2(f_1+f_2+f_3)] = O(s).
$
Hence $f_1+f_2+f_3$ is bounded. 
Combining this estimate with \eqref{diag} yields 
\begin{align}\label{diag2} \frac{d}{ds}(\frac{f_j}{s}) = O(s^{-2}).
\end{align}
Hence each $f_j$ is bounded.

Next, we prove that the $\gamma_j$ are bounded as well. Set  $\rho_+ = \frac{\rho_1+\ii \rho_2}{2}$ and $\rho_-=\rho_+^\dagger=-\frac{\rho_1-\ii\rho_2}{2}$.  Then, together with $\ii\rho_3$, they  form the standard basis of $sl(2)$:
\begin{align}
[\ii\rho_3,\rho_+]&=2\rho_+,& [\ii\rho_3,\rho_-]&=-2\rho_-,&[\rho_+,\rho_-]=\ii\rho_3.
\end{align}
Recalling some standard representation theory, the $d$-dimensional irreducible representation of $sl(2)$ is given, e.g. by homogeneous  polynomials of degree $d-1$ of two variables: $V^{d}=\mathbb{C}^{d-1}[x,y]= \mathrm{Sym}^{d-1}\mathbb{C}^2.$ The defining two-dimensional representation has $\ii\rho_3 x=x, \ii\rho_3y=-y, \rho_+x=0,$ and $\rho_+y=x.$  A vector in $V^d$ given by a monomial $v_\beta=x^a y^b$ (has $a+b=d-1$ and) has weight $\beta=a-b$ as 
$\ii\rho_3 x^ay^b=(a-b)x^ay^b$ (via the product rule). And the raising and lowering operators act as follows: 
$\rho_+ v_\beta=\rho_+ x^ay^b=b x^{a+1}y^{b-1}=b v_{\beta+2},$ where  $a+b=d-1$ and $a-b=\beta$.  Thus
\begin{align}
\rho_+ v_\beta=\frac{d-\beta-1}{2}v_{\beta+2}.
\end{align}
Also,  
$\rho_- v_{\beta+2}=\rho_- x^a y^b=a x^{a-1}y^{b+1}=a v_{\beta},$ 
where $a+b=d-1$ and $a-b=\beta+2;$ so,
\begin{align}
\rho_- v_{\beta+2}=\frac{d+\beta+1}{2}v_\beta.
\end{align}

We recast Eqs.~\eqref{jacobi4} and \eqref{fineb1} as follows
\begin{align}\label{Eq:jacgamma}
[\rho_1,\gamma_1]+[\rho_2,\gamma_2]+[\rho_3,\gamma_3]=O(1),
\end{align}
\begin{align}
\label{Eq:gamma1}
\frac{d}{ds}\gamma_1+\frac{[\rho_2,\gamma_3]}{2s}+\frac{[\gamma_2,\rho_3]}{2s}&=O(s^{-1}),\\
\label{Eq:gamma2}
\frac{d}{ds}\gamma_2+\frac{[\rho_3,\gamma_1]}{2s}+\frac{[\gamma_3,\rho_1]}{2s}&=O(s^{-1}),\\
\label{Eq:gamma3}
\frac{d}{ds}\gamma_3+\frac{[\rho_1,\gamma_2]}{2s}+\frac{[\gamma_1,\rho_2]}{2s}&=O(s^{-1}).
\end{align}
In terms of $\Gamma:= \gamma_1+\ii \gamma_2$,  Eqs.~\eqref{Eq:gamma1} and \eqref{Eq:gamma2}, combine into
\begin{align}
\frac{d}{ds}\Gamma+\frac{[\ii\rho_3,\Gamma]}{2s}-\frac{\ii[\rho_+,\gamma_3]}{s}&=O(s^{-1}),
\end{align}
while Eq.~\eqref{Eq:gamma3} with the help of Eq.~\eqref{Eq:jacgamma} becomes
\begin{align}
\frac{d}{ds}\gamma_3+\frac{\ii[\rho_-,\Gamma]}{s}-\frac{[\ii\rho_3,\gamma_3]}{2s}&=O(s^{-1}).
\end{align}

Decomposing $\Gamma$ and $\gamma_3$ into $\ii\rho_3$ weight vectors of  irreducible $sl(2)$ representations, we have
\begin{align}
\frac{d}{ds}
\begin{pmatrix}
\Gamma^{\beta+2}\\ \gamma_3^\beta
\end{pmatrix}
+\frac{1}{s}
\begin{pmatrix}
\frac{\beta+2}{2} & -\ii\frac{d-1-\beta}{2}\\[3pt]
\ii\frac{d+1+\beta}{2} & -\frac{\beta}{2}
\end{pmatrix}
\begin{pmatrix}
\Gamma^{\beta+2}\\ \gamma_3^\beta
\end{pmatrix}=O(s^{-1}).
\end{align}
This has the form of Eq.~\eqref{Eq:minoreq} of Lemma~\ref{minor}, and the eigenvalues of its  coefficient matrix $A$ are $\frac{1}{2}(1\pm d)$.  Thus, for $d>1$ Lemma~\ref{minor} applies and both $\Gamma^\beta$ and  $\gamma_3^\beta$ are   bounded. Note, this conclusion holds with some notational modifications for the extreme cases when $\beta+2$ is the lowest weight or $\beta$ is the highest weight. 

Now, with all modes $\Gamma^\beta$ and  $\gamma_3^\beta$ of nontrivial $sl(2)$ representations bounded, we are left to consider the modes $\Gamma^{(0)}$ and  $\gamma_3^{(0)}$ in the trivial representation.  In that case Eqs.~(\ref{Eq:gamma1}-\ref{Eq:gamma3}) by direct integration imply that $\Gamma^{(0)}$ and  $\gamma_3^{(0)}$ are $O(\ln\frac{1}{s})$.  With this estimate and all other modes bounded, Eq.~\eqref{fineb1} (and its cyclic index permutations) imply (\ref{Eq:gamma1}-\ref{Eq:gamma3}) but with $O(s^{-\frac12}\ln\frac{1}{s})$ on the right-hand side, which now by direct integration implies that 
$\Gamma^{(0)}$ and  $\gamma_3^{(0)}$ are bounded.

\end{proof}

\subsection{The Moduli Spaces of the Small Bow Representation }\label{Sec:SmallRep}

An important example of a  bow representation is the {\em small representation}  denoted by $\mathfrak{s}.$  It has no $\lambda$-points, and therefore its fundamental data is empty. For the $A$-type  bow  its rank function  has constant rank 1: $R(s)\equiv 1$.  We denote the line bundle of $\mathfrak{s}$ over the bow by $\Fe,$ to distinguish it from the bundle $\FE$ associated to a general representation $\mathfrak{R}.$ We similarly  denote the bifundamental data and the Nahm data by the lower case symbols $b_\sigma$ and $t$ respectively to distinguish them from the data of a general representation.  

For the small bow representation, the Nahm equations \eqref{nahmeq0} - which follow from the moment map condition - imply the $t_j,$ $j=1,2,3$ are real constants, $t_j^\sigma$, on each $I_\sigma$. The moment map  conditions \eqref{pEnd} at, respectively, $p_\sigma-$ and $p_\sigma+$ read:  
\begin{align}\label{onemoment}
\rr_\sigma:=- \mathrm{Im}\, \ii b_\sigma b_\sigma^\dagger&=t^\sigma_1\e_1+ t^\sigma_2 \e_2+ t^\sigma_3 \e_3-\nu_\sigma\\
\rr_\sigma
=\mathrm{Im}\, \ii b_\sigma^c (b_\sigma^c)^\dagger&=t^{\sigma+1}_1\e_1+ t^{\sigma+1}_2 \e_2+ t^{\sigma+1}_3 \e_3 -\nu_\sigma.
\end{align}
If follows that $\vec{t}^\sigma=\vec{t}^{\sigma+1}$ and are  independent of $\sigma$; thus we introduce $\vec{t}:=\vec{t}^\sigma$ and observe that $\vec{r}_\sigma=\vec{t}-\vec{\nu}_\sigma.$ 

The restriction of the $L^2$ metric to the tangent space of the Nahm space is 
\begin{align}\sum_{j=1}^3ldt_j\otimes dt_j.
\end{align}
The level set, however, is still infinite-dimensional. Since we need a few details about this hyperk\"ahler quotient construction, in order to avoid analytical issues associated to quotients by infinite dimensional groups, we construct the quotient explicitly in two steps. 

\begin{figure}[htb]
\begin{center}
    \includegraphics[width=0.75\textwidth]{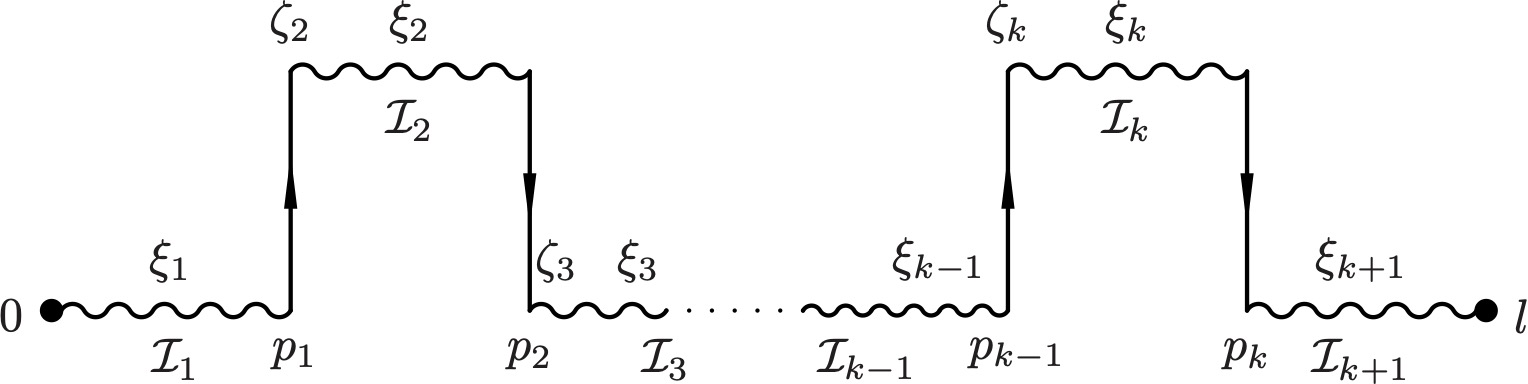}
\caption{The A-type bow over an interval $[0,l]$.}\label{Fig:Cut}
\end{center}
\end{figure}

Recall that the endomorphism bundle of a line bundle is canonically trivial, and all Hermitian connections on the line bundle induce the same connection on its endomorphism bundle, which is therefore a canonically trivial connection. With this observation, we note that the constant gauge transformations define a $U(1)$ subgroup of $\mathcal{G}(\mathfrak{s})$ that acts trivially on the bow data of $\mathfrak{s}.$   We choose to identify  $\mathcal{G}(\mathfrak{s})/U(1)$ (non-canonically) with 
 \begin{align}
 \mathcal{G}^0(\mathfrak{s}):=\left\{g\in\mathcal{G}(\mathfrak{s})\, |\, g(0)=1 \right\}.
 \end{align}
 As illustrated in Figure~\ref{fig:Rep}, the intervals of an A-type bow are identified with subarcs of a circle $S^1=\mathbb{R}/l\mathbb{Z}$. Since the $\mathcal{G}^0(\mathfrak{s})$ gauge group action is trivial at $s=0$, it is convenient to divide the interval containing that point into two sub-intervals. If needed, relabel the $p$-points so that $0<p_1<\ldots<p_k<l,$
and define intervals 
$\mathcal{I}_1=[0,p_1-], \mathcal{I}_{k+1}=[p_{k}+,l]$, and $\mathcal{I}_\sigma=I_\sigma=[p_{\sigma-1}+,p_\sigma-],$ for $ \sigma=2,3,\ldots,k,$ with  respective lengths $\ell_0=p_1, \ell_{k+1}=l-p_k$ and $\ell_\sigma=l_\sigma,$ for $\sigma=2,3,\ldots,k,$
 as illustrated in Figure~\ref{Fig:Cut}.

The gauge group $\mathcal{G}^0(\mathfrak{s})$ is a direct product
\begin{align}
 \mathcal{G}^0(\mathfrak{s})= \mathcal{G}^0_1\times\mathcal{G}_2\times\mathcal{G}_3\times\ldots\times\mathcal{G}_k \times\mathcal{G}^0_{k+1},
\end{align}
with each factor acting on the restriction of the bundle $\Fe$ to the corresponding interval $\mathcal{I}_\sigma.$ The $0$ superscript indicates the constraint $g(0) = 1= g(l)$. 
Consider the subgroup $\mathcal{H}_\sigma$ of the $\sigma$-th factor that acts trivially at the ends of the corresponding interval and has trivial winding:  
$$\mathcal{H}_\sigma=\left\{e^{\ii f_\sigma(s)}:\mathcal{I}_\sigma=[a,b]\rightarrow U(1)\, \Big|\, f_\sigma(a)=0= f_\sigma(b), \int_{\mathcal{I}_\sigma} f(s)ds=0 \right\}.$$ 
Then there are the following group isomorphisms:
\begin{align*}
&\mathbb{R}\times\mathcal{H}_0\xrightarrow{\sim}\mathcal{G}^0_1:&
&(\xi_1,e^{\ii f_1(s)})\mapsto e^{\ii s\xi_1/\ell_1}e^{\ii f_1(s)},\\
&\mathbb{R}\times\mathcal{H}_{k+1}\xrightarrow{\sim}\mathcal{G}^0_{k+1}:&
&(\xi_{k+1},e^{\ii f_{k+1}(s)})\mapsto e^{\ii (s-l)\xi_{k+1}/\ell_{k+1}}e^{\ii f_{k+1}(s)},\\
&U(1)\times\mathbb{R}\times\mathcal{H}_\sigma\xrightarrow{\sim}\mathcal{G}_\sigma:&
&(e^{\ii\zeta_\sigma},\xi_\sigma,e^{\ii f_\sigma(s)})\mapsto e^{\ii\zeta_\sigma}e^{\ii (s-p_{\sigma-1})\xi_\sigma/\ell_\sigma}e^{\ii f_\sigma(s)}, \\
&&&\qquad\qquad\qquad\qquad\qquad\qquad\qquad  \sigma=2, \ldots , k.
\end{align*}   
 
The first step in our quotient construction is the hyperk\"ahler reduction of Dat$(\mathfrak{s})=\mathcal{N}\times\mathcal{B}$ by $\mathcal{H}:=\times_{\sigma=1}^{k+1}\mathcal{H}_\sigma$, which acts only on the first factor $\mathcal{N}.$ Its hyperk\"ahler moment map is $(\frac{d t^\sigma_j}{ds}\e_j)_{\sigma}$; thus, on level sets of the moment map, the  $t^\sigma_j$ are constant on each $\mathcal{I}_\sigma$. 

The equivalence class of connections modulo the action of $\mathcal{H}$ has a unique representative connection in each equivalence class of the form $\nabla=\nabla^0_{\frac{d}{ds}}+it_0^\sigma$, with $t_0^\sigma$ a real constant. 
Moreover, the tangent space to the connections of this form is orthogonal  
to the orbits of $\mathcal{H}$. Thus, we may gauge fix to connections of this form and simply restrict the $L^2$ metric to this subset in order to obtain the hyperk\"ahler reduction with respect to $\mathcal{H}$. Thus $\mathcal{N}/\!/\!/\mathcal{H}=\mathbb{H}^{k+1}$ with the  
metric  descending from the $L^2$ metric on the Nahm data 
\begin{align}
\sum_{\sigma=1}^{k+1} \ell_\sigma  (dt_0^\sigma\otimes dt_0^\sigma+
dt_1^\sigma\otimes dt_1^\sigma+dt_2^\sigma\otimes dt_2^\sigma+dt_3^\sigma\otimes dt_3^\sigma).
\end{align}

Step two is the finite dimensional hyperk\"ahler reduction of the resulting hyperk\"ahler space 
$\mathrm{Dat}(\mathfrak{s})/\!/\!/\mathcal{H}=\mathbb{H}^{k+1}\times\mathcal{B} = \mathbb{H}^{k+1}\times(\mathbb{C}^2)^{k}$  by the action of the remaining $U(1)^{k-1}\times \IR^{k+1}$ gauge group. 
The Euclidean metric on each bifundamental $\mathbb{C}^2$ factor is 
$\sum_\sigma db_\sigma^\dagger db_\sigma=\sum_\sigma\frac{d\vec{r}_\sigma{}^{\,2}}{2r_\sigma}+2r_\sigma(d\varphi_\sigma {+}\eta_\sigma)^2.$  Here $\rr_\sigma:=-\mathrm{Im}\, \ii b_\sigma b_\sigma^\dagger\in \mathrm{Im}\, \mathbb{H}$,  $\varphi_\sigma$ is a local choice of phase angle for the circle fiber of the map $b_\sigma\to \rr_\sigma$, $\hat{\eta}_\sigma
=\ii \frac{b_\sigma^\dagger db_\sigma-db_\sigma^\dagger b_\sigma}{4 r_\sigma}
=d\varphi_\sigma+\eta_\sigma,$ and the one-form $\eta_\sigma$ satisfies  $d\eta_\sigma=*_3d\frac{1}{2r_\sigma}.$  (See Appendix~\ref{App:R4} for details.)
  Thus, the bow moduli space is
\begin{align}
\mathcal{M}_\mathfrak{s}^\nu:=\bm{\mu}_\mathfrak{s}^{-1}(\nu)/\mathcal{G}^0(\mathfrak{s})=(\mathbb{H}^{k+1}\times\mathbb{H}^k)/\!/\!/(\mathbb{R}^{k+1}\times U(1)^{k-1}),
\end{align}
 where $\mathbb{H}^{k+1}\times\mathbb{H}^k$ is equipped with the metric 
\begin{align}
\sum_{\sigma=1}^{k+1} \ell_\sigma (dt_0^\sigma\otimes dt_0^\sigma+d\vec{t}{}^{\,\sigma}\otimes d\vec{t}{}^{\,\sigma})+\sum_{\sigma=1}^{k+1}\frac{d\vec{r}_\sigma{}^{\,2}}{2r_\sigma}+2r_\sigma(d\varphi_\sigma  {+}\eta_\sigma)^2.
\end{align}

This presentation of the quotient is similar to but not quite in the form of the standard references \cite{Gibbons:1996nt} and  \cite{Witten:2009xu}; so, for the convenience of the reader, we will derive the metric, showing that the quotient is, in fact,  TN$_k^\nu$ (the multi-Taub-NUT space with $k$ NUTs at $\nu_1,\ldots,\nu_k$).
 
 The group element $(\xi_1,\ldots\xi_{k+1})\times(e^{\ii\zeta_2},\ldots,e^{\ii\zeta_{k}})\in\mathbb{R}^{k+1}\times U(1)^{k-1}$ acts by   
 \begin{align}\label{Eq:GrAction}
 t^\sigma_0&\mapsto t^\sigma_0+\frac{\xi_\sigma}{\ell_\sigma},& \sigma&=1,\ldots,k+1;\\
\varphi_1&\mapsto\varphi_1  {+} \xi_1  {-} \zeta_2,&&\nonumber\\
\varphi_\sigma&\mapsto\varphi_\sigma  {+} \xi_\sigma  {+} \zeta_\sigma  {-} \zeta_{\sigma+1},& \sigma&=2,\ldots,k-1,\nonumber\\
\varphi_k&\mapsto\varphi_k  {+} \xi_k  {+} \zeta_k  {+} \xi_{k+1} &&\nonumber
 \end{align}
 keeping $\vec{t}^\sigma$ and $\vec{r}_\sigma$ inert. On level $i\nu$ of the moment map with $\nu$ of the form \eqref{Level}, the
$\xi_\sigma$ action imposes the moment map condition 
$\vec{t}^\sigma-\vec{r}_\sigma = \vec{\nu}_\sigma$, 
while the $\zeta_\sigma$ action imposes the moment map condition  $\vec{r}_{\sigma-1}-\vec{r}_\sigma=\vec{\nu}_\sigma-\vec{\nu}_{\sigma-1}.$  It enforces $\vec{t}^\sigma=:\vec{t}$  equal across all intervals, and $\vec{r}_\sigma=\vec{t}-\vec{\nu}_\sigma.$  Thus, the metric on the level set is
 \begin{align}\label{Eq:LevelFin}
g^{lev}=Vd\vec{t}{\,}^2
+\ell_{k+1} (dt^{k+1}_0)^2
 +\sum_{\sigma=1}^k\left(\ell_\sigma (dt^\sigma_0)^2 
 +  2 r_\sigma(d\varphi_\sigma  {+} \eta_\sigma)^2\right),
 \end{align}
where we have set
$$V:= l+\sum_{\sigma=1}^k\frac{1}{2r_\sigma}.$$
This determines a metric on the quotient space by projecting onto the orthogonal complement of the gauge orbits.  The tangent space to the gauge orbit is spanned by 
$\{\frac{\p}{\p t_0^\sigma}  {+} l_\sigma \frac{\p}{\p \varphi^\sigma}\}_{\sigma=1}^{k}\cup 
\{\frac{\p}{\p t_0^{k+1}}+l_{k+1}\frac{\p}{\p \varphi^k}\}\cup
\{ \frac{\p}{\p \varphi^\sigma}-\frac{\p }{\p \varphi^{\sigma+1}} \}_{\sigma=1}^{k-1}$. 
Hence, the orthogonal complement of the tangent space to the gauge orbits is spanned by 
\begin{align}
&W_0:= \sum_{\sigma=1}^{k}\frac{1}{2r_\sigma }\frac{\p}{\p \varphi^\sigma}  {-} \sum_{\sigma=1}^{k+1} \frac{\p}{\p t_0^\sigma},&
&\text{ and  }& 
&\left\{W_j:= \frac{\p}{\p t_j}  {-} \eta_\sigma^j\frac{\p}{\p\varphi_\sigma} \right\}_{j=1}^3.
\end{align}

Define the local gauge invariant function\footnote{The reason for the negative signs appearing in this formula is that our chosen orientations were $(t_0,t_1,t_2,t_3)$ and $(r_1,r_2,r_3,\varphi).$} 
 $$\tau={-} \ell_{k+1}t^{k+1}_0+\sum_{\sigma=1}^k(\varphi_\sigma  {-} \ell_\sigma t^\sigma_0),$$ 
and set 
$$\eta:= \sum_{\sigma = 1}^k\eta_\sigma.$$
Then
$$(d\tau {+} \eta)(W_0)= V=|W_0|^2=|W_1|^2=|W_2|^2=|W_3|^2,$$
and 
$$(d\tau + \eta)(W_j)= 0, \forall j>0.$$

Hence $V^{-1}W_0$ is the horizontal lift of $\frac{\p}{\p\tau}$, and the quotient metric is 
\begin{align}
V d\vec{t}\,^{\,2}+\frac{(d\tau+\eta)^2}{V},
\end{align}
with $\eta$ satisfying $d\eta=*_3dV.$

Therefore the moduli space of the small representation, 
 $$\mathcal{M}_\mathfrak{s}^\nu:=\bm{\mu}_\mathfrak{s}^{-1}(\ii\nu)/\mathcal{G}^0(\mathfrak{s})$$ 
 at level $\ii\nu$ of Eq.~\eqref{Level}, is the $k$-centered Taub-NUT space, $\mathrm{TN}_k^\nu$, with the Taub-NUT centers at $\vec{\nu}_\sigma\in\mathrm{Im} \,\mathbb{H}=\mathbb{R}^3$.   (For $k=1$ this was demonstrated in \cite[Sec.3.2]{Cherkis:2009jm}.) 
This space is hyperk\"ahler with a triholomorphic isometric circle action, which has exactly $k$ fixed points, $\{\nu_1,\ldots,\nu_k\}$.  The quotient of TN$_k^\nu$ by this circle action is $\IR^3$, and the quotient map  
$\pi_k: \mathrm{TN}_k^\nu\to \IR^3$ defines a principal circle fibration in the complement of the fixed points:
\begin{equation}\label{fiber}
S^1\rightarrow\mathrm{TN}_k^\nu\setminus\{ \nu_1,\ldots, \nu_k\}\xrightarrow{\pi_k} \mathbb{R}^3\setminus\{\nu_1,\ldots,\nu_k\},
\end{equation}
with a globally defined Ehresmann connection one-form $\hat{\eta}$.  
A choice of a local section defines a fiber coordinate $\tau$, so $\hat{\eta}=d\tau+\pi_k^*\eta$. In turn, a connection one-form $\eta$ has  curvature $d\eta=*_3 dV.$ 
  This is a smooth Riemannian manifold, so long as all of the points $\vec{\nu}_\sigma$ are distinct.
 
\subsection{The Tautological Bundles}
Every hyperk\"ahler reduction produces a natural principal bundle over the quotient, so long as the group action on the level set is free.  In our case we have a family of associated line bundles $L_s$ over the quotient, which we now describe. See \cite{Witten:2009xu} and \cite{Cherkis:2008ip} for more details. 
The level set $\bm{\mu}_\mathfrak{s}^{-1}(\ii\nu)$  forms a principal $\mathcal{G}^0(\mathfrak{s})$ bundle over the moduli space $\mathcal{M}_\mathfrak{s}^\nu=\mathrm{TN}_k^\nu.$  We denote the quotient projection by 
\begin{align}\label{HKQ}
\rho: \bm{\mu}_\mathfrak{s}^{-1}(\ii\nu)\rightarrow \bm{\mu}_\mathfrak{s}^{-1}(\ii\nu)/\mathcal{G}^0(\mathfrak{s})=\mathrm{TN}_k^\nu.
\end{align}

Consider a point on the bow with coordinate $v$, contained in an interval $\mathcal{I}_\beta.$ Let 
$$\mathcal{G}^0_v(\mathfrak{s})=\{g\in\mathcal{G}(\mathfrak{s})\, |\, g(0)=1=g(v) \}$$
be the group of gauge transformations acting by identity at $s=0$ and at $s=v.$ Note, that the  quotient group $U(1)_v=\mathcal{G}^0(\mathfrak{s})/\mathcal{G}^0_v(\mathfrak{s}),$ naturally acts on the fiber $\Fe_v.$  
 Then, the partial quotient space 
\begin{align}\label{partq}\mathcal{P}_v:=\bm{\mu}_\mathfrak{s}^{-1}(\ii\nu)/\mathcal{G}^0_v(\mathfrak{s})
\end{align}
 is a $U(1)_v$ principal bundle over the moduli space 
 $$\mathrm{TN}_k^\nu=\bm{\mu}_\mathfrak{s}^{-1}(\ii\nu)/\mathcal{G}^0(\mathfrak{s})=\mathcal{P}_v/U(1)_v.$$  
 The trivial line bundle $\mathcal{P}_v\times\Fe_v\rightarrow\mathcal{P}_v$ is $U(1)_v$-equivariant; its quotient is the associated line bundle  $L_v.$   It changes continuously with $v$ except at $p$-points.  At a $p$-point we let  
\begin{align}K_\sigma:=L_{p_\sigma+}\otimes L_{p_\sigma-}^{-1}
\end{align} denote  the discontinuity, so that $L_{p_\sigma+}=   K_\sigma\otimes L_{p_\sigma-}.$  We can view the one parameter family of line bundles $\{L_u\}_{u\in \sqcup_\sigma I_\sigma}$ as a bundle $\mathcal{L}$ over 
$\mathcal{M}_\mathfrak{s}^\nu\times \sqcup_\sigma I_\sigma$ in a natural manner with
\begin{align}
\mathcal{L}:=\bm{\mu}_\mathfrak{s}^{-1}(\ii\nu){\times \sqcup_\sigma I_\sigma}\times_{\mathcal{G}^0(\mathfrak{s})} \Fe = (\bm{\mu}_\mathfrak{s}^{-1}(\ii\nu){\times \sqcup_\sigma I_\sigma}\times \Fe)/ \mathcal{G}^0(\mathfrak{s}),
\end{align}
where $\mathcal{G}^0(\mathfrak{s})$ 
acts as $g(p,v,z) = ( pg,v, g^{-1}(v)z).$
This bundle will play a role in our work similar to the Poincar\'e bundle in the Nahm-Fourier-Mukai transform. The restriction $\mathcal{L}_s$ of this bundle to $\mathcal{M}_\mathfrak{s}^\nu\times \{s\}$ is the line bundle $L_s\rightarrow \mathcal{M}_\mathfrak{s}^\nu$ discussed above. We next derive the natural connection on $\mathcal{L}_v$.

\begin{figure}[htb]
\begin{center}
    \includegraphics[width=0.75\textwidth]{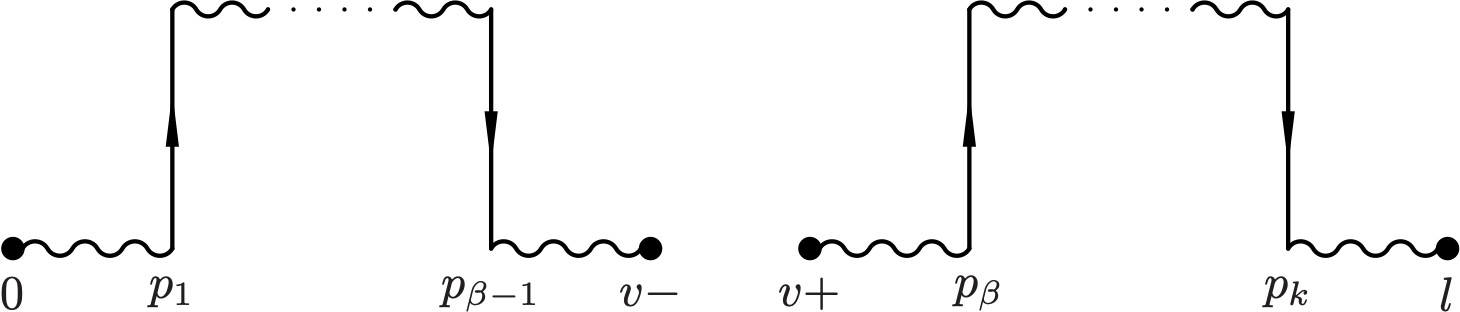}
\caption{An A-type bow cut at $s=0$ and at $s=v.$}
\label{DoubleCut}
\end{center}
\end{figure}

We cut the bow at two points $s=0$ and $s=v$, as in Figure~\ref{DoubleCut}. The result is two bows, left and right, each with its own affine space of data, Dat$_l$ and Dat$_r$, respectively. 
The group factors as $\mathcal{G}^0_v(\mathfrak{s})=\mathcal{G}_l\times\mathcal{G}_r,$ with each factor acting on the data of the corresponding bow. 
As computed in Sec.~\ref{Sec:SmallRep}, the metric on this partial quotient $(\mathrm{Dat}_l\times\mathrm{Dat}_r)/\!/\!/\mathcal{G}^0_v=(\mathrm{Dat}_l/\!/\!/\mathcal{G}_l)\times(\mathrm{Dat}_r/\!/\!/\mathcal{G}_r)$ is
 \begin{align}\label{Eq:Partial}
V_l d\vec{t}_l^{\,\,2}+\frac{(d\tau_l{+}\eta_l)^2}{V_l}
 + V_r d\vec{t}_r^{\,\,2}+\frac{(d\tau_r{+}\eta_r)^2}{V_r},
 \end{align}
 with \begin{align}V_l=v+\sum_{\sigma=1}^{\beta-1}\frac{1}{2|\vec{t}_l-\vec{\nu}_\sigma|}\text{  and }V_r=l-v+\sum_{\sigma=\beta}^{k}\frac{1}{2|\vec{t}_r-\vec{\nu}_\sigma|},
\end{align} and corresponding $\eta_l$ and $\eta_r$ one-forms.

This partial hyperk\"ahler reduction provides a hyperk\"ahler quotient space $\bm{\mu}^{-1}(\ii\nu)/\!/\!/\mathcal{G}^0_v(\mathfrak{s})$ which inherits a triholomorphic action of 
 $U(1)_v=\mathcal{G}^0(\mathfrak{s})/\mathcal{G}^0_v(\mathfrak{s}),$  which acts by $(\tau_l,\tau_r)\mapsto(\tau_l  {-} \epsilon,\tau_r  {+} \epsilon)$, while keeping $\vec{t}_l$ and $\vec{t}_r$ inert.  The resulting action leaves invariant the function  $\tau:=\tau_l+\tau_r$; its moment map is $\vec{t}_l-\vec{t}_r$.  The vanishing moment map condition provides a five-dimensional level set $\mathcal{P}_v$ with metric (which we write in terms of the invariant coordinates $\vec{t}:=\vec{t}_l=\vec{t}_r$,  $\tau$, and $\tau_l$):
\begin{align}
&Vd\vec{t}^2+\frac{(d\tau_l  {+} \eta_l)^2}{V_l}+\frac{(d\tau_r  {+} \eta_r)^2}{V_r}\nonumber\\
=&Vd\vec{t}^2 +\frac{(d\tau  {+} \eta)^2}{V}+\frac{V}{V_l V_r}\left(d\tau_l  {+} \eta_l-\frac{V_l}{V}(d\tau  {+} \eta)\right)^2,
\label{Eq:Taut}
\end{align}
with $V=V_l+V_r$ and $\eta=\eta_l+\eta_r.$  
Hence 
in the trivialization given by the local coordinate $\tau_l$, the  $U(1)$ bundle $\mathcal{P}_v$ over the $k$-centered Taub-NUT space has the Ehresmann connection form $d\tau_l  {+} \eta_l-\frac{V_l}{V}\hat{\eta}.$ The associated line bundle $\Fe_v$ inherits the $d+\ii a_v$ with the connection one-form
\begin{align}\label{Eq:Conn-e}
a_{v} :=V_l(v)\frac{d\tau+\eta}{V}-\sum_{\sigma=1}^{\beta-1}\eta_\sigma.
\end{align}
For example, the horizontal lift of $\frac{\p}{\p\tau}$ to $\mathcal{P}_v$ is therefore
\begin{align}(\frac{\p}{\p\tau})^H= \frac{\p}{\p\tau}+ \frac{V_l(v)}{V}\frac{\p}{\p\tau_l}.
\end{align}

\section{The Up Transform}\label{Sec:Up}
Let $(\nabla,\tilde\T,Q,B)$ be a solution  of the moment map equations for the large bow representation $\mathfrak{R}$. We say $(\nabla,\tilde\T,Q,B)$ is Without Abelian  Factors  (WAF), if  there is no  subbundle $\FE_{ab}$ of  $\FE$ that is preserved by $\nabla$ and restricted to which, the $T_j$ satisfy the trivial  Nahm equation  $\dot T_1=0= [T_2,T_3]$, and cyclic permutations.  For any WAF moment map solution of $\mathfrak{R}$, we now associate a Dirac operator $\D_\mathfrak{R}.$   Moreover, we twist it to produce a family of Dirac operators  $\D_{t,b}$ parameterized by the small bow representation level set $\bm{\mu}_\mathfrak{s}^{-1}(\ii\nu)\ni(t_0,\vec{t},b).$ The index bundle of this family is equivariant under the action of the gauge group $\mathcal{G}^0(\mathfrak{s})$ of the small representation and therefore descends to a bundle $\mathcal{E}\rightarrow\mathrm{TN}_k^\nu=\bm{\mu}_\mathfrak{s}^{-1}(\ii\nu)/\mathcal{G}^0(\mathfrak{s})$. 
This is the instanton bundle. In this section we define the bundle $\mathcal{E}$ and compute its rank, leaving the examination of its induced connection, which we denote $d_A$, to 
Section~\ref{Sec:AsymInst}.  Anti-self-duality of the curvature of $d_A$ is proved in 
Appendix~\ref{ASDproof}. Note that  we orient TN$_k^\nu$ using the volume form $Vdt^1\wedge dt^2\wedge dt^3\wedge d\tau$. 
\subsection{The Bow Dirac Operator}\label{BDDef}

Writing the connection $\nabla$ as $\nabla=\nabla^0_{\frac{d}{ds}}+\ii T_0,$ it is convenient to assemble the Nahm data into an endomorphism valued  quaternion 
$$\T:=1\otimes T_0+\e_1\otimes {T}_1+ \e_2\otimes {T}_2+\e_3\otimes {T}_3.$$
 Its adjoint is 
  $\T^\dagger=1\otimes T_0^*+\e^\dagger_1\otimes {T}^*_1+ \e^\dagger_2\otimes {T}^*_2+\e^\dagger_3\otimes {T}^*_3=1\otimes T_0-\sum_j\e_j\otimes T_j$. 
We associate to any element $(\nabla,\tilde\T,Q,B)$ of the large bow representation $\mathfrak{R}$ a bow Dirac operator  as follows.
\begin{define} Let 
$H^1_{term}(\S\otimes \FE)$ denote the subspace of the direct sum  $\oplus_\sigma H^1(I_\sigma,\S\otimes \FE_\sigma)$ with vanishing terminating components. Define 
the {\em bow Dirac operator}
\begin{align}\D_\mathfrak{R}: H^1_{term}(\S\otimes \FE)\to L^2(\S\otimes \FE)\oplus\mathop{\oplus}\limits_{\lambda\in\Lambda^0} W_\lambda\oplus\mathop{\oplus}\limits_{\sigma=1}^k \left(\FE_{p_{\sigma}-}\oplus \FE_{p_{\sigma}+}\right)
\end{align}
by
\begin{align}
\D_\mathfrak{R}:=
\begin{pmatrix}
-\nabla^0_{\frac{d}{ds}}-{ \ii \T^\dagger}
\\ -Q_\lambda^*\, \ev_\lambda \\ B_\sigma^\dagger\, \ev_{p_{\sigma}-}\\  -{B_\sigma^c}^\dagger\, \ev_{p_{\sigma}+}
\end{pmatrix} ,
\end{align}
where $\ev$ denotes the evaluation map, $\ev_yf:= f(y).$ 
\end{define}
The formal adjoint of $\D_\mathfrak{R}$ is given by
\begin{align}
\D_\mathfrak{R}^\dagger=\nabla^0_{\frac{d}{ds}}{+\ii \T}-\sum_{\lambda\in\Lambda^0}\delta_{\lambda}Q_\lambda+\sum_\sigma\Big(\delta_{p_{\sigma}-}{B}_\sigma-\delta_{p_{\sigma}+}B^c_\sigma\Big),
\end{align}
 which maps  $L^2(\S\otimes \FE)\oplus\mathop{\oplus}\limits_{\lambda\in\Lambda^0} W_\lambda\oplus\mathop{\oplus}\limits_{\sigma=1}^k \left(\FE_{p_{\sigma}-}\oplus \FE_{p_{\sigma}+}\right)$ to $\mathop{\oplus}\limits_{\sigma=1}^k H^{-1}(I_\sigma,\S\otimes \FE_\sigma).$
 
A direct calculation readily verifies that in the notation of \eqref{imdef}, the moment map  density of a bow representation \eqref{Moment} can be expressed as  
\begin{align}
\mu_\mathfrak{R}=-\mathrm{Im}\,\D_\mathfrak{R}^\dagger\D_\mathfrak{R}.
\end{align}

In the same fashion, any  $(t_0,\tilde{\t},b)\in \bm{\mu}^{-1}_{\mathfrak{s}}(\ii\nu)$ also has an associated bow Dirac operator  
${\mathpzc{d}_\mathfrak{s}}.$ Since $\Lambda =\emptyset$ for $\mathfrak{s}$, the fundamental   terms containing $Q$ are absent in this case.  Hence the formal adjoint of ${\mathpzc{d}_\mathfrak{s}}$ is  $\mathpzc{d}_\mathfrak{s}^\dagger=\nabla^0_{\frac{d}{ds}}+\ii\, t_0+\ii\, \tilde{\t}+\sum_\sigma(\delta_{p_\sigma-}b_\sigma-\delta_{p_\sigma+}b_\sigma^c)$ acting on  $L^2(\S\otimes \Fe)\oplus\oplus_\sigma  (\Fe_{p_\sigma+}\oplus\Fe_{p_\sigma-}).$ This defines a family of operators parameterized by the level set $\bm{\mu}^{-1}_{\mathfrak{s}}(\ii\nu)$ of the small bow representation at level $\ii\nu.$  Note, that the gauge group of the small bow representation acts naturally on the kernel of this family.

We use the family $ \mathpzc{d}_{\mathfrak{s}}^\dagger$ parameterized by  $\bm{\mu}^{-1}_{\mathfrak{s}}(\ii\nu)$ to twist the large bow Dirac operator $\D_\mathfrak{R}$. To be exact, we need the charge conjugate family $\mathpzc{d}_{\mathfrak{s}^*}^\dagger$ with domain  $L^2(\S\otimes \Fe^*)\oplus\oplus_\sigma  (\Fe_{p_\sigma+}^*\oplus\Fe_{p_\sigma-}^*)$ and range
$H^{-1}(\S\otimes \Fe^*)$ which we now define. 
(The reason for using $\Fe^*$ rather than $\Fe$ is analogous to the usual Fourier transform where the kernel $e^{-\ii p x}$ of the inverse Fourier transform is conjugate to the kernel  $e^{\ii p x}$ of the direct Fourier transform.  
Similarly, $\Fe^*$ serves a role analogous to an integral kernel for the Up transform,  which we are formulating now, while $\Fe$ serves as the kernel for the Down transform described in \cite{Third}.)  

Let $\mathbf{f}=\left(\begin{smallmatrix} 
 f  \\  
n_+\\
n_-
\end{smallmatrix} \right)\in L^2(\S\otimes \Fe)\oplus\oplus_\sigma  (\Fe_{p_\sigma+}\oplus\Fe_{p_\sigma-})$ and write its charge conjugate as ${\mathbf{f}}^c=\left(\begin{smallmatrix} 
 f^c  \\  
-n_-^*\\
n_+^*
\end{smallmatrix} \right)\in L^2(\S\otimes \Fe^*)\oplus\oplus_\sigma  (\Fe_{p_\sigma-}^*\oplus\Fe_{p_\sigma+}^*).$ 
Then, as 
\begin{align}
 \mathpzc{d}_{\mathfrak{s}}^\dagger \mathbf{f}=(\frac{d}{ds}+\ii\t)f+\sum_\sigma (\delta_{p_\sigma-} b_\sigma n_+-\delta_{p_\sigma+}b^c n_-),
\end{align}
its charge conjugate operator is defined (using the notation of Sec.~\ref{Sec:ChargeConj}) via 
\begin{multline}
 \mathpzc{d}_{\mathfrak{s}^*}^\dagger \mathbf{f}^c:
 =\left( \mathpzc{d}_{\mathfrak{s}}^\dagger \mathbf{f}\right)^c
  =\epsilon\left( \mathpzc{d}_{\mathfrak{s}}^\dagger \mathbf{f}\right)^*\\
 =(\frac{d}{ds}-\ii\,^r\!\t)f^c+\sum_\sigma (\delta_{p_\sigma-}\, ^r\!b_\sigma^c n_+^*-\delta_{p_\sigma+}\,^r\!(b^c)^c n_-^*)\\
=\frac{d}{ds}f^c - \ii\,f^c t_0 - \sum_{j=1}^3\ii\e_j f^c t_j +\sum_\sigma (\delta_{p_\sigma+}\,n_-^* b 
+\delta_{p_\sigma-}\, n_+^* b_\sigma^c )\\
=\left((\frac{d}{ds}-\ii\,^r\!\t)+\sum_\sigma (- \delta_{p_\sigma+}\,^r\!b
+\delta_{p_\sigma-}\, ^r\!b_\sigma^c)\right)\mathbf{f}^c.
\end{multline}
Thus, the charge conjugate operator family\footnote{Since for the A-type bow considered here the small representation is one-dimensional the transposition $^r$ above can be ignored.  We  keep it nevertheless to have our formulas apply for all bows.} 
is
\begin{equation}\label{twist0}
  \mathpzc{d}_{\mathfrak{s}^*}^\dagger=\dds-\ii \,^r\!t_0 - \ii\,^r\!\tilde{\t}+\delta_{p_{\sigma -}}\, { ^r\!b^c_\sigma}
  -\,  
  \delta_{p_{\sigma +}}\, ^r\!{b}_\sigma.
 \end{equation} 
 Significantly, if $(\t, b)$ is a solution of the small bow representation at level $\ii\nu$, then 
  since $\ii\nu=\mu_\mathfrak{s} =-\mathrm{Im}\,  \mathpzc{d}_{\mathfrak{s}}^\dagger \mathpzc{d}_{\mathfrak{s}}$,  its charge conjugate relation is 
  $-\ii\nu =(\ii \nu)^c=\mu_\mathfrak{s}^c=-\mathrm{Im}\,  \mathpzc{d}_{\mathfrak{s}^*}^\dagger \mathpzc{d}_{\mathfrak{s}^*}.$

Now we are ready to consider the family of {\em twisted bow Dirac operators} 
$$\D_{t,b}:=\D_\mathfrak{R} +  \mathpzc{d}_{\mathfrak{s}^*} $$ 
acting on sections in $H^1_{term}(\S\otimes\FE\otimes\Fe^*)$. 
The adjoint operator, $\D_{t,b}^\dagger=\D_\mathfrak{R}^\dagger +  \mathpzc{d}_{\mathfrak{s}^*}^\dagger$,  acts on 
$$L^2(\S\otimes\FE\otimes\Fe^*)\oplus\mathop{\oplus}\limits_{\lambda\in\Lambda^0}W_\lambda\otimes \Fe^*_\lambda\oplus\mathop{\oplus}\limits_{\sigma=1}^k N_\sigma(\FE,\Fe^*),$$ 
where we have set 
\begin{align}\label{Eq:Ndef}
N_\sigma(\FE,\Fe^*):= (\FE_{p_{\sigma}+}\otimes\Fe^*_{p_{\sigma}-})\oplus(\FE_{p_{\sigma}-}\otimes\Fe^*_{p_{\sigma}+}).
\end{align}
We have used the fact that an element of $\mathrm{Hom}(W_\lambda,\FE_\lambda)$ also defines an element of  $\mathrm{Hom}(W_\lambda\otimes \Fe^*_\lambda,\FE_\lambda\otimes \Fe^*_\lambda)$, since $\Fe^*_\lambda$ is one dimensional.

More explicitly, $\D_{t,b}^\dagger=\D_\mathfrak{R}^\dagger+ \mathpzc{d}_{\mathfrak{s}^*}^\dagger$ acts on 
$$\D_{t,b}^\dagger:L^2(\S\otimes \FE\otimes\Fe^*)\oplus \mathop{\oplus}_{\lambda\in\Lambda^0}W_\lambda\otimes \Fe^*_\lambda\oplus\mathop{\oplus}_\sigma N_\sigma(\FE,\Fe^*)\rightarrow H^{-1}(S\otimes \FE\otimes\Fe^*),$$
as 
\begin{multline}\label{twistedDirac}
\D_{t,b}^\dagger
 \left(\begin{smallmatrix}
 \chi\\ w_\lambda\\ \nu^+_\sigma \\ \nu^-_\sigma
 \end{smallmatrix}\right)
 =\dds\chi+\ii (T^0\chi-\chi t^0)+\ii \e_j(T^j\chi-\chi t^j)
 -\mathop{\sum}\limits_{\lambda\in\Lambda^0}\delta_\lambda Q_\lambda w_\lambda\\
 +\mathop{\sum}\limits_{\sigma}(\delta_{p_{\sigma -}}B_\sigma \nu^+_\sigma
-\,
 \delta_{p_{\sigma +}}b_\sigma\nu^+_\sigma 
 +\delta_{p_{\sigma -}}b^c_\sigma\nu^-_\sigma -\delta_{p_{\sigma +}}B^c \nu^-_\sigma),
\end{multline}
where $\chi\in L^2(\S\otimes \FE\otimes\Fe^*), w_\lambda\in W_\lambda\otimes \Fe^*_\lambda,$ and $n^+_\sigma\in \FE_{p_\sigma+}\otimes\Fe_{p_\sigma-}$ and $n^-_\sigma\in \FE_{p_\sigma-}\otimes\Fe_{p_\sigma+}.$ 
 
Because we have chosen the same level sets $\ii\nu\otimes 1_{\Fe}$ of  ${\mu}_{\mathfrak{s}}$ and  $\ii\nu\otimes 1_{\FE}$ of ${\mu}_{\mathfrak{R}}$, the twisted Dirac operator satisfies the crucial property : 
\begin{equation}\label{lapl}
-\mathrm{Im}\,  \D_{t,b}^\dagger \D_{t,b}=\mu_\mathfrak{R}\otimes 1_{\Fe^*}+(1_\FE\otimes \mu_{\mathfrak{s}})^c =\ii\nu\otimes 1_\FE\otimes 1_{\Fe^*} +(\ii \nu\otimes 1_\FE\otimes 1_{\Fe^*})^c =0.
\end{equation}
Thus, $\D_{t,b}^\dagger \D_{t,b}$ is 
  quaternionic real; i.e. it commutes with the action of the quaternionic units. 
The WAF condition we imposed on the solution $(\nabla,\tilde\T,Q,B)$ implies that $\D_{t,b}$ has zero $L^2$ kernel. (Lemma \ref{WAF} of Sec.~\ref{Sec:GreenFn}). The index bundle of the family is therefore ${\rm Ker}\, \D_{t,b}^\dagger\rightarrow \mu_{\mathfrak{s}}^{-1}(\ii\nu).$ This bundle is equivariant under the action of $\mathcal{G}^0(\mathfrak{s});$  thus it descends to a bundle on the quotient space: $\mathcal{E}\rightarrow \textrm{TN}_k^\nu.$  This bundle has a connection $A$ induced by the natural connection on the $\Fe^*$ family. As proved in Appendix~\ref{ASDproof} and in \cite[Sec.7]{Cherkis:2009jm}, this connection $A$ has anti-self-dual curvature.

\subsection{An Index Theorem on the Bow}\label{BowIndex}


To lighten our notation, since the exact form of certain linear operators is not significant for this section, let $\TT:=\T-{}^r\!\t,$ and work in the gauge in which $T_0=0$ and $t_0=0,$ so that $\ii\TT$ is Hermitian.  

We will compute the dimension of $\mathrm{Ker}\, \D^\dagger_{t,b},$ which is the space  of  solutions $(x,w,\nu)$ of
\begin{align}
&\left(\dds+\ii\TT\right)x=0,\label{sol1}\\
&x(\lambda+)-x(\lambda-)=Q_\lambda w_\lambda,\,\forall\lambda\in \Lambda^0\label{sol2}\\
&x(p_\sigma+)=\mathcal{A}_\sigma\nu_\sigma,\label{sol3}\\
&x(p_\sigma-)=\mathcal{B}_\sigma\nu_\sigma.\label{sol4}
\end{align}
Here $x$ is an $L^2$ section of $\S\otimes \FE\otimes\Fe^*$, $w_\lambda\in W_\lambda\otimes \Fe^*_\lambda,$ and $\nu_\sigma\in N_\sigma(\FE,\Fe^*).$  

Importantly, near a $\lambda$-point, with $\lambda\in\Lambda^-\cup\Lambda^+,$ 
$$\ii\TT(s)=\frac{1}{2(s-\lambda)}
\begin{pmatrix}
\sum_{j=1}^3\e_j\otimes{\bm\rho}_j &0\\
0 & 0
\end{pmatrix}
+O\left((s-\lambda)^0\right),$$ 
with $\{\e_j\}_{j=1}^3$ and $\{{\bm\rho}_j\}_{j=1}^3$ the images of a unitary basis of $su(2)$ under the  irreducible representations of $su(2)$ of  dimensions 2 and $m:=|R(\lambda+)-R(\lambda-)|>0$, respectively, acting on spaces $V^2$ and $V^m$.   The representation space $V^2\otimes V^m$ of terminating components decomposes  as $V^2\otimes V^m=V^{m+1}\oplus V^{m-1}$. Expressing $\sum_{j=1}^3 \e_j\otimes{\bm\rho}_j$ as a sum of Casimirs via  $2\sum_{j=1}^3 \e_j\otimes{\bm\rho}_j=\sum_{j=1}^3(\e_j\otimes 1+1\otimes{\bm\rho}_j)^2-\e_j^2\otimes 1-1\otimes{\bm\rho}_j^2$, with the Casimir in $m$ dimensional representation $\sum_j{\bm\rho}_j^2=-m^2+1$, it follows that 
$\sum_{j=1}^3 \e_j\otimes{\bm\rho}_j$ takes the value $1-m$ on the $V^{m+1}$ summand and $1+m$ on the $V^{m-1}$ summand. (See, for example \cite[(2.7)]{Hitchin:1983ay}.) We correspondingly write 
\begin{align}\label{zerodirac} 0 &= \left(\dds+\ii\TT\right)x \nonumber\\
&=\left(\begin{array}{ccc}\dds+\frac{m+1}{2(s-\lambda)}+b_{11}&b_{12}&b_{13}\\
b_{21}&\dds-\frac{m-1}{2(s-\lambda)}+b_{22}&b_{23}\\b_{31}&b_{32}&\frac{d}{ds}+b_{33}\end{array}\right)\left(\begin{array}{c}x_{m-1}\\
x_{m+1}\\x_c\end{array}\right),
\end{align}
where the $b_{ij}$ are bounded, $x_{m\pm 1}\in V^{m\pm 1}$, and $x_c\in (V^2\otimes V^m)^\perp.$ 
	If we denote by $Q_0$ the orthogonal projection of a fiber of $\S\otimes \FE\otimes \Fe^*$ onto the subspace $\S\otimes \FE^{cont}\otimes \Fe$ of its continuous components, then $x_c=Q_0 x$. In turn, we denote by $Q_+$ the projection of this fiber on the $m+1$ eigenspace of $\sum_j \e_j \otimes\bm{\rho}_j$ and by $Q_-$ the projection on its $1-m$ eigenspace.  Then, $x_{m+1}=Q_+x$ and $x_{m-1}=Q_- x. $ 
These projections are used extensively below.

The following proposition is implicit in the literature. (See \cite[(2.7)]{Hitchin:1983ay}.) We sketch a proof for completeness.  
\begin{proposition}\label{impthm} Let 
$\left(\begin{smallmatrix}x_{m-1}\\x_{m+1}\\x_c\end{smallmatrix}\right)$ be a solution to \eqref{zerodirac}. Then $y_{m-1}:= 
\lim_{s\to\lambda}|s-\lambda|^{\frac{m+1}{2}}x_{m-1}$ exists. Moreover, if $y_{m-1} = 0$, 
then the limit $\lim_{s\to \lambda}(x_{m-1}(s),x_{m+1}(s)) = (0,0)$,   $x_c$ is bounded,  and $y_{c}:= \lim_{s\to\lambda} x_{c}$ exists. If $y_c = 0$, then $|s-\lambda|^{\frac{1-m}{2}}(|x_{m-1}|+|x_{m+1}|+|x_c|)$ is bounded and $\lim_{s\to\lambda}|s-\lambda|^{\frac{1-m}{2}}x_{m+1}$ exists. 
\end{proposition}
\begin{proof}
  Without loss of generality, consider the case $\lambda = 0\in \Lambda^-$ and $s\leq 0$. Then for some $B<\infty$, we have 
\begin{align}& \frac{d}{ds}|x_{m-1}|\leq \frac{m+1}{2|s|}|x_{m-1}|+B(|x_{m+1}|+|x_{m-1}|+|x_c|),\label{23}\\
& \frac{d}{ds}|x_{m+1}| +\frac{m-1}{2|s|}|x_{m+1}|\leq  B(|x_{m+1}|+|x_{m-1}|+|x_c|),\text{ and}\label{24}\\
& \frac{d}{ds}|x_{c}| \leq  B(|x_{m+1}|+|x_{m-1}|+|x_c|).\label{25}
\end{align}
Hence, $\frac{d}{ds}(|x_{m+1}|+|x_{m-1}|+|x_c|)\leq (\frac{m+1}{2|s|}+B)(|x_{m+1}|+|x_{m-1}|+|x_c|),$ and integration from some nonzero $s_0$ to $s$ yields for some $C_1>0$
\begin{align}
|x_{m+1}|+|x_{m-1}|+|x_c|\leq C_1|s|^{-\frac{m+1}{2}}.
\end{align}
 Suppose now that we have shown for some $1\leq j\leq \frac{m+1}{2}$ and some $c_j>0$
\begin{align}
|x_{m+1}|+|x_{m-1}|+|x_c|\leq c_j|s|^{-j}.
\end{align}
Then integrating \eqref{24} plus \eqref{25} yields for some $c_{j,2} >0$
\begin{align}
|x_{m+1}| +|x_c|\leq c_{j,2}\left(|s|^{1-j}+\delta_{j1}\ln(\frac{1}{|s|})\right).
\end{align}
 Replacing \eqref{23} by the sharper 
\begin{align}\label{26}|\frac{d}{ds} \left(|s|^{\frac{m+1}{2}}e^{-|s|B}|x_{m-1}|\right)|\leq  B|s|^{\frac{m+1}{2}}(|x_{m+1}|+|x_c|),\end{align}
and integrating, yields $|s|^{\frac{m+1}{2}}|x_{m-1}|$ is bounded and the limit  $y_{m-1}:= \lim_{s\to 0}|s|^{\frac{m+1}{2}}x_{m-1}$ exists. If this limit is zero, then integrating from $0$ to $s$  yields $|x_{m-1}|\leq c_{3,j}(|s|^{1-j}+\delta_{j1}\ln(\frac{1}{|s|}))$ and therefore $|x_{m+1}|+|x_{m-1}|+|x_c| \leq c_{4,j}(|s|^{1-j}+\delta_{j1}\ln(\frac{1}{|s|})).$ By induction, we deduce that if $y_{m-1} =0$, 
then 
$$|x_{m+1}|+|x_{m-1}|+|x_c| \leq c_{5}\ln(2+|s|^{-1}).$$
Iterating the estimates one more time  yields $|x_{m+1}|+|x_{m-1}|+|x_c|$ is bounded. Now  integrate   \eqref{26}  from $0$ to $s$ instead of from $s_0$ to $s$ to obtain  $\lim_{s\to 0}x_{m-1}(s) = 0$. Integrate   \eqref{23}    to obtain  $\frac{|x_{m+1}|}{ |s|}\in L^1((s_0,0]) $ and $\lim_{s\to 0}x_{m+1}(s) = 0$. Equation \eqref{zerodirac} then implies  $y_{c}:= \lim_{s\to 0} x_{c}(s)$ exists.  

Suppose that for some $j\in( -\frac{1}{4},\frac{m-2}{2}]$, we have 
\begin{align}
 |x_{m+1}|+|x_{m-1}|+|x_c| \leq \tilde c_j|s|^{j}. 
\end{align}
Then integrating inequalities \eqref{25} and \eqref{26} from $0$ to $s$ imply 
$|x_{m-1}| +|x_c| =O(|s|^{1+j}).$ Replace \eqref{24} by the equivalent inequality 
\begin{multline}\label{modeq}
\frac{d}{ds}(|s|^{-j-\frac{1}{2}}|x_{m+1}|) +\frac{m-2-2j}{2|s|}|s|^{-j-\frac{1}{2}}|x_{m+1}|\\
\leq  B|s|^{-j-\frac{1}{2}}(|x_{m+1}|+|x_{m-1}|+|x_c|).
\end{multline}
Integrating \eqref{modeq} from $s_0\not = 0$ to $s$, yields $|x_{m+1}| = O(|s|^{\frac{1}{2}+j}).$ Hence 
$|x_{m+1}|+|x_{m-1}| +|x_c| =O(|s|^{\frac{1}{2}+j}).$ Iterating, we deduce 
$|x_{m+1}|+|x_{m-1}| +|x_c| =O(|s|^{\frac{m-1}{2}}),$ and $\lim_{s\to 0}|s|^{\frac{1-m}{2}}x_{m+1}$ exists. 
\end{proof} 
Reversing some signs in the preceding proof, we have the following proposition. 
\begin{proposition}\label{impthmc} Let 
$\tilde x = \left(\begin{smallmatrix}\tilde x_{m-1}\\\tilde x_{m+1}\\\tilde x_c\end{smallmatrix}\right)$ be a solution to 
\begin{align}\label{zerodiracstar}0  = \left( \dds-\ii\TT\right)\tilde x.
\end{align} 
Then $\tilde y_{m+1}:= 
\lim_{s\to\lambda}|s-\lambda|^{\frac{m-1}{2}}\tilde x_{m+1}$ exists. If $\tilde y_{m+1} = 0$,
 then the limit $\lim_{s\to 0}\tilde x_{m-1}(s)=0$,   and $\tilde x_c$ is bounded,  and $\tilde y_{c}:= \lim_{s\to\lambda} \tilde x_{c}$ exists. 
 If $\tilde y_c = 0$, then $|s-\lambda|^{\frac{-1-m}{2}}(|\tilde x_{m-1}|+|\tilde x_{m+1}|+|\tilde x_c|)$ is bounded and $\lim_{s\to\lambda}|s-\lambda|^{\frac{-1-m}{2}}\tilde x_{m-1}$ exists. 
\end{proposition}
\begin{corollary}\label{termlife} There exists an $m-1$ dimensional space of unbounded solutions to \eqref{zerodirac}. 
There exists an $m+1$ dimensional space of unbounded solutions to \eqref{zerodiracstar}. 
\end{corollary}
\begin{proof}
The dimension of the space of solutions to \eqref{zerodiracstar} which vanish at $\lambda$ is greater than or equal to $m-1$ by Proposition \ref{impthmc}. On the other hand, given a solution $x$ to \eqref{zerodirac} and $\tilde x$ to
 \eqref{zerodiracstar}, we have 
\begin{align}\label{invar}  \frac{d}{ds}\langle x,\tilde x\rangle = \langle -\ii\TT x,\tilde x\rangle + \langle x,\ii\TT\tilde x\rangle
=0.
\end{align}
Since there exists an $m-1$ dimensional space of solutions to \eqref{zerodiracstar} vanishing at $\lambda$, there must exist a space of solutions of \eqref{zerodirac} of the same dimension which blow up.  Similarly there exists an $m+1$ dimensional space of solutions to \eqref{zerodiracstar} which blow up.
\end{proof}

\subsubsection{Matching Spaces}\label{Sec:Ys}
Each $\lambda\in\Lambda$ is the right end of some subinterval $I_\sigma^\alpha=[a(\lambda),\lambda]$ and the left end of some other subinterval $I_\sigma^{\alpha+1}=[\lambda,b(\lambda)].$  In other words, $a(\lambda)$ is the $\lambda$- or $p$-point immediately preceding $\lambda$, while  $b(\lambda)$ is the $\lambda$- or $p$-point immediately following $\lambda$ (along the right circle of Fig.~\ref{fig:Rep}).

For $\lambda\in\Lambda\setminus\Lambda^0=\Lambda^-\cup\Lambda^+$, define the space of continuous sections
\begin{align}\label{Eq:Ylambdapole}
Y_\lambda=\left\{y\in \Gamma((a(\lambda),b(\lambda)),\S\otimes \FE\otimes \Fe^*)\, \Big|\, (\dds-\ii\TT)y=0,  y^{term}(\lambda)=0\right\}.
\end{align}
The space of all solutions to $(\dds-i\TT)y=0$ on the disjoint union of the two intervals $(a(\lambda),\lambda)$ and $(\lambda, b(\lambda))$ has dimension $2R(\lambda-)+2R(\lambda+)$. By Corollary \ref{termlife}, the conditions required for such solutions to lie in $Y_\lambda$  impose  $m+1$ conditions by requiring the terminal components to vanish and an additional $2R_{min}(\lambda):=2\min\{R(\lambda+),R(\lambda -)\}$ conditions by requiring the  continuing components to be continuous. Hence the dimension of $Y_\lambda$ is  $2R(\lambda+) + 2R(\lambda-) - (m+1)-2R_{min}(\lambda) =  (m-1)+2R_{min}(\lambda)=  R(\lambda+) +  R(\lambda-)-1.$

For $\lambda\in\Lambda^0$ define the space of continuous sections
\begin{align}\label{Eq:YlambdaQ}
Y_\lambda=\left\{y\in \Gamma\left(a(\lambda),b(\lambda)\right),\S\otimes \FE\otimes \Fe^*)\, \Big|\, (\dds-\ii\TT)y=0,\  Q_\lambda^\dagger y(\lambda)=0\right\}.
\end{align}
Its dimension is $\mathrm{dim}\,(S\otimes E_\lambda) -1=R(\lambda+)+R(\lambda-)-1,$ so long as $Q\neq0.$ This gives a universal formula applicable to all $\lambda$-points:
\begin{align}\label{dimYs}
\dim Y_\lambda=R(\lambda+)+R(\lambda-)-1.
\end{align}

For a $p$-point, similarly, let $a(p)$ denote the $\lambda$- or $p$-point immediately preceding $p$, and let $b(p)$ be the $\lambda$- or $p$-point immediately following $p$.  Let $Y_p=Y_{p-}\oplus Y_{p+},$ with 
\begin{align}\label{Eq:Ypm}
Y_{p-}=\left\{y\in \Gamma \left((a(p),p],\S\otimes \FE\otimes \Fe^*\right)\, \Big|\, (\dds-\ii\TT)y=0\right\},\\
\label{Eq:Ypp}
Y_{p+}=\left\{y\in \Gamma \left([p,b(p)),\S\otimes \FE\otimes \Fe^*\right) \, \Big|\,  (\dds-\ii\TT)y=0\right\}.
\end{align}
The dimension of $Y_p$ is $2 R(p-)+2 R(p+).$
\begin{lemma}\label{rep} Let 
$\lambda\in \Lambda^+\cup \Lambda^-$. The map  $v:Y_\lambda\to (Q_0+Q_+)(\S\otimes \FE\otimes \Fe^*)_\lambda $ given by $v(y)= Q_0y(\lambda) + \lim_{s\to\lambda}(s-\lambda)^{\frac{-1-m}{2}}Q_+y(s)$ is an isomorphism. 
\end{lemma}
\begin{proof}
 Since the two vector spaces have the same dimension, it suffices to show the map is injective. So, suppose $v(y) = 0.$ Then Proposition \ref{impthmc} implies that $|(s-\lambda)^{\frac{-1-m}{2}}y(s)|$ is bounded. Let $x\in Ker((\dds+\ii\TT)).$ Then 
\begin{align}\langle x,y\rangle(s_0) &= \lim_{s\to\lambda}\langle x,y\rangle(s)\nonumber\\
& = \langle \lim_{s\to\lambda}(s-\lambda)^{\frac{1+m}{2}}x(s),\lim_{s\to\lambda}(s-\lambda)^{\frac{-1-m}{2}}y(s)\rangle =0,
\end{align}
because $\lim_{s\to\lambda}(s-\lambda)^{\frac{1+m}{2}}x(s)\in \mathrm{Im}\, Q_+,$ and $\lim_{s\to\lambda}(s-\lambda)^{\frac{-1-m}{2}}Q_+y(s)=0$ by hypothesis. 
\end{proof}
\subsubsection{Interior Spaces}\label{Sec:Spaces}
In the interior $\mathring{I}_\sigma^\alpha$ of a subinterval $I_\sigma^\alpha$, define  the space 
\begin{align}\label{Eq:Xsp}
X_\sigma^\alpha=\left\{x\in\Gamma(\mathring{I}_\sigma^\alpha,\S\otimes \FE\otimes \Fe^*) \, \Big|\,  (\dds+\ii\TT)x=0 \right\},
\end{align}
with {\em no boundary conditions} imposed. This space has dimension $2 R_\sigma^\alpha.$  Let $X:=\oplus_{\alpha,\sigma}X_\sigma^\alpha,$ $N:=\oplus_\sigma N_\sigma(\FE,\Fe^*)$ and $Y:=\oplus_{\lambda\in\Lambda}Y_\lambda\oplus\oplus_p Y_p.$  
Define maps $F_\lambda: X\rightarrow Y_\lambda^\vee$   by
\begin{align}\label{Eq:Flam}
F_\lambda(x)(y_\lambda):= \langle x,y_\lambda\rangle_{\frac{\lambda+b(\lambda)}{2}} - \langle x,y_\lambda\rangle_{\frac{\lambda+a(\lambda)}{2}}
=\langle x,y_\lambda\rangle_{\lambda+\epsilon} - \langle x,y_\lambda\rangle_{\lambda-\epsilon}.
\end{align}
Due to Eq.~\eqref{invar} the inner product, $\langle x,y\rangle_s := \langle x(s),y(s)\rangle$,  is independent of $s$ in the interior $\mathring{I}_\sigma^\alpha$ of each subinterval. 

Define additional linear maps $F_{p\pm}:X\oplus N\rightarrow Y_{p\pm}^\vee$  by 
\begin{align}\label{Eq:Fpm}
F_{p-}(x,\nu)(y_p)&=
\langle\mathcal{B}\nu,y_p\rangle_{p-}-\langle x,y_p\rangle_{\frac{p+a(p)}{2}}
=\langle \mathcal{B}\nu,y_p\rangle_{p-}-\langle x,y_p\rangle_{p-},\\
\label{Eq:Fpp}
F_{p+}(x,\nu)(y_p)&=
\langle x,y_p\rangle_{\frac{p+b(p)}{2}}-\langle \mathcal{A}\nu,y_p\rangle_{p+}
=\langle x,y_p\rangle_{p+}-\langle \mathcal{A}\nu,y_p\rangle_{p+}.
\end{align}
Add these to obtain a map 
\begin{align}\label{Delta}
\Delta:X\oplus N&\rightarrow Y^\vee=\mathop{\oplus}_\lambda Y_\lambda^\vee\oplus\mathop{\oplus}_p (Y_{p-}^\vee\oplus  Y_{p+}^\vee),
\end{align}
with
\begin{multline}
 \Delta(x,\nu)(y): =\sum_\lambda F_\lambda(x)(y_\lambda)+\sum_p (F_{p-}(y_p)+F_{p+}(y_p))\\
 =\sum_{\lambda\in\Lambda}(\langle x,y_\lambda\rangle_{\lambda+\epsilon}-\langle x,y_\lambda\rangle_{\lambda-\epsilon})+\sum_{p\in P}(\langle x-\mathcal{A}\nu,y_p\rangle_{p+}-\langle x-\mathcal{B}\nu,y_p\rangle_{p-}).
\end{multline}
Observe that by construction, for each open interval $\mathring{I}_\sigma^\alpha$ there are 2 components of $y$ whose support contains $\mathring{I}_\sigma^\alpha$. For example if $x$ is supported in $(\lambda, \beta)$, then 
$$\Delta(x,0)(y) =  \langle x,y_\lambda\rangle_{\frac{\lambda+\beta}{2}}-\langle x,y_{\beta}\rangle_{\frac{\lambda+\beta}{2}}.$$
\begin{lemma}\label{kers}
$\mathrm{Ker}(\Delta)=\mathrm{Ker}(\D^\dagger_{t,b}).$
\end{lemma}
\begin{proof}
Every element of $X$ satisfies $(\dds+i\TT)x=0$. Hence we only need to check the explicit and implicit boundary conditions \eqref{sol2} - \eqref{sol4}. For $y\in Y$, $y(p+)$ and $y(p-)$ are independent. Hence $\Delta(x,\nu)(y)=0, \forall y\in Y$, implies 
$x(p+) = \mathcal{A}\nu_p$ and $x(p-) = \mathcal{B}\nu_p.$ 

 Consider $\lambda\in \Lambda^+$. Because the space of solutions to $(\dds+i\TT)x=0$ and $(\dds-i\TT)\tilde x=0$ have the same dimension on $(a(\lambda),\lambda)$, and because the inner product is locally constant, $\langle x,\tilde x\rangle_s$ is a perfect pairing between  the two spaces, for $s\not = \lambda$.  Consider the subspace $Z_1$ of solutions to $(\dds-i\TT)\tilde x=0$ satisfying  $\lim_{s\to\lambda} \tilde x_c(s)=0$ and $\lim_{s\to\lambda}|s-\lambda|^{\frac{m-1}{2}}\tilde x_m(s))=0.$ By Proposition \ref{impthmc}, $Z_1\subset Y_\lambda$ and $\dim(Z_1) = m+1$. Hence for all $y\in Z_1$, $\langle x,y\rangle_s = 0$ for all $x\in \mathrm{Ker}(\Delta).$ The annihilator of $Z_1$ is comprised of the bounded solutions of $(\dds+i\TT) x=0$. In particular, their terminating components vanish. For such elements taking the limit as $\epsilon\to 0$ gives for $y\in Y_\lambda$,  
\begin{multline}
0 = \langle x,y\rangle_{\lambda+\epsilon}-\langle x,y\rangle_{\lambda-\epsilon} = \langle x_c(\lambda +),y_c(\lambda+)\rangle -\langle x_c(\lambda -),y_c(\lambda -)\rangle\\
 = \langle x_c(\lambda +) -  x_c(\lambda -),y_c(\lambda )\rangle,
 \end{multline}
since $y_c$ is assumed continuous. Therefore $x_c(\lambda +) =  x_c(\lambda -)$ and $x_c$ is continuous. Similarly, for $\lambda\in \Lambda^-$, the terminating components vanish and $x_c$ is continuous. 

We are left to consider $\lambda\in \Lambda^0$. For $x\in \mathrm{Ker}(\Delta)$, we have for all $y\in \mathrm{Ker}(Q_\lambda^\dagger) $, 
$$\langle   x(\lambda+)  -      x(\lambda-) ,y \rangle = 0.$$
Hence $x(\lambda+)  -      x(\lambda-)\in \mathrm{Ker}(Q_\lambda^\dagger)^\perp,$ and $x(\lambda+)  -      x(\lambda-)\in  \mathrm{Im}(Q_\lambda),$ as desired. 
\end{proof}
 
\begin{lemma}\label{cokers}
$\mathrm{CoKer}(\Delta)$ is isomorphic to $\mathrm{Ker}(\D_{t,b}).$
\end{lemma}
\begin{proof}Because 
$Y$ is finite dimensional, the dimension of $\mathrm{CoKer}(\Delta)$ is equal to the dimension of the annihilator of 
$\mathrm{Im}(\Delta)$ in $Y^{\vee\vee}\simeq Y.$ Thus we consider $y\in Y$ which satisfies $\Delta(x,\nu)(y) = 0$, $\forall (x,\nu)\in X$. 

The condition $0= \Delta(0,\nu)(y)$ $\forall \nu$ implies $\mathcal{B}_p^\dagger y(p-)=\mathcal{A}_p^\dagger y(p+).$
Arguing as in the preceding lemma, we see that the terminating component of $y$ vanishes and $y_c$ is bounded. Hence we are left to check that on overlapping domains the different components of $y$ agree. 
 Consider all $x$ supported in a single sub interval of the form $( \lambda ,\beta(\lambda))$.  Then 
$$0= \Delta(x,0)(y) = \langle x(\frac{\lambda+\beta(\lambda)}{2}),y_\lambda(\frac{\lambda+\beta(\lambda)}{2})-y_{\beta(\lambda)}(\frac{\lambda+\beta(\lambda)}{2})\rangle.$$   
Hence $y_\lambda = y_{\beta(\lambda)}$ in this subinterval and similarly agree on all overlapping domains. Thus the elements of $Y\cap \text{annihilator}\,\mathrm{Im}(\Delta)$ define a single element $y$ of $H^1_{term}(\S\otimes \FE\otimes\Fe^*)$ which satisfies $(\dds-i\TT)y=0$ on each subinterval. Thus 
$Y\cap \text{annihilator}\,\mathrm{Im}(\Delta) = \mathrm{Ker}(\D_{t,b}).$
\end{proof}
\begin{corollary}\label{equalindex}
\end{corollary}
 $\Delta$ is surjective, and $\text{Index}(\Delta) = -\text{Index}(\D_{t,b})=\text{dim }\mathrm{Ker}(\D_{t,b}^*).$
\begin{proof}
By Lemma \ref{fred}, $\D_{t,b}$ is Fredholm, and its cokernel is isomorphic to the kernel of $\D_{t,b}^*$. 
The equality of the indices now follows from Lemmas \ref{kers} and \ref{cokers}. 
As we observed following Eq.~\eqref{lapl}, $\mathrm{Ker}(\D_{t,b})=\{0\}$; therefore, $\Delta$ is surjective, and  
$\text{Index}(\D_{t,b}^*)=\text{dim }\mathrm{Ker}(\D_{t,b}^*).$
\end{proof}

\begin{theo}\label{Thm:Index} 
$$\dim\mathrm{Ker}(\D^\dagger_{t,b})=|\Lambda|.$$
\end{theo}
\begin{proof}
By Corollary \ref{equalindex}, $\text{dim }\mathrm{Ker}(\D_{t,b}^*)=\text{Index}(\Delta)  =\dim (X\oplus N) - \dim (Y^\vee).$ Computing the dimensions yields: 
\begin{align*}
\dim (X\oplus N)&=\sum_{\sigma,\alpha} 2R_\sigma^\alpha+\sum_\sigma(R(p_\sigma+)+R(p_\sigma-)),\\
\dim (Y^\vee)&= \dim(Y)= \sum_\lambda(R(\lambda+) +R(\lambda-)-1)\\ 
&\qquad\qquad\qquad  +2\sum_\sigma(R(p_\sigma+)+ R(p_\sigma-)).
\end{align*}
Since each interval $I_\sigma^\alpha$ has two ends and the rank is constant on the interior of each interval, 
$$2\sum_{\sigma,\alpha} R_\sigma^\alpha=\sum_\lambda(R(\lambda+)+R(\lambda-))+\sum_\sigma(R(p_\sigma+)+ R(p_\sigma-)).$$
Thus the contribution of the rank functions cancel and the difference in dimensions  is  $|\Lambda|.$
\end{proof}

\section{The  Instanton Connection}\label{Sec:AsymInst}
In this section, we recall the definition of the induced connection on the index bundle $\mathcal{E}$ and study its asymptotic behavior as $|t\,|\to\infty.$ 

The vector fields generating the  action of the group $\mathcal{G}^0(\mathfrak{s})$ form a vertical distribution $\mathcal{V}$ in the tangent space of the level set $\bm{\mu}_{\mathfrak{s}}^{-1}(\ii\nu).$  As a subset of a hyperk\"ahler affine space, $\bm{\mu}_{\mathfrak{s}}^{-1}(\ii\nu)$ inherits a metric.  The horizontal distribution $\mathcal{H}$ consists of the orthogonal complement in  $T\bm{\mu}_{\mathfrak{s}}^{-1}(\ii\nu)$ of the vertical distribution: $\mathcal{H=V^\perp}.$ Since the normal space of the level set is spanned by $\e_1\mathcal{V}+\e_2\mathcal{V}+\e_3\mathcal{V},$ one can also view $\mathcal{H}$ as the orthogonal complement of the quaternionic subbundle $\mathcal{V}+\e_1\mathcal{V}+\e_2\mathcal{V}+\e_3\mathcal{V}$ of the tangent space to Dat$(\mathfrak{s}).$ Therefore, $\mathcal{H}$ is also quaternionic, and the tangent space to the quotient inherits this  quaternionic structure.  In fact, this is the very reason why the quotient $\bm{\mu}_{\mathfrak{s}}^{-1}(\ii\nu)/\mathcal{G}^0(\mathfrak{s})$ is hyperk\"ahler.

The kernel of $\D^\dagger_{t,b}$ defines  a smoothly varying  family of constant rank subspaces of $L^2(\S\otimes \FE\otimes\Fe^*) \oplus\,\oplus_{\lambda\in \Lambda^0}W_\lambda\otimes\Fe^*_\lambda\oplus\,\oplus_\sigma N_\sigma(\FE,\Fe^*) $, parameterized by $(t,b)\in\mu_{\mathfrak{s}}^{-1}(\ii\nu)$. It therefore defines a Hermitian subbundle $\mathcal{\tilde E}=\mathrm{Ker}\,\D^\dagger_{t,b}$ of the trivial Hilbert bundle 
$\mu_{\mathfrak{s}}^{-1}(\ii\nu)\times L^2(\S\otimes \FE\otimes\Fe^*) \oplus\,\oplus_{\lambda\in \Lambda^0}W_\lambda\otimes\Fe^*_\lambda\oplus\,\oplus_\sigma N_\sigma(\FE,\Fe^*).$ Moreover, both this subbundle and the trivial bundle are $\mathcal{G}^0(\mathfrak{s})$-equivariant.  We use the natural action of 
$\mathcal{G}^0(\mathfrak{s})$ on the base $\mu_{\mathfrak{s}}^{-1}(\ii\nu)$ and the action on $L^2(\S\otimes \FE\otimes\Fe^*) \oplus\,\oplus_{\lambda\in \Lambda^0}W_\lambda\otimes\Fe^*_\lambda\oplus\,\oplus_\sigma N_\sigma(\FE,\Fe^*) $ inherited from $\Fe^*.$

As a subbundle of a trivial Hilbert  bundle, $\mathcal{\tilde E}$ inherits a connection $\tilde{\nabla}$.  Then the covariant derivative $\nabla$ of a section $f$ of the quotient bundle $\mathcal{E}$ with respect to a vector field $X$ on TN$_k^\nu$ is defined to be the covariant derivative $\tilde{\nabla}$ of the corresponding  equivariant representative $\tilde{f}$ with respect to the horizontal lift $X^H$ of $X:$ 
$\nabla_X f (p)=\tilde{\nabla}_{X^H}\tilde{f} (x)$ with $x$ in the $\mathcal{G}^0(\mathfrak{s})$-orbit of $p$.
We will henceforth drop the notational distinction between $f$ and its equivariant representative, denoting both $f$. 

Let $\left\{\psi_i(t,b) = \left(\begin{smallmatrix}\hat\psi_i\\w_i\\n_{-,i}\\n_{+,i}\end{smallmatrix}\right)\right\}_{i=1}^n,$ with $n=|\Lambda|$, be a smoothly varying $\mathcal{G}^0(\mathfrak{s})$-equivariant family of unitary bases of $\mathrm{Ker}(\D^\dagger_{t,b}),$ providing a local unitary frame of   $\tilde{\mathcal{E}}$ over $\mu_{\mathfrak{s}}^{-1}(\ii\nu)$. Here 
$\hat \psi_i\in L^2(\S\otimes \FE\otimes\Fe^*)$, $w_i\in \oplus_{\lambda\in \Lambda^0}W_\lambda\otimes\Fe^*_\lambda$,  $n_{-,i}\in \oplus_\sigma\FE_{p_{\sigma+}}\otimes \Fe^*_{p_{\sigma -}}$ and $n_{+,i}\in \oplus_\sigma\FE_{p_{\sigma-}}\otimes \Fe^*_{p_{\sigma +}}$.

In this frame, the connection matrix $A$ is 
$$A_{ii'}(X)=\langle  \psi_i, X^H\psi_{i'}\rangle,$$
where $X^H$ is the horizontal lift of $X$, and  we write $X^H\psi_i$ for $\nabla_{X^H}^{\text{triv}}\psi_i$, where $\nabla^{\text{triv}}$ denotes the trivial connection on the Hilbert bundle.

Consider
$ (\frac{\p}{\p\tau})^H  \psi_i=  \left(\begin{smallmatrix} 
(\frac{\p}{\p\tau})^H_s \hat\psi_i(s)\\ 
(\frac{\p}{\p\tau})^H_\lambda w_i\\ 
(\frac{\p}{\p\tau})^H_{p-} n_{-,i}\\ 
(\frac{\p}{\p\tau})^H_{p+} n_{+,i}
\end{smallmatrix}\right)$. In order to compute $ (\frac{\p}{\p\tau})^H \hat\psi_i(v),$ we observe that $\mathcal{G}_v^0(\mathfrak{s})$ leaves $\hat\psi_i(v)$ fixed. Hence we can be computing $ X^H \hat\psi_i(v)$ on the partial quotient $\mathcal{P}_v$ (defined in \eqref{partq}) instead of lifting to $\bm{\mu}_{\mathfrak{s}}^{-1}(\ii\nu)$. From \eqref{Eq:Taut}, the horizontal lift of $\frac{\p}{\p\tau}$ to $\mathcal{P}_v$ is 
\begin{align}\label{Eq:Vell}
\left(\frac{\p}{\p\tau}\right)^H_v &= \frac{\p}{\p\tau} + \frac{V_l(v)}{V}\frac{\p}{\p\tau_l},& &\text{with}&
  V_l(v)&:=v+\mathop{\sum}\limits_{\substack{\sigma\\ p_\sigma<v}}\frac{1}{2|\vec{t}-\nu_\sigma|}.
  \end{align}     
Hence 
\begin{align}\label{conn1}
A_{ii'}(\frac{\partial}{\partial\tau})
=&\left\langle{  \psi}_i(s), \left(\frac{\partial}{\partial\tau}+\frac{V_l(s)}{V}\frac{\partial}{\partial\tau_l}\right)
{  \psi}_{i'}(s)\right\rangle_{L^2},\end{align}
where the $L^2$ inner product is with respect to the Lebesgue measure on  each $I_\sigma$ plus the Dirac measure on each $\lambda-$ and $p-$point. We also have 
\begin{align}\label{tjh}
\left(\frac{\p}{\p t^j}\right)^H_v  = \frac{\p}{\p t^j}
-\left(\sum\limits_{\substack{\sigma\\p_\sigma<v}}\eta_\sigma^j
            -\frac{V_l(v)}{V}\eta^j\right)\frac{\p}{\p\tau_l}.
\end{align}
Here $\eta_\sigma^j=\eta_\sigma(\frac{\p}{\p t^j})$ and $\eta^j=\eta(\frac{\p}{\p t^j})$ are the components of the respective forms.

The main goal of this section is proving the following theorem
\begin{theo}\label{Thm:Asymp}
There is a local frame of the bundle $\mathcal{E}={\rm Ker}_{L^2}D^\dagger_{t,b}\to \mathrm{TN}_k$ in which 
\begin{equation}\label{goal_again}
A(\frac{\p}{\p\tau} )= -\ii\, {\rm diag}\frac{{\lambda_{i}+\frac{\hat{m}_{i}}{2|\vec{t}\,|}}}{V}+O(|\vec{t}\,|^{-2})_{C^1},
\end{equation}
where $\hat{m}_{i}=R(\lambda_i+)-R(\lambda_i-)+\#\{p_\sigma | 0\leq p_\sigma<\lambda_i\},$ and  
\begin{equation}\label{goal-1_again}
A(\frac{\p}{\p t^j})=  O(|\vec{t}\,|^{-2})_{C^1}.
\end{equation}
\end{theo}
 The proof of Theorem~\ref{goal} constitutes Section \ref{Sec:AsymCon}. 
Note, the integer $\hat{m}_i$ associated with $\lambda_i$   coincides with the D3-brane linking number, as defined in \cite{Witten:2009xu}. It equals  the rank change plus the number of $p$-points to the left of $\lambda_i.$

As an immediate corollary, we have 
\begin{corollary}\label{Coro:quad}
$|F_A|=  O(|\vec{t}\,|^{-2})$, and therefore $F_A\in L^2$. 
\end{corollary}
Thus the connection induced on $\mathcal{E}$ is, indeed, an instanton, i.e. it has an $L^2$ anti-self-dual curvature.

 In fact, as proved in \cite[Thm.21]{First}, any connection on TN$_k^\nu$ with $L^2$ anti-self-dual curvature has the asymptotic form of \eqref{goal}, and, as we will prove in \cite{Third}, any $L^2$ anti-self-dual connection on TN$_k^\nu$ with generic asymptotic holonomy can be constructed in this way.
 
\subsection{The Green's Function of the Bow Laplacian}\label{Sec:GreenFn}
Let $\mathfrak{G}_{t,b}$ denote the Green's function of $\D^\dagger_{t,b}\D_{t,b}$. In this subsection we prove that the operator norm of $\mathfrak{G}_{t,b}$ has quadratic decay in $t$. 

Let $f\in {\rm Dom}\,\D_{t,b}^\dagger\D_{t,b}$, then $f=\oplus f^\alpha_\sigma$ with $f^\alpha_\sigma\in H^1_{term}(S\otimes E^\alpha_\sigma)$. Recall that functions in ${\rm Dom}\,\D_{t,b}$ have vanishing terminating components and continuous continuing components at $\lambda$-points. Denote by $\langle\cdot,\cdot\rangle$ the $L^2$-pairing. We use a gauge in which all $T_{0}=0$ and the $t_{0}=0$. In this gauge we have
\begin{multline}\label{Green2}
\|\D_{t,b} f\|^2=\|\dds f\|^2+\sum_{j=1}^3 \|(T_{j } -t_{j}) f  \|^2
+\| {B}_\sigma^\dagger f(p_{\sigma}-) - f(p_{\sigma}+) {b}^\dagger\|^2_{\FE_{\sigma+}\otimes \Fe^*_{\sigma-}} \\
+\| {B^c}_\sigma^\dagger f(p_{\sigma}+) - f(p_{\sigma}-) {b^c}^\dagger\|^2_{\FE_{\sigma-}\otimes \Fe^*_{\sigma+}}
+\sum_{\lambda\in\Lambda^0}\| Q^\dagger_\lambda f(\lambda)\|^2_{\FE_{\lambda}\otimes \Fe^*_{\lambda}}.
\end{multline}
\begin{lemma}\label{fred}
$\D_{t,b}$ is Fredholm.
\end{lemma}
\begin{proof}Suppose 
$\{f_n\}\subset {\rm Dom}\,\D_{t,b}^\dagger\D_{t,b}$ is an $H^1$ orthonormal sequence such that 
$\lim_{n\to\infty}\|\D_{t,b}f_n\|^2_{L^2} =0$. Then there is a subsequence $\{f_{1,n}\}$ which converges in $L^2$ to a limit $F$ which is covariant constant and has $L^2$ norm $1$. The sequence $\{f_{1,n}-f_{1,n+1}\}$ then converges to zero, contradicting orthonormality. Hence no such orthonormal sequence exists. Therefore, $\text{Ker}(\D_{t,b})$ is finite dimensional, and $\exists \lambda_1>0$ such that 
\begin{align}\label{poincare}\|\D_{t,b}f\|^2_{L^2} \geq \lambda_1 \| f\|^2_{H^1}, \forall f\perp
\text{Ker}(D_{t,b}).
\end{align} 
By \eqref{poincare}, the image of $\D_{t,b}$ is closed. Hence the cokernel is isomorphic to the kernel of the adjoint operator. The latter is finite dimensional since it is an ordinary differential operator, and the result follows. 
\end{proof}
We now examine the kernel of $\D_{t,b}$.
\begin{lemma}\label{WAF}If $(\nabla,\tilde\T,Q,B)$ is WAF, then 
$\text{Ker}(\D_{t,b})=0.$
\end{lemma}
\begin{proof}Let 
$f\in \text{Ker}(D_{t,b}).$ Then, from \eqref{Green2}, on each interval $I_\sigma$,  $f_\sigma:=f|_{I_\sigma}$ is covariant constant, and 
$(T_j-t_j)f_\sigma=0$, for $ j=1,2,3$. If $f_\sigma\not = 0,$ then there is a $\nabla$-stable subbundle over $I_\sigma$ on which $T_j = t_j$, contradicting the WAF condition.
\end{proof}
Set $$C_\Lambda:= \max_{\lambda\in\Lambda}\{ 2+2|R(\lambda+)-R(\lambda-)|\}.$$
\begin{lemma}\label{LemSup} 
Let $f\in \mathrm{Dom}\D_{t,b}^\dagger\D_{t,b}$, then for $|\vec{t}\,|$ sufficiently large, 
\begin{equation}\label{super}
\|\D_{t,b}f\|^2\geq  C_\Lambda^{-1}|\vec{t}\,|^2\|f\|^2.
\end{equation}
\end{lemma}
\begin{proof}
In order to find a lower bound for $\|\D_{t,b} f\|^2$, it suffices, by (\ref{Green2}), to find a lower bound for $\sum_{j=1}^3 \|(T_{j } -t_{j}) f  \|^2$. The definition of the Nahm affine space implies that, for any given $(T_1,T_2,T_3)$, there exist  constants $c_1>0$ and $\epsilon>0$, so that 
\begin{itemize}
\item$\sum_j|T_j(s)|<c_1$, if $|s-\lambda|\geq \epsilon, \forall\lambda\in \Lambda\setminus \Lambda^0$, 
\item  for each $\lambda\in \Lambda\setminus \Lambda_0$, if $|s-\lambda|< \epsilon,$ conditions~(\ref{Pole1}) and (\ref{Pole2}) provide  $\{{\bm\rho}_j\}_{j=1}^3$ (depending on $\lambda$) forming an $su(2)$ irreducible representation of dimension $|R(\lambda+)-R(\lambda-)|$ and highest weight $\mu=|R(\lambda+)-R(\lambda-)|-1,$ so that 
$\sum_j|T_j(\lambda+s) + \frac{\ii{\rho}_j}{2s}|<c_1$. The Casimir of this representation is $\sum_{j=1}^3({\ii\bm\rho}_j)^2=\mu^2+2\mu.$
\end{itemize}
 When $\sum_j|T_j(s)|<c_1$,  we have for $|\vec{t}|$ large:
\begin{align}\label{dumbest1}\sum_{j=1}^3  |(T_{j } -t_{j}) f   |^2 \geq  (|\vec{t}|-c_1)^2|f|^2\geq \frac{1}{2}|t|^2|f|^2.\end{align}

Consider $s$ in a  neighborhood where $\sum_j|T_j(\lambda+s)+ \frac{\ii{\rho}_j}{2s}|^2<c_1^2.$ 
Then we have
\begin{align}\label{preup}
& \sum_{j=1}^3  | \frac{\ii{\rho}_j+t_j}{2s}   f   |^2\nonumber\\
&= \sum_{a}   (\frac{\sqrt{\mu^2+2\mu}}{2s}-\frac{w_a|\vec{t}\,|}{\sqrt{\mu^2+2\mu}})^2+(1-\frac{w_a^2}{\mu^2+2\mu})|\vec{t}\,|^2)|\pi_af|^2 \nonumber\\
&= \sum_{a}   (\frac{ \mu^2+2\mu -w_a^2}{4s^2}+(\frac{w_a }{2s}-|\vec{t}\,|)^2 )|\pi_af|^2.
\end{align}
where we have used the fact that the Casimir of an $su(2)$ representation with highest weight $\mu$ is $\mu^2+2\mu$. 

Next, using  \eqref{preup}  for $|\vec{t}\,|$ sufficiently large,
\begin{align}\label{finest2}
\sum_{j=1}^3  |(T_j(s) -t_{j}) f   |^2&= \sum_{j=1}^3|(\frac{\ii{\rho}_j}{2s} + t_{j}) f - (\frac{\ii{\rho}_j}{2s} + T_{j}(s)) f |^2\\\nonumber
&\geq   ( \sqrt{\frac{2 }{\mu+2}}|\vec{t}\,|-  c_1)^2 |f|^2.
\end{align}
Inequalities \eqref{dumbest1}, \eqref{finest2}, and the obvious extension to the $\mu=0$ case imply the desired result. 
\end{proof}
This lemma implies an upper bound on the norm of the Green's function. 
\begin{lemma}\label{Greendecay}
 For $|\vec{t}\,|$ sufficiently large, $\|\mathfrak{G}_{t,b}\|\leq \frac{C_\Lambda}{|\vec{t}\,|^2}$.
\end{lemma}
\begin{proof} 
{$$\frac{|\vec{t}\,|^2}{C_\Lambda}\|\mathfrak{G}_{t,b}F\|^2\leq \|\D_{t,b} \mathfrak{G}_{t,b}F\|^2=\langle F,\mathfrak{G}_{t,b}F\rangle\leq \|F\|\|\mathfrak{G}_{t,b}F\|.$$
Hence  
$$\|\mathfrak{G}_{t,b}F\| \leq  \frac{C_\Lambda}{|\vec{t}\,|^2}\|F\| .$$}
\end{proof}


\subsection{Asymptotic Form of the Induced Connection}\label{Sec:AsymCon}
In this subsection, we prove Theorem \ref{Thm:Asymp}.  

  We first define approximate solutions to the bow Dirac equation that are localized near $\lambda\in \Lambda$. We then take the $L^2$ orthogonal projection of these approximate solutions onto the kernel of $\D^\dagger_{t,b}$ and analyze the difference between the approximation and its projection. Next we estimate the $L^2$ inner products of these distinguished solutions. Finally, we use this information to construct a local unitary frame for the index bundle and compute the connection in this frame.

We write an element of the domain of $\D^\dagger_{t,b} $ as a triple $(f,w,n)$, where $f$ denotes the component in  $\oplus   {L^2}(I_\sigma)$, $w$ denotes the component in $\oplus_{\lambda\in \Lambda^0}W_\lambda  {\otimes\Fe^*_\lambda}$, and $n$ denotes the component in $\oplus_\sigma N_\sigma$. The latter component will play no role in the following computations. 

Let  $\eta\in C_c^\infty ((-2,2)),$ with $\eta = 1$ on $[-1,1]$ and $|d\eta|\leq 2.$ Fix $\lambda_i\in\Lambda^{+}$, and let $m_i=R(\lambda_i+)-R(\lambda_i-)>0$ be the rank discontinuity. Let $w_+$ be a unit spinor which is an   eigenvector of $\t:=t_1\e_1+t_2\e_2+t_3\e_3,$ 
with  $ \t\,w_{+}= -\ii |\vec{t}\,|w_{+}$. Hence, $w_+$ is a   {highest} weight vector with respect to the Cartan and Weyl chamber determined by $  {\ii}\t$. Similarly, let   $\xi_{i}$ be a unit covariant constant section of  $\FE\otimes\Fe^*$ which is a   {highest} weight vector with respect to  $  {\ii} t_j \rho_j$.  
Define for $\delta$   {small}:
\begin{align}f_{i,\delta}(s) := (s-\lambda_i)^{\frac{m_i-1}{2}}e^{-|\vec{t}\,|(s-\lambda_i)}\eta(\frac{s-\lambda_i}{\delta})w_{+}\otimes \xi_{i},\text{ for }s\geq \lambda_i\end{align}
and $f_{i,\delta}(s) = 0$ for $s\leq \lambda_i$. Similarly define  $f_{i,\delta}(s)$ for $\lambda_i\in \Lambda^-$. 

When $\lambda_c\in \Lambda^0$, let $P_+$ and $P_-$ denote the orthogonal projections onto the $\pm |\vec{t}|$ eigenspaces of  $\ii\t$. Let $w_c$ denote a unit vector in $  {W_\lambda\otimes}\Fe^*_\lambda$. 
Near $\lambda_c\in\Lambda^0$, consider the section
\begin{align}f_{c,\delta}=\begin{cases}
e^{-|\vec{t}\,|(s-\lambda_c)}\eta\big(\frac{s-\lambda_c}{\delta}\big)P_+Q_{\lambda_c}(w_c)\,& \text{ for } s>\lambda\\
e^{|\vec{t}\,|(s-\lambda_c)}\eta\big(\frac{s-\lambda_c}{\delta}\big)P_- Q_{\lambda_c}(w_c),& \text{ for } s<\lambda,
\end{cases}    
\end{align}
with $ Q_{\lambda_c}(w_c)(s)$ the covariant constant extension of $Q_{\lambda_c}(w_c)$ from $s=\lambda_c$.  Then $(f_{c,\delta},w_+,0)^T$ satisfies the discontinuity condition of Eq.~\eqref{sol2} at $\lambda_c$. 

 We next show that the triples $(f_{i,\delta},0,0)^T$ and  $(f_c,w_c,0)^T$  are approximate solutions to $\D^\dagger_{t,b}\left(\begin{smallmatrix}x\\ w\\ \nu\end{smallmatrix}\right)=0$. We will only compute in detail for $\lambda_i\in \Lambda^+$. The other cases are similar. Let $\Pi_{t,b}$ denote the orthogonal projection on $\mathrm{Ker}\,\D_{t,b}^\dagger.$
\begin{lemma}
\begin{align}\| (I-\Pi_{t,b})f_{i,\delta}\|_{L^2}^2 = O(|\vec{t}\,|^{-2}\| f_{i,\delta}\|_{L^2}^2).\end{align}
\end{lemma}
\begin{proof}
We have the following norm estimates.    
\begin{align*}
|\vec{t}\,|^{-{ m_i}}\int_{0}^{|\vec{t}\,|\delta} s^{ m_i-1 }e^{-2s}ds&\leq \|f_{i,\delta}\|_{L^2}^2  \leq   |\vec{t}\,|^{-{ m_i}}\int_{0}^{ 2|\vec{t}\,|\delta} s^{ m_i-1 }e^{-2s}ds.
\end{align*}
Acting on sections supported on $(\lambda_i,\lambda_i+2\delta)$, $\delta$ small, we write  $\D_{t,b}^\dagger = \nabla_{\frac{d}{ds}}+ \e_j(\frac{\rho_j}{2(s-\lambda_i)}-\ii t^j)+\mathfrak{B}$, where ${\mathfrak{B}}$ is a uniformly bounded endomorphism.  Since 
 $(\sum_{j=1}^3\e_j\otimes \frac{\rho_j}{2s} )\xi_{i} =\frac{(-m_i+1)}{2s}\xi_{i},$ we have 
\begin{multline}
\|\D_{t,b}^\dagger f_{i,\delta}\|_{L^2}^2 = \|\delta^{-1}(s-\lambda_i)^{\frac{m_i-1}{2}}e^{-|\vec{t}\,|(s-\lambda_i)}\eta'(\frac{s-\lambda_i}{\delta})w_{+}\otimes \xi_{i} + {\mathfrak{B}}f_{i,\delta}\|_{L^2}^2\\
\leq 8\delta^{-2}|\vec{t}\,|^{-m}\int_{|\vec{t}\,|\delta}^{ 2|\vec{t}\,|\delta} s^{ m_i-1 }e^{-2s}ds + 2\|{\mathfrak{B}}\|_{L^\infty}\|  f_{i,\delta}\|_{L^2}^2,
\end{multline}
  Hence
\begin{align}\frac{\|\D_{t,b}^\dagger f_{i,\delta}\|_{L^2}^2}{\| f_{i,\delta}\|_{L^2}^2}\leq 
\frac{ 8|\vec{t}\,|^{-m}\delta^{-2}\int_{|\vec{t}\,|\delta}^{ 2|\vec{t}\,|\delta} s^{ m_i-1 }e^{-2s}ds }{\| f_{i,\delta}\|_{L^2}^2}+
2\|{\mathfrak{B}}\|_{L^\infty}. \end{align}
For sufficiently large $|\vec{t}\,|$, so that $|\vec{t}\,|\delta$ is large, we have
\begin{align}\label{anotherest}\frac{\|\D_{t,b}^\dagger f_{i,\delta}\|_{L^2}^2}{\| f_{i,\delta}\|_{L^2}^2}\leq 
1+2\|{\mathfrak{B}}\|_{L^\infty}. \end{align}
Write 
$$f_{i,\delta} = \Pi_{t,b}f_{i,\delta}+(I-\Pi_{t,b})f_{i,\delta}.$$
Then  
\begin{align}1+2\|{\mathfrak{B}}\|_{L^\infty}\geq \frac{\|\D_{t,b}^\dagger(I-\Pi_{t,b})f_{i,\delta}\|_{L^2}^2}{\| f_{i,\delta}\|_{L^2}^2}\geq 
C_\Lambda^{-1}|\vec{t}\,|^2\frac{\| (I-\Pi_{t,b})f_{i,\delta}\|_{L^2}^2}{\| f_{i,\delta}\|_{L^2}^2}. \end{align}
{  The last inequality follows since $\D_{t,b}\D_{t,b}^\dagger$ and $\D_{t,b}^\dagger\D_{t,b}$ have the same nonzero spectrum. }Hence 
\begin{align}\label{compprojbnd}\| \Pi_{t,b}f_{i,\delta}\|_{L^2}^2\geq  \frac{(C_\Lambda^{-1}|\vec{t}\,|^2-1-2\|{\mathfrak{B}}\|_{L^\infty})}{1+2\|{\mathfrak{B}}\|_{L^\infty}}\| (I-\Pi_{t,b})f_{i,\delta}\|_{L^2}^2 . \end{align}
{ This implies the statement of the lemma. }
\end{proof}
Set $H_{ij}:= \langle \Pi_{t,b}f_{i,\delta},\Pi_{t,b}f_{j,\delta}\rangle_{L^2}.$ Let $(h^{ij})$ denote the inverse of the self-adjoint positive square root of $H$. We now choose our unitary frame, \{$\psi^j\}_j$ { of the index bundle.} 
$$\psi^j(t,b):= h^{jk}\Pi_{t,b}f_{k,\delta}.$$
Let $\beta(\lambda_i)$ denote the   $\lambda$-point or a $p$-point immediately following $\lambda_i$ along the circle. Choose $\delta$ small enough so that the support of $f_{i,\delta}$ is contained in $[\lambda_i,\frac{\lambda_i+\beta(\lambda_i)}{2}].$ This guarantees $\langle f_{i,\delta},f_{j,\delta}\rangle_{L^2} = 0$ for $i\not = j$. Hence for $i\not = j$,
\begin{align}H_{ij}
&= \langle  f_{i,\delta},f_{j,\delta}\rangle_{L^2}-\langle (I-\Pi_{t,b})f_{i,\delta},(I-\Pi_{t,b})f_{j,\delta}\rangle_{L^2}\nonumber\\
&= O(|\vec{t}\,|^{-2}\| f_{i,\delta}\|_{L^2}\| f_{j,\delta}\|_{L^2}).
\end{align}
Therefore, $$\psi^j = \frac{\Pi_{t,b}f_{j,\delta}}{\|\Pi_{t,b}f_{j,\delta}\|_{L^2}}+O(|\vec{t}\,|^{-2}).$$
We will also need to estimate derivatives of $H_{ij}$. 
\begin{lemma}
\begin{align}\label{Hbound}(\frac{\p}{\p\tau})^HH_{ij} =  O(|\vec{t}\,|^{-2} \|f_{i,\delta}\|_{L^2}\|f_{j,\delta}\|_{L^2}).
\end{align}
\end{lemma}
\begin{proof}
First, we need to estimate the derivatives of $\Pi_{t,b}$ and $f_{i,\delta}$.  Let $P$ be any curve of projection operators depending differentiably on a parameter $\tau$. Then differentiation of the expression $PP=P$ followed by left or right composition with $P$ yields
\begin{align}\label{PID}\frac{\p P}{\p\tau} = P\frac{\p P}{\p\tau}(I-P) + (I-P)\frac{\p P}{\p\tau}P.
\end{align} 
Write 
$$\Pi_{t,b} = I-\D_{t,b} \mathfrak{G}_{t,b}\D_{t,b}^\dagger.$$
Then 
$$\Pi_{t,b}(\frac{\p}{\p\tau})^H(\Pi_{t,b})(I-\Pi_{t,b}) = -\Pi_{t,b}(\frac{\p}{\p\tau})^H(\D_{t,b})\mathfrak{G}_{t,b}\D_{t,b}^\dagger .$$
By Lemma \ref{Greendecay} and \eqref{anotherest} respectively,  $ \mathfrak{G}_{t,b}  = O(|\vec{t}\,|^{-2}),$ and $\|\D_{t,b}^\dagger f_{i,\delta}\|_{L^2}=O(\| f_{i,\delta}\|_{L^2}).$ Therefore $ \mathfrak{G}_{t,b}\D_{t,b}^\dagger f_{i,\delta}   = O(|\vec{t}\,|^{-2} \|f_{i,\delta}\|_{L^2}).$ Since $(\frac{\p}{\p\tau})^H(\D_{t,b})$ is bounded,  we have 
\begin{align}\label{projdbound}\Pi_{t,b}(\frac{\p}{\p\tau})^H(\Pi_{t,b})(I-\Pi_{t,b}) f_{i,\delta}  = O(|\vec{t}\,|^{-2} \|f_{i,\delta}\|_{L^2}).
\end{align}

To compute  $(\frac{\p}{\p\tau})^H  f_{j,\delta}$,   observe that by construction 
\begin{align}\frac{\p}{\p\tau_l}  f_{j,\delta} = -i  f_{j,\delta}\end{align}  and $  f_{j,\delta}$ is independent of $\tau$.   Hence,   using the horizontal lift \eqref{Eq:Vell}, 
\begin{align}(\frac{\p}{\p\tau})^H  f_{j,\delta}=-\frac{\ii{ V_l(s)}}{V}  f_{j,\delta}.\end{align}

We compute 
\begin{align}(\frac{\p}{\p\tau})^HH_{ij}
&= \langle (\frac{\p}{\p\tau})^H( \Pi_{t,b})f_{i,\delta} ,\Pi_{t,b}f_{j,\delta}\rangle_{L^2} +\langle  \Pi_{t,b}( -\frac{\ii s}{V}) f_{i,\delta} ,\Pi_{t,b}f_{j,\delta}\rangle_{L^2}\nonumber\\
&+\langle   \Pi_{t,b}f_{i,\delta} , (\frac{\p}{\p\tau})^H(\Pi_{t,b})f_{j,\delta}\rangle_{L^2} +\langle \Pi_{t,b}f_{i,\delta} ,\Pi_{t,b}(-\frac{\ii s}{V}) f_{j,\delta}\rangle_{L^2}\nonumber\\
&=  \ii\langle  f_{i,\delta} ,[\Pi_{t,b},\frac{s}{V}] f_{j,\delta}\rangle_{L^2}+O(|\vec{t}\,|^{-2} \|f_{i,\delta}\|_{L^2}\|f_{j,\delta}\|_{L^2}).
\end{align}
The commutator with multiplication by $\frac{s}{V}$ is a derivation, and $[\D_{t,b}, \frac{s}{V}]$ and $[\D_{t,b},\frac{s}{V}]$ are bounded. Hence by the same argument as the derivation of \eqref{projdbound} (see also \eqref{142}), we have 
$$\langle  f_{i,\delta} ,[\Pi_{t,b},\frac{s}{V}] f_{j,\delta}\rangle_{L^2}=O(|\vec{t}\,|^{-2} \|f_{i,\delta}\|_{L^2}\|f_{j,\delta}\|_{L^2}).$$ Consequently 
\begin{align*}(\frac{\p}{\p\tau})^HH_{ij} =  O(|\vec{t}\,|^{-2} \|f_{i,\delta}\|_{L^2}\|f_{j,\delta}\|_{L^2}).
\end{align*}
\end{proof}
{  We now proceed to the proof of Theorem \ref{Thm:Asymp}}. We start analyzing inner products of our distinguished solutions $\Pi_{t,b}f_{i,\delta}$ that arise in the computation of the connection matrix. 
\
\begin{align}&\frac{\int_{\lambda_i}^{\beta(\lambda_i)}s|\Pi_{t,b}f_{i,\delta}|^2ds}{\int_{\lambda_i}^{\beta(\lambda_i)} |\Pi_{t,b}f_{i,\delta}|^2ds} = \frac{\int_{\lambda_i}^{\beta(\lambda_i)}s| f_{i,\delta}-(I-\Pi_{t,b})f_{i,\delta}|^2ds}{\int_{\lambda_i}^{\beta(\lambda_i)} |f_{i,\delta}-(I-\Pi_{t,b})f_{i,\delta}|^2ds}\nonumber\\
 \label{keyasymp}
 =&\frac{\int_{\lambda_i}^{\beta(\lambda_i)}s| f_{i,\delta} |^2ds}{\int_{\lambda_i}^{\beta(\lambda_i)} |f_{i,\delta} |^2ds}-2\mathrm{Re}\,\frac{\int_{\lambda_i}^{\beta(\lambda_i)}s \langle \Pi_{t,b}f_{i,\delta},(I-\Pi_{t,b})f_{i,\delta}\rangle ds}{\int_{\lambda_i}^{\beta(\lambda_i)} |f_{i,\delta} |^2ds} + O(|\vec{t}\,|^{-2})_{C^1}.\end{align}
Observe now that 
\begin{multline}\label{142}
\langle  \Pi_{t,b}f_{i,\delta},s(I-\Pi_{t,b})f_{i,\delta}\rangle_{L^2} =  
\langle  \Pi_{t,b}f_{i,\delta},s\D_{t,b}\mathfrak{G}_{t,b}\D_{t,b}^\dagger f_{i,\delta}\rangle_{L^2}\\
= \langle  \Pi_{t,b}f_{i,\delta},[s,\D_{t,b}]\mathfrak{G}_{t,b}\D_{t,b}^* f_{i,\delta}\rangle_{L^2}
= -\langle  \Pi_{t,b}f_{i,\delta},\mathfrak{G}_{t,b}\D_{t,b}^\dagger f_{i,\delta}\rangle_{L^2}.
\end{multline}
But, since $\|\D_{t,b}^\dagger f_{i,\delta}\|_{L^2}=\|\mathfrak{B}f_{i,\delta}\|_{L^2}\leq C\|  f_{i,\delta}\|_{L^2}$ and $\|\mathfrak{G}_{t,b}\|\leq C|\vec{t}\,|^{-2},$ this inner product is $O(|\vec{t}\,|^{-2})$, and \eqref{142} implies
  
 \begin{equation}\label{final}
\frac{\int_{\lambda_i}^{\beta(\lambda_i)}s|\Pi_{t,b}f_{i,\delta}|^2ds}{\int_{\lambda_i}^{\beta(\lambda_i)} |\Pi_{t,b}f_{i,\delta}|^2ds}= \frac{\int_{\lambda_i}^{\beta(\lambda_i)}s|f_{i,\delta}|^2ds}{\int_{\lambda_i}^{\beta(\lambda_i)} |f_{i,\delta}|^2ds}+O(|\vec{t}\,|^{-2})=\lambda_i+\frac{m_i}{2|\vec{t}\,|}+O(|\vec{t}\,|^{-2})_{C^1}.
 \end{equation}
 For $i\not = j$, a similar computation to \eqref{142} using the disjoint supports of $f_{i,\delta}$ and $f_{j,\delta}$ shows 
\begin{equation}\label{finalij}
\frac{\int_{\lambda_i}^{\beta(\lambda_i)}s\langle \Pi_{t,b}f_{i,\delta},\Pi_{t,b}f_{j,\delta}\rangle ds}{\|\Pi_{t,b}f_{i,\delta}\|_{L^2}\|\Pi_{t,b}f_{j,\delta}\|_{L^2}}= O(|\vec{t}\,|^{-2})_{C^1}.
 \end{equation}

The same computation gives the same asymptotics near $\lambda_i\in\Lambda^-$ and $\lambda_i\in \Lambda^0.$

Finally, we show that in the orthonormal frame $\{\psi^j\}_j$, $A(\frac{\p}{\p\tau})$ has the form claimed in Theorem \ref{Thm:Asymp}. We compute
\begin{align}\label{leibniz3}&A_{ij}(\frac{\p}{\p\tau})= \langle \psi^i,(\frac{\p}{\p\tau})^H\psi^j\rangle \nonumber\\
&=
\langle \psi^i, (\frac{\p}{\p\tau})^H(h^{jk}) f_{k,\delta} + h^{jk}[(\frac{\p}{\p\tau})^H,\Pi_{t,b}]f_{k,\delta}+h^{jk}(\frac{\p}{\p\tau})^Hf_{k,\delta}\rangle.
\end{align}

By \eqref{projdbound} and \eqref{Hbound}, 
$$\langle \psi^i, (\frac{\p}{\p\tau})^H(h^{jk}) f_{k,\delta} + h^{jk}[(\frac{\p}{\p\tau})^H,\Pi_{t,b}]f_{k,\delta} \rangle =O(|\vec{t}\,|^{-2}),$$
and we have 
\begin{align}\label{leibniz4} 
A_{ij}(\frac{\p}{\p\tau})& =
\langle \psi^i,  h^{jk}(\frac{\p}{\p\tau})^Hf_{k,\delta}\rangle+O(|\vec{t}\,|^{-2})\nonumber\\
& =
-\ii\delta_{ij}\frac{\int_{\lambda_i}^{\beta(\lambda_i)}s|f_{i,\delta}|^2ds}{V\int_{\lambda_i}^{\beta(\lambda_i)} |f_{i,\delta}|^2ds}+O(|\vec{t}\,|^{-2})-\frac{\ii (V_l(s)-s)}{V}\nonumber\\
& =
-\ii\delta_{ij}\frac{ \lambda_i+\frac{\hat m_i}{2|\vec{t}\,|}}{V} +O(|\vec{t}\,|^{-2}). 
\end{align}
 Computing additional derivatives, we see 
 \begin{align}\label{leibnizf} A_{ij}(\frac{\p}{\p\tau})& = 
-\ii\delta_{ij}\frac{ \lambda_i+\frac{\hat m_i}{2|\vec{t}\,|}}{V} +O(|\vec{t}\,|^{-2})_{C^1}. 
\end{align}
The computation of $A(\frac{\p}{\p t^j})$ is similar. This concludes the proof of Theorem \ref{Thm:Asymp}. 
(See also \cite[Sec.~10]{Cherkis:2009jm}.)

\section{Topology}
Our goal now is to express the Chern character values of $A$ in terms of the bow representation.  We achieve this by compactifying the multi-Taub-NUT space and extending the short exact sequence  \eqref{Eq:Exact1} of section \ref{BowIndex} to this compact space.

\subsection{Compactification \texorpdfstring{$\Mbar$}{Mbar} of the multi-Taub-NUT}

For the multi-Taub-NUT circle fibration
$\mathrm{TN}^\nu_k=\mathcal{M}\xrightarrow{\pi_k}\mathbb{R}^3$ we employ
the Hausel-Hunsicker-Mazzeo (HHM) compactification
$\Mbar=\mathcal{M}\bigsqcup S^2_\infty$ of \cite{HHM} which adds
a point at infinity for each direction in the base $\mathbb{R}^3.$
The result is a smooth compact manifold $\Mbar$.  We denote
the sphere at infinity by $\Cbar_\infty=S^2_\infty.$ 

Let us discuss $\Mbar$ in some more detail.  Consider a large closed ball $B_R=\left\{ \vec{t} : |\vec{t}|<R \right\}$ in $\mathbb{R}^3$ containing all NUTs, and denote its preimage by $\mathcal{M}_R=\pi^{-1}_k(B_R) \subset\mathcal{M}\subset\Mbar$. Its complement  in $\mathcal{M}$ is 
\begin{align}
\mathcal{M}\setminus\mathcal{M}_R\simeq\mathbb{R}_+\times S^3/\mathbb{Z}_k\simeq \mathrm{Tot}\, (\mathcal{H}^{-k})^0,
\end{align}
where the last expression stands for the total space of the degree $-k$ Hopf bundle $\mathcal{H}^{-k}\to S^2$ with its zero section deleted.  Let $\kappa=p\circ \pi_k:\mathcal{M}\setminus\mathcal{M}_R=(\mathrm{Tot}\, \mathcal{H}^{-k})^0\to S^2.$  The HHM compactification is defined by setting 
\begin{align}
\mathcal{M}_R^c=\Mbar\setminus\mathcal{M}_R=\mathrm{Tot}\, \mathcal{H}^{-k},
\end{align}
with the fiber radial coordinate $\rho=1/|\vec{t}|.$ In this description, $\Cbar_\infty$ is the zero section of this line bundle. The projection $\bar{\kappa}: 
\Mbar\setminus\mathcal{M}_R=\mathrm{Tot}\, \mathcal{H}^{-k}\to S^2$ allows us to pull back bundles from $S^2.$

Choosing a direction vector $\hat{n}$, for any NUT $\nu_\rho$, the preimage of the ray $ray_\rho:=\{\nu_\rho+\hat{n}u | u\in[0,+\infty) \}$ has a shape of a infinite cigar $C_\rho=\pi_k^{-1}(ray_\rho)$. (We use a generic direction $\hat{n}$, so that no ray  $ray_\rho$  passes through any other point $\nu_\sigma.$)   Note,
that for each cigar $C_\rho\subset\mathcal{M},$  its
compactification $\Cbar_\rho\subset\Mbar$ is diffeomorphic
to a two-sphere.  
The intersection numbers are 
\begin{align}
[\Cbar_\sigma]\cdot[\Cbar_\rho]&=-\delta_{\sigma\rho},&
[\Cbar_\sigma]\cdot[\Cbar_\infty]&=-1,&
[\Cbar_\infty]\cdot[\Cbar_\infty]&=-k.
\end{align}
The set $\left\{[\Cbar_\sigma\right]\}_{\sigma=1}^k$ forms a
basis of $H_2(\Mbar,\mathbb{Z})=\mathbb{Z}^k$  and
$[\Cbar_\infty]=\sum_\sigma[\Cbar_\sigma].$ 
Indeed, the Mayer-Vietoris sequence for  $\Mbar=\mathcal{M}\cup \mathcal{M}_R^c$  is
\[
\xymatrix{
	H_2(\mathcal{M}\cap \mathcal{M}_R^c) \ar[r]^-{(i_*,j_*)} & 	H_2(\mathcal{M})\oplus H_2(\mathcal{M}_R^c) \ar[r]^-{k_*-l_*} & 
	H_2(\Mbar) \ar `r[d]`[l]  `[ld]+<-5em,0em>_-{\partial_*} `[dl]
	[dl] \\
							    & H_1(\mathcal{M}\cap \mathcal{M}_R^c)\ar[r]^-{(i_*,j_*)} & H_1(\mathcal{M})\oplus H_1(\mathcal{M}_R^c).    
} 
\] 
Since $\mathcal{M}\cap\mathcal{M}_R^c$ deformation retracts to $S^3/\mathbb{Z}_k$ its second homology vanishes and its first homology is $\mathbb{Z}_k$.  
The neighborhood of infinity  $\mathcal{M}_R^c$ is contractible to $\overline{C}_\infty=S^2.$ And the space $\mathcal{M}$ deformation retracts to a bouquet of $k-1$ spheres $\{\Delta_\sigma:=\pi_k^*([\nu_{\sigma+1},\nu_\sigma]\}_{\sigma=1}^{k-1} $;   $[\Delta_\sigma]=[C_{\sigma+1}]-[C_\sigma].$ 
Thus the Mayer-Vietoris sequence becomes
\begin{align}
	0\to \mathbb{Z}^{k-1}\oplus \mathbb{Z}\to H_2(\Mbar)\to \mathbb{Z}_k\to 0
.\end{align}
In fact, the generator of $H_1(\mathcal{M}\cap\mathcal{M}_R^c)=\mathbb{Z}_k$ is the image $\partial_*\Cbar$ of any cigar $\Cbar_\sigma\subset \Mbar.$

If a line bundle $\Fe\rightarrow\mathcal{M}$  is presented near infinity (outside of $\mathcal{M}_R$) as a pullback of the degree $-m$ Hopf bundle: $\Fe|_{\mathcal{M}\setminus \mathcal{M}_R}=\pi_k^*p^* \mathcal{H}^{-m},$ then it naturally extends to $\ebar$ with
$\ebar|_{\Mbar\setminus\mathcal{M}_R}=\bar{\kappa}^* \mathcal{H}^{-m}$ 
 and $\ebar|_\mathcal{M}=\Fe$. 
Note, that since $\kappa^* \mathcal{H}^{-k}$ is trivial over $\mathcal{M}\setminus\mathcal{M}_R$ (since a pullback of a principal bundle to itself is trivial), $\Fe$ by itself determines $m$ up to a multiple of $k$, so defining $\ebar$ requires choosing a representative of $m \mod k$ (by specifying a particular bundle isomorphism between $\Fe_{\mathcal{M}\setminus \mathcal{M}_R}$ and $\pi_k^*p^* \mathcal{H}^{-m}$).  
 We return  to discuss this choice in Section~\ref{sec:distinguished}.  
 Our first task is to relate the Chern characters of $\ebar\rightarrow\Mbar$ with integrals of the Chern-Weil character forms over $\mathcal{M}.$

\subsection{Abelian Instanton}\label{Sec:AbIn}
We obtained the multi-Taub-NUT metric in the Gibbons-Hawking form $V d\vec{t}\,^{2}+\hat{\eta}^2/V$ in Section~\ref{Sec:SmallRep}. 
Abelian instantons $(\Fe,a)$ on this space enjoy a rather simple expression for their  connections, which allows us to explicitly evaluate their Chern-Weil forms and illustrate the relation between these values and the Chern characters of the {bundle $\bar{\Fe}$ defined over $\Mbar$.}  This example immediately specializes to the case of $\Fe_s$ and $\bar{\Fe}_s.$
In particular, we shall evaluate characteristic classes of $\Fe_s$ in  Lemma~\ref{lem:charFe}.  

The connection of any Abelian instanton on a Hermitian line bundle $\Fe\rightarrow\mathcal{M}$ has the form
\begin{align}\label{Eq:AbCon}
	d_a&=\pi_k^*d_\omega-\ii\frac{H}{V}\hat{\eta},&
	H&=\lambda+\sum_{\sigma=1}^k \frac{v_\sigma}{2 r_\sigma},&
	F_\omega=\ii *_3 dH,
\end{align}
where  $\omega$ is a connection 
on a line bundle $\Ff$ over the base $\mathbb{R}^3\setminus\{\nu_\sigma\}_\sigma$ with curvature $F_\omega$, {$\lambda\in\mathbb{R}$   
and $v_\sigma\in\mathbb{Z}.$   Indeed, this connection is anti-self-dual (with the orientation given by the volume form $V dt^1\wedge dt^2\wedge dt^3\wedge \hat{\eta})$, and, as we are about to verify, has $L^2$ curvature. 
Since for sufficiently large radius $R$  
$\frac{\ii}{2\pi}\int_{S^2_R} F_\omega=-\frac{1}{2\pi}\int_{S^2_R} *_3 dH=\sum_\sigma v_\sigma,$
one has $\Ff=\mathcal{H}^{\sum_\sigma v_\sigma}$ near infinity.   Therefore,
outside of a large ball  
$\Ff=p^*(\mathcal{H}^{\sum_\sigma v_\sigma})$, and  outside of $\mathcal{M}_R$, 
  $\Fe|_{\mathcal{M}\setminus\mathcal{M}_R}=\kappa^*(\mathcal{H}^{\sum_\sigma v_\sigma}).$ 
Integrals of the Chern-Weil character forms can be explicitly evaluated:
\begin{align}\label{Eq:ch1Exp}
ch_1^a[C_\rho]=c_1^a[C_\rho]=\frac{\ii}{2\pi}\int_{C_\rho} F_a=\frac{1}{2\pi}\int_{C_\rho} d\frac{H}{V}\wedge d\tau=\frac{\lambda}{\ell} - v_\rho,
\end{align}
\begin{multline}\label{Eq:ch2Exp}
ch_2^a[\mathcal{M}]=\frac12\left(\frac{\ii}{2\pi}\right)^2\int_\mathcal{M}F_a\wedge F_a\\
=\left(\frac{\ii}{2\pi}\right)^2\int_{\mathcal{M}}  (-\ii) \left(d\frac{H}{V} \right)\wedge (F_\omega-\frac{H}{V}F_\eta)
\wedge\hat{\eta}\\
=\frac{\ii}{2\pi} \lim_{R\rightarrow\infty}\left(\int_{S^2_R}-\sum_\sigma\int_{S^2_{\frac{1}{R}}(\nu_\sigma)}\right)\left(\frac{H}{V} F_\omega -\frac12\left(\frac{H}{V}\right)^2 F_\eta\right)\\
=-\frac12\left(\frac{\lambda}{\ell}\right)^2 k+\frac{\lambda}{\ell}\sum_\sigma v_\sigma
-\frac12\sum_\sigma v_\sigma^2.
\end{multline}

The connection \eqref{Eq:AbCon} defines a singular connection on $\ebar\to \Mbar$ (with singularity along $\overline{C}_\infty$, and we can use it to compute the characteristic classes of $\ebar$.  
Singular connections for {\em  holomorphic} bundles were used to express Chern classes in \cite{Atiyah57}, and relations similar to those we are about to obtain appear, for example, in  \cite[Prop.~13]{Atiyah57},  Biquard's \cite[Prop.7.2]{Biquard97}, and Esnault and Viehweg \cite[App. B, Cor.~B.3]{EV86}.  We emphasize, however, that the compactification $\Mbar$ is not compatible with any of the complex structures of the hyperk\"ahler manifold $\mathcal{M}$; so, we cannot apply these results immediately, since one cannot regard an instanton connection $A$ as a (singular) connection on a holomorphic bundle over $\Mbar.$  (These results could be applied, however, to a  different, holomorphic, compactification of $\mathcal{M}$ used in  \cite{Cherkis:2017pop}.)

The connection \eqref{Eq:AbCon} has logarithmic singularity along $\Cbar_\infty$, with singular part $-\ii \frac{\lambda}{\ell}\hat{\eta}$. One can use instead a family of regular connections
\begin{align}\label{Eq:regC}
	d_{\tilde{a}}&=d_a+\ii \Psi \frac{\lambda}{\ell}\hat{\eta}
,\end{align}
with $\Psi$ a bump function near infinity: $\Psi\equiv 1$ for $t\ge R+1$ and $\Psi\equiv 0$ for $t\le R.$

Noting that $-\frac{1}{2\pi}d(\Psi \hat{\eta})$ is the Thom class of $\Cbar_\infty,$ evaluating Chern classes becomes straightforward.
The integral of $-d(\Psi \hat{\eta})$ over a two-cycle equals the intersection number of that cycle with $\Cbar_\infty$ and integrating any two-form $\alpha$ against it equals  $\int_{\Cbar_\infty}\alpha.$  
Therefore, 
$\ch_1^{\ebar}[\Cbar_\rho]=\ch_1(d_a)[C_\rho]+\frac{\lambda}{\ell}[\Cbar_\rho]\cdot[\Cbar_\infty]=\frac{\lambda}{\ell}-v_\rho + \frac{\lambda}{\ell}(-1)=-v_\rho.$ In some detail, using the regularized connection \eqref{Eq:regC}:
\begin{multline}\label{Eq:ch_1e}
\ch_1^{\ebar}[\Cbar_\rho]=c_1^{\ebar}[\Cbar_\rho]=\frac{\ii}{2\pi}\int_{\Cbar_\rho}F_{\tilde{a}}\\
				 =\lim_{R\to \infty}\frac{\ii}{2\pi} \left\{ \int_{C_\rho^R} F_{a} 
				 + \int_{\Cbar_\rho\setminus  C_\rho^R} \left(F_\omega - \ii d \big((\frac{H}{V}  - \frac{\lambda}{\ell}\Psi)\hat{\eta}\big)  \right)\right\}\\
	=\frac{\ii}{2\pi}  \int_{C_\rho} F_{a}- \frac{\lambda}{\ell}=-v_\rho
.\end{multline}
Here $C_\rho^R=\pi_k^{-1}{(\{\nu_\sigma+u\hat{n}:0\leq u\leq R\})}$ is the truncated cigar, and its complement $\Cbar_\rho\setminus C_\rho^R$ forms a transverse disk to $\Cbar_\infty$ at its point of intersection with $\Cbar_\rho.$ We also use the crucial fact that $(\frac{\lambda}{\ell}\Psi-\frac{H}{V}  )\hat{\eta}$ is a one-form on $\mathcal{M}_R^c:=\Mbar\setminus \mathcal{M}_R$, allowing the use of the Stokes theorem in \eqref{Eq:ch_1e}.

Since $\ebar$ is a line bundle, its second Chern character is 
\begin{align}\label{Eq:ch_2e}
\ch_2^{\ebar}[\mathcal{M}]=\frac{1}{2}(c_1^{\ebar})^2[\mathcal{M}]=\frac{1}{2}\sum_{\sigma\rho} c_1^{\ebar}[\Cbar_\sigma] c_1^{\ebar}[\Cbar_\rho]\ \   [\Cbar_\sigma]\cdot [\Cbar_\rho]=-\frac{1}{2}\sum_\sigma v_\sigma^2. 
\end{align}
We  can also evaluate it explicitly, using the regularized connection. Observe that $F_{\tilde{a}}=F_a$ on $\mathcal{M}_R$ with $F_{\tilde{a}}=F_\omega-\ii d\left((\frac{H}{V}-\Psi \frac{\lambda}{\ell})\hat{\eta}\right)$, and split the integration over $\Mbar$ into an integral over $\mathcal{M}_R$ and an integral over $\Mbar\setminus \mathcal{M}_R$.  The latter is an integral of a closed form. 
\begin{align*}
	\ch_2^{\ebar}[\Mbar]=\frac{1}{2}\left( \frac{\mathrm{i}}{2\pi} \right)^2 \int_{\mathcal{M}_R} F_a\wedge F_a +\frac{1}{2}\left( \frac{\mathrm{i}}{2\pi} \right)^2 \int_{\mathcal{M}_R^c} F_\omega\wedge F_\omega\\ 
			    + \frac{\ii}{(2\pi)^2}\int_{\partial \mathcal{M}_R^c }(\frac{H}{V}-\Psi \frac{\lambda}{\ell}) F_\omega\wedge \hat{\eta}  \\
		    - \frac{1}{2}\left( \frac{\mathrm{i}}{2\pi} \right)^2 \int_{\partial \mathcal{M}_R^c }(\frac{H}{V}-\Psi \frac{\lambda}{\ell})\hat{\eta} d \left( (\frac{H}{V}- \Psi \frac{\lambda}{\ell})\hat{\eta} \right) 
.\end{align*}
As $R\to \infty$, the last term on the first line vanishes (since $|F_\omega|(\vec{t})<\frac{C}{|\vec{t}|^2}$).  Noting that $\partial M_R^c=-\partial M_R$, $\left(\frac{H}{V}-\Psi \frac{\lambda}{\ell}\right)|_{\partial \mathcal{M}_R}=\frac{H}{V}\vert_{|\vec{t}|=R}\xrightarrow{R\to \infty} \frac{\lambda}{\ell}$, we have
\begin{align*}
	\ch_2^{\ebar}[\Mbar]&=\frac{1}{2}\left( \frac{\mathrm{i}}{2\pi} \right)^2\int_{\mathcal{M}} F_a\wedge F_a\\
			    &-\frac{\ii}{2\pi}\frac{\lambda}{\ell} \int_{S^2_\infty}F_\omega 
	+\frac{1}{2}\left( \frac{\mathrm{i}}{2\pi} \right)^2 \left( \frac{\lambda}{\ell} \right)^2\lim_{R\to \infty}\int_{\partial \mathcal{M}_R} \hat{\eta}d \hat{\eta}\\
			    &=\frac{1}{2}\left( \frac{\mathrm{i}}{2\pi} \right)^2\int_{\mathcal{M}} F_a\wedge F_a
-\frac{\lambda}{\ell}\sum v_\sigma+\frac{1}{2}\left( \frac{\lambda}{\ell} \right)^2 k
.\end{align*}

Now we know the values of Chern characters of $\ebar$ and their relation to the integrals of the Chern-Weil forms over $\mathcal{M}$: 
\begin{align}
\frac{\ii}{2\pi}&\int_{C_\sigma} F_a=c_1^{\ebar}[\bar{C}_\sigma]-\frac{\lambda}{\ell}\underbrace{[\Cbar_\sigma]\cdot[\Cbar_\infty]}_{-1},
\end{align}
\begin{multline}
\frac12 \left(\frac{\ii}{2\pi}\right)^2 \int_{\mathcal{M}} F_a\wedge F_a=\frac12 (c_1^{\ebar})^2(\Mbar)+\frac{\lambda}{\ell}\mathrm{deg}\, \ebar|_{\Cbar_\infty}
+\frac12\left(\frac{\lambda}{\ell}\right)^2\mathrm{deg}\, N_{\Cbar_\infty}\\
\nonumber
=\frac12 c_1^{\ebar}(\bar{C}_\sigma) c_1^{\ebar}(\bar{C}_\rho) \underbrace{[\Cbar_\sigma]\cdot[\Cbar_\rho]}_{-\delta_{\sigma\rho}}
+\frac{\lambda}{\ell}\sum_\sigma v_\sigma
+\frac12\left(\frac{\lambda}{\ell}\right)^2 \underbrace{[\Cbar_\infty]\cdot[\Cbar_\infty]}_{-k}.
\end{multline}
This agrees of course with our direct evaluation \eqref{Eq:ch1Exp} and \eqref{Eq:ch2Exp}.

\begin{lemma}\label{lem:charFe}
		\begin{align}
			\ch_1^{\bar{\Fe}_s}[C_\rho]&=
			\begin{cases}
				1,& p_\rho<s\\
				0, & p_\rho>s
			\end{cases},\\
			\ch_2^{\bar{\Fe}_s}[\mathcal{M}]&=-\frac{1}{2}|\{p_\sigma | p_\sigma<s\}|  
		.\end{align}
\end{lemma}

\begin{proof}
The connection on $\Fe_s$ was computed in Eq.~\eqref{Eq:Conn-e}
\begin{align*}
	d_a= \pi_k^*d_{\eta_s}+\ii \frac{V_s}{V}\hat{\eta}
,\end{align*}
with $F_{\eta_s}=-\ii d \eta_s=-\ii *_3 dV_s.$  
Thus, comparing to \eqref{Eq:AbCon}, $\lambda=-s$, $v_\sigma=-1$ for  $p_\sigma<s$ and  $v_\sigma=0$ for  $p_\sigma>s$ and the lemma immediately follows from \eqref{Eq:ch_1e} and \eqref{Eq:ch_2e}:  $\ch_1^{\ebar}[\Cbar_\rho]=-v_\rho=1$ and $\ch_2^{\ebar}(\Mbar)=-\frac{1}{2}\sum_\sigma v_\sigma^2=-\frac12\sum_{p_\sigma<s}(-1)^2.$
\end{proof}

\subsection{Extending the Instanton Bundle to \texorpdfstring{$\Mbar$}{Mbar}}\label{Sec:AsympGen}
By Theorem~\ref{Thm:Asymp}, an instanton connection $d_A$ has asymptotic form
\begin{equation}\label{Eq:SplitA}
d_A= \oplus_{j=1}^n \left(\pi_k^* d_{\eta_j} - \ii
\frac{{\lambda_{j}+\frac{\hat{m}_{j}}{2|\vec{t}\,|}}}{V}\hat{\eta}
\right)+ O(|\vec{t}\,|^{-2}).
\end{equation}
Here 
$\hat{m}_{j}=R_{\lambda_j+}-R_{\lambda_j-} + \mathfrak{l}_{\lambda_j},$
 with $\mathfrak{l}_s:=\#\{p | 0\leq p <s\}$  denoting the number of $p$-points to the left of $s\in[0,\ell).$ 
Let $\mathfrak{r}_p:=\#\{\lambda | p<\lambda<\ell \}$  denote the number of $\lambda$-points to the right of $p,$  
and use $\Delta R_p:=R_{p+}-R_{p-}$ and $\Delta R_\lambda:=R_{\lambda+}-R_{\lambda-}$ for the discontinuity of the bow representation rank at the points $p$ and $\lambda$, respectively. 
 
The form \eqref{Eq:SplitA} of the connection implies that over the complement
of a compact set $\mathcal{M}_R$, $\mathcal{E}|_{\mathcal{M}\setminus\mathcal{M}_R}$ is isomorphic to the Whitney sum of pullbacks of $n$ Hopf bundles of degrees $\hat{m}_{j}$. Therefore, we extend it to the bundle $\Ebar\rightarrow\Mbar$ which coincides with $\mathcal{E}$ over $\mathcal{M}$ and with the pullback of the sum of Hopf bundles of degrees $\hat{m}_{j}$ over $\mathcal{M}_R^c=\Mbar\setminus\mathcal{M}_R:$
\begin{align}\label{Eq:EbComp}
\Ebar|_{\mathcal{M}}&=\mathcal{E},&
\Ebar|_{\mathcal{M}_R^c}&=\oplus_\lambda \bar{\kappa}^*(\mathcal{H}^{\hat{m}_\lambda}).
\end{align}
From this point of view \eqref{Eq:SplitA} defines a singular connection on $\Ebar$. The Chern character of $\Ebar$ is then expressed in terms of integrals of the Chern-Weil forms of this singular connection, understood as  distributions.   Equivalently, one can use exactly the same regularization as in the previous section as follows.

Let $\Psi$ be a bump function near infinity $(\Psi\equiv 1$ for  $t>R$ and $\Psi\equiv 0$ for  $t<R-1)$ and consider a regular connection
\begin{align*}
	d_{\tilde{A}}:= d_A+\ii\, \mathrm{diag}\,\left( \frac{\lambda_j}{\ell}\Psi \hat{\eta}\right)
.\end{align*}
Its curvature is $F_{\tilde{A}}=F_A+\ii\,\mathrm{diag}\,\left( \frac{\lambda_j}{\ell} (d\Psi\wedge \hat{\eta}+\Psi d\eta) \right)+\Psi O(t^{-2}).$ We note that $-\frac{1}{2\pi}d(\Psi \hat{\eta})$ is the Thom class of $\overline{C}_\infty$.  In particular, integrating any closed two-form against it equals  the integral of that form over $\overline{C}_\infty.$

Sending $R$ to infinity, we then have 
\begin{align}\label{Eq:ch1bar}
\frac{\ii}{2\pi}\int_{C_\sigma}\mathrm{tr} F_A&=
ch_1^{\Ebar}(\Cbar_\sigma) +
\sum_{\lambda} \frac{\lambda}{\ell},\\
\label{Eq:ch2bar}
\frac12\left(\frac{\ii}{2\pi}\right)^2\int_{\mathcal{M}}\mathrm{tr} F_A\wedge F_A&=
ch_2^{\Ebar}(\Mbar) 
+ 
\sum_\lambda\frac{\lambda}{\ell}\hat{m}_\lambda 
-\frac{k}{2}\sum_\lambda  \left(\frac{\lambda}{\ell}\right)^2.
\end{align}
To achieve our goal of evaluating the integrals of the Chern-Weil forms on the left-hand sides all that remains is to compute the Chern numbers of the bundle $\Ebar\rightarrow\Mbar.$

\subsection{Splitting Principle}
\label{sec:split}
Using the fact that the bundle $\mathcal{E}$ emerged via the short  exact sequence
\begin{align}\label{Eq:Exact1}
0\rightarrow\mathcal{E}\rightarrow X\oplus N \xrightarrow{\Delta} Y^\vee\rightarrow 0,
\end{align}
as the kernel of $\Delta$, we
 slightly modify and extend this sequence to the compactification $\Mbar$:
\begin{align}\label{Eq:SplitBar}
0\rightarrow\Ebar\rightarrow \bar{X}\oplus \bar{N}\oplus \overline{\W} \xrightarrow{\bar{\Delta}} \bar{Y}^\vee\rightarrow 0,
\end{align}
so that $\Ebar|_\mathcal{M}=\mathcal{E}.$ 
The first step in extending $X, N,$ and $Y$ to $\Mbar$
is  identifying them with certain bow sub-fibers of $\S\otimes \FE\otimes \Fe^*|_\mu$. Each of these bow fibers forms a bundle over $\mathcal{M}$ that  extends to   $\Mbar.$  
In particular, we extend the line bundles $\Fe_s$ to $\bar{\Fe}_s,$ using the pullback of the degree $\mathfrak{l}_s$
Hopf bundle.\footnote{
	Note, for $s\in I_\sigma=[a+,b-]$ we have  $\mathfrak{l}_s=\mathfrak{l}_{a+}=\mathfrak{l}_{b-}$ is the number of $p$-points to the left of $s$.  The only thing that determines the compactification of $\Fe_s$ is the number of $p$-points to the left of $s$. For any given bow interval $I_\sigma$ the bundles $\bar{\Fe}_s$ are isomorphic for all $s\in I_\sigma.$
}
Extension of $\S\otimes\FE$ compatible with our asymptotic decomposition \eqref{Eq:SplitA} is less straightforward.  
Since  $X$ and $Y$ are defined by sections of bundles that solve linear ODE, we can identify them with values these sections take in  specific fibers.  To simplify our exposition, we introduce an additional subdivision of the bow intervals, which we now describe.

\begin{figure}[htb!]
\centering
%
%
\begin{tikzpicture}
    \draw (0,0) ellipse (4 and 2);
    \node(p1) at ($(-155:4 and 2)$) [label={[shift={(0.4,-0.1)}]:$p_1$}] {\tikz\draw[red,fill=white] (0,0) circle (.5ex);}; 
           \node(p2) at ($(-55:4 and 2)$) {};
           \node(p3) at ($(75:4 and 2)$)  [label=below:$p_3$] {\tikz\draw[red,fill=white] (0,0) circle (.5ex);};
\coordinate (lam1) at ($(-185:4 and 2)$);
\coordinate (lam2) at ($(-130:4 and 2)$);
\coordinate (lam3) at ($(-85:4 and 2)$);
\coordinate (lam4) at ($(-5:4 and 2)$);
\coordinate (lam5) at ($(105:4 and 2)$);
           \node(l1) at (lam1) [label={[label distance=1pt]0:$\lambda_1$}] {\tiny$\bullet$};
           \node(l2) at (lam2) [label={[label distance=1]60:$\lambda_2$}] {\tiny$\bullet$};
           \node(l3) at (lam3)  [label={[label distance=0]90:$\lambda_3$}] {\tiny$\bullet$};
           \node(l4) at (lam4) [label={[label distance=0]175:$\lambda_4$}] {\tiny$\bullet$};
           \node(l5) at (lam5)  [label=-85:$\lambda_5$] {\tiny$\bullet$};
\coordinate (mu1) at ($(-179:4 and 2)$);
\coordinate (mu2) at ($(-142:4 and 2)$);
\coordinate (mu3) at ($(-125:4 and 2)$);
\coordinate (mu4) at ($(-90:4 and 2)$);
\coordinate (mu5) at ($(-70:4 and 2)$);
\coordinate (mu6) at ($(-13:4 and 2)$);
\coordinate (mu7) at ($(35:4 and 2)$);
\coordinate (mu8) at ($(90:4 and 2)$);
\coordinate (mu9) at ($(140:4 and 2)$);
\coordinate (c) at ($(-107:4 and 2)$);

           \node(m1) at (mu1) [label={[shift={(0.4,0)}]185:$\mu_1$},color=blue] {$\circledcirc$};
	   \node(m2) at (mu2) [label={[shift={(-0.1,0.2)}]-90:$\mu_2$},color=blue] {$\circledcirc$};
	   \node(m3) at (mu3) [label={[shift={(0,-0.9)}]:$\mu_3$},color=blue] {$\circledcirc$};
           \node(m4) at (mu4)  [label={[shift={(0,-0.9)}]:$\mu_4$},color=blue] {$\circledcirc$};
	   \node(m5) at (mu5)  [label={[shift={(0.3,-0.8)}]:$\mu_5$}, color=blue] {$\circledcirc$};
   \node(m6) at (mu6) [label={[shift={(0.4,-0.8)}]:$\mu_6$},color=blue] {$\circledcirc$}; 
   \node(m7) at (mu7) [label={[shift={(0.4,-0.5)}]:$\mu_7$},color=blue] {$\circledcirc$};
   \node(m8) at (mu8)  [label={[shift={(-0.3,-0.3)}]:$\mu_8$} ,color=blue] {$\circledcirc$};
   \node(m9) at (mu9)  [label={[shift={(0.3,-0.1)}]:$\mu_9$} ,color=blue] {$\circledcirc$};

           \node(c1) at (c)  [label=below:$c$, color=blue] {\tiny$\boxdot$};

    \draw[thick] ($(-185:5 and 3)$) node(o1){} arc (-185:-130:5 and 3);
    \draw[thick] ($(-130:6 and 4)$) node(o2){}  arc (-130:-85:6 and 4) node(o3){};
    \draw[thick] ($(-85:5 and 3)$)  arc (-85:-55:5 and 3);
    \draw[thick] ($(-55:5.5 and 3.5)$) node(o35){}  arc (-55:-5:5.5 and 3.5) node(o4){};
    \draw[thick] ($(-5:4.5 and 2.5)$)  arc (-5:105:4.5 and 2.5);
    \draw[thick] ($(105:4.5 and 2.5)$) node(o5){} arc (105:175:4.5 and 2.5);
   \draw (lam1) -- ($(-185:5 and 3)$);
   \draw (lam2)--(o2.center);
   \draw (lam3)--(o3.center);
   \draw (p2.center)--(o35.center);
   \draw (lam4)--(o4.center);
   \draw (lam5)--(o5.center);
       \draw[very thin,dotted] (0,0) ellipse (4.5 and 2.5);
       \draw[very thin,dotted] (o1.center) arc (-185:-5:5 and 3);
              \draw[very thin,dotted] (-85:5.5 and 3.5) arc (-85:-130:5.5 and 3.5);
   \node at (p2) [label=above:$p_2$] {\tikz\draw[red,fill=white] (0,0) circle (.5ex);};
\end{tikzpicture}%
\caption{New $\mu$-points, shown as \textcolor{blue}{$\circledcirc$}, and new $c$-point, shown as \textcolor{blue}{$\boxdot$}, result from a three step process. 1. Near the changing rank $\lambda$-points new $\mu$-points are $\mu_1,\mu_3,\mu_4,$ and $\mu_6$. 2. Points $\mu_2, \mu_5, \mu_7,\mu_8,$ and $\mu_9$ are introduced next. 3. A $c$ point is needed to separate $\mu_3$ and $\mu_4$.   }\label{fig:muPts}
\end{figure}
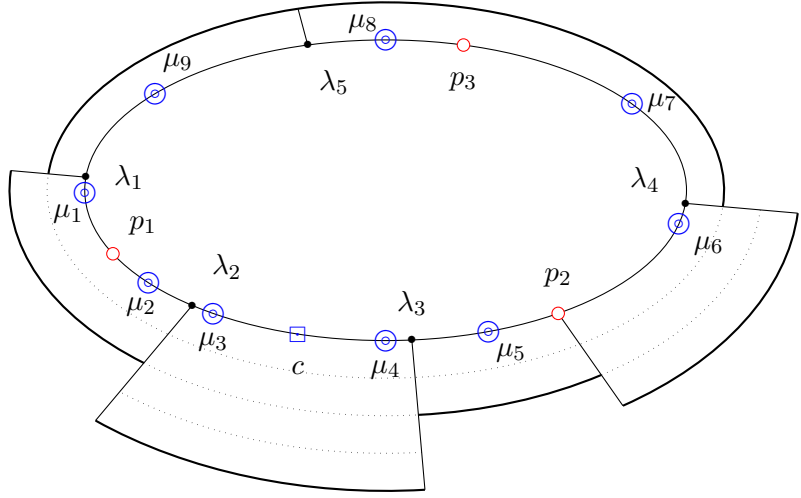
The bow intervals $I_\sigma$ emerged from dividing the circle at $p$-points.  Next, these were subdivided into subintervals $I_\sigma^\alpha$ by $\lambda$-points.  For  each sufficiently large $t:= |\vec{t}|$, we introduce a more refined subdivision as a result of a three-step process illustrated in Fig.~\ref{fig:muPts}:
\begin{enumerate}
\item
For each $\lambda\in\Lambda^+$ let $\mu_\lambda=\lambda+ t^{-\frac13}$ and for each $\lambda\in\Lambda^-$ let $\mu_\lambda=\lambda- t^{-\frac13}.$  This gives one nearby $\mu$-point on the larger rank side of the corresponding $\lambda$-point. 
\item
For any subinterval containing no  $\mu$-points, introduce one $\mu$-point in the middle of that subinterval.
\item
If a subinterval contains two $\mu$-points at his stage, it must have the form $[\lambda_1,\lambda_2]$ with the rank increasing at $\lambda_1$ and decreasing at $\lambda_2.$  For each such subinterval introduce a point $  c_{\lambda_1\lambda_2}$ (or simply $c$ when $\lambda_1,\lambda_2$  are clear) in its middle.  
\end{enumerate}
Now, to define $\tilde{X}:=\oplus_\mu X_\mu$, we cut the circle at $p-, \lambda-$, and $c-$points into {\em elementary intervals}, each containing a single $\mu$-point, and for each such interval, we define  $X_\mu$ to be the space of solutions to $(\frac{d}{ds}+\ii\T-\ii \t)x=0$ in its interior. 

As a result of our procedure, every $\lambda$-point is flanked by two $\mu$-points; denote them $\mu^\lambda_-$ and $\mu^\lambda_+$, so that $\mu^\lambda_-<\lambda<\mu^\lambda_+$. The same holds for every $c$-point: $\mu^c_-<c<\mu^c_+$.  
A  point $p-$ has a $\mu^p_-$ to its left and a point $p+$ has a $\mu^p_+$ to its right. Cutting the bow intervals $I_\sigma$ at $\mu$ points, we form $\tilde{Y}:=\oplus_\lambda \tilde{Y}_\lambda\oplus \oplus_c \tilde{Y}_c\oplus\oplus_p(\tilde{Y}_{p-}\oplus \tilde{Y}_{p+})$, where
\footnote{
	Clearly, these spaces are isomorphic to those of Eqs.~(\ref{Eq:Ylambdapole},\ref{Eq:YlambdaQ},\ref{Eq:Ypm},\ref{Eq:Ypp})
, except for $\tilde{Y}_\lambda$ with $\lambda\in\Lambda_0$. This space is now one dimension higher that the corresponding $Y_\lambda$.  We compensate for this by re-including the $\W=\W_\lambda$ space in the defining sequence, with $\W_\lambda=W_\lambda\otimes \Fe_\lambda^*$.
}:  
\begin{align*}
	\tilde{Y}_\lambda&=\left\{y\in \Gamma((\mu^\lambda_-,\mu^\lambda_+))\, \Big|\, (\frac{d}{ds}-\ii\T+\ii \t)y=0,  y^{term}(\lambda)=0\right\},
\lambda\in\Lambda\setminus\Lambda^0,\\
	\tilde{Y}_\lambda&=\left\{y\in \Gamma((\mu^\lambda_-,\mu^\lambda_+))\, \Big|\, (\frac{d}{ds}-\ii\T+\ii \t)y=0\right\},\  
\lambda\in\Lambda^0,\\
	\tilde{Y}_{p-}&=\left\{y\in \Gamma \left((\mu^p_-,p-]\right)\, \Big|\, (\frac{d}{ds}-\ii\T+\ii \t)y=0\right\},\\
	\tilde{Y}_{p+}&=\left\{y\in \Gamma \left([p+,\mu^p_+)\right) \, \Big|\,  (\frac{d}{ds}-\ii\T+\ii \t)y=0\right\}.
\end{align*}
In addition we introduce 
\begin{align}
	\tilde{Y}_{c}&=\left\{y\in \Gamma \left(\mu^c_-,\mu^c_+)\right) \, \Big|\,  (\frac{d}{ds}-\ii\T+\ii \t)y=0\right\}.
\end{align}
Most components of the map $\Delta: \tilde{X}\oplus N\oplus W\to \tilde{Y}^\vee$
are defined as earlier by pairing the elements of  $X_\mu$ space
with the two $\tilde{Y}$ spaces  associated with the intervals
containing $\mu$. The only modification involves $F_\lambda$
for $\lambda\in\Lambda_0$, which is now defined as  (compare to Eq.~\eqref{Eq:Flam})
\begin{align}
	F_\lambda(x)(y_\lambda)=\langle x,y_\lambda \rangle_{\mu^\lambda_-} - \left<x, y_\lambda \right>_{\mu^\lambda_+} + \left<Qw, y_\lambda \right>_\lambda
.\end{align}

As  $Y_\lambda$ was introduced to enforce the desired boundary
conditions on $x\in \tilde X$ at the $\lambda$-point,  the  purpose of $Y_c$ is to enforce continuity of $x$  at the $c$-point.   Indeed, if a subinterval $I_\sigma^\alpha=[\lambda_1,\lambda_2]$ is subdivided by a point $c$, then the space $X_\sigma^\alpha$ defined in Eq.~\eqref{Eq:Xsp}
is the kernel of the map $X_{[\lambda_1,c]}\oplus X_{[c,\lambda_2]}\xrightarrow{\Delta} \tilde{Y}_c.$  Thus, the same bundle $\mathcal{E}$ fits into the short exact sequence $0\rightarrow\mathcal{E}\rightarrow \tilde{X}\oplus N\oplus \W \xrightarrow{\Delta} \tilde{Y}^\vee\rightarrow 0,$ with the spaces $\tilde{X}$ and $\tilde{Y}$ being the direct sum of all the spaces in our new subdivision and $\W=\oplus_{\lambda\in \Lambda_0} W_\lambda\otimes \Fe_\lambda^*$.

Now, since each elementary interval has a unique $\mu$-point, we identify the corresponding $X_\mu$ space with the fiber at $\mu$:  $X_\mu\simeq(\S\otimes\FE\otimes\Fe^*)|_\mu.$  Using this identification we extend $ X_\mu$ to $\Mbar$ as the bundle $(\S\otimes\FE\otimes\bar{\Fe}^*)|_\mu$.  Here $\S$ and $\FE$ are trivial bundles\footnote{
For trivial bundles we use the same symbol to denote the bundle over $\mathcal{M}$ and its trivial extension to $\Mbar.$
} 
over $\mathcal{M}$ and are extended to trivial bundles over $\Mbar$, while $\Fe_\mu$ is extended to the previously defined $\overline{\Fe}_\mu$. 
(There is one significant exception to this rule, which we address in Section~\ref{sec:distinguished}. 

To compress notation, we denote the bundle $\S\otimes\FE\otimes\Fe^*$ by $\F$. 
The space $\tilde{Y}_c$ can be identified with the fiber $\F_c$, and, for $\lambda\in\Lambda^0$, the space $\tilde{Y}_\lambda$ can be identified with the fiber  $\F_\lambda.$ Thus each of these bundles extends easily to $\Mbar.$ 
The other $\tilde{Y}$ spaces  do not immediately admit such an identification.  Instead, on $\mathcal{M}_R^c$, we  identify each $Y_a$   with the direct sum of appropriate subspaces of $P_+ \F^*_{\mu^a_-}$ and $P_- \F^*_{\mu^a_+}$.  The isomorphism between $Y_a$ and these fibers allows us to extend  the bundle from $\mathcal{M}$ to a bundle over $\Mbar.$
Moreover, these extensions respect the map $\Delta$, so that it extends to all of $\Mbar.$ 

For any fiber, let $P_+$ be the orthogonal projection on the $t$ eigenspace of $\ii\t$ and let $P_-$ be such a projection on the $-t$ eigenspace of $\ii\t:$
\begin{align}
P_+&=\frac{t+\ii\t}{2t},&
P_-&=\frac{t-\ii\t}{2t}.
\end{align}
This gives a decomposition of $\S= {\kappa^*}\mathcal{H}\oplus  {\kappa^*}\mathcal{H}^{-1}$ into  {pullbacks of} Hopf line bundles with $P_+\S= {\kappa^*}\mathcal{H}$ and $P_-\S= {\kappa^*}\mathcal{H}^{-1}$.  

At a $\lambda$-point with changing rank, the fiber of $\S\otimes\FE$ decomposes into the space of continuous components and the two Casimir eigenspaces of dimensions $m-1$ and $m+1$, as discussed in Sec.~\ref{BowIndex}. 
As we already trivialized  our bundle over the subintervals, this decomposition extends to the whole interval. Let $E_m: \S \otimes \FE\to \S \otimes \FE$ denote the orthogonal projection onto the lowest weight vector of $\ii \t \otimes 1_{\FE} + \ii 1_\S\otimes t^j{\rho}_j.$

Each $\lambda$- or $c$-point is flanked by a pair of $\mu$-points: $\mu^\lambda_-<\lambda<\mu^\lambda_+$ or $\mu^c_-<c<\mu^c_+$.  For large $t$ the following linear maps are isomorphisms of linear spaces:  
for $\lambda\in\Lambda^+$
\begin{align}\label{Eq:IsoLamP}
	&\psi_\lambda:\tilde{Y}_\lambda \xrightarrow{\sim}P_+ \F_{\mu^\lambda_-}\oplus P_-(1-E_m)\F_{\mu^\lambda_+}\nonumber\\
&\psi_\lambda(y):= (P_+ y(\mu^\lambda_-), P_- (1-E_m)y(\mu^\lambda_+)),
\end{align}
for $\lambda\in\Lambda^-$
\begin{align}\label{Eq:IsoLamM}
	&\psi_\lambda:\tilde{Y}_\lambda \xrightarrow{\sim}P_+ (1-E_m) \F_{\mu^\lambda_-}\oplus P_-\F_{\mu^\lambda_+}\nonumber\\ 
&\psi_\lambda(y):= (P_+ (1-E_m)y(\mu^\lambda_-), P_- y(\mu^\lambda_+)),
 \end{align}
for $\lambda\in\Lambda^0$
\begin{align}\label{Eq:IsoLam0}
	&\psi_\lambda:\tilde{Y}_\lambda \xrightarrow{\sim} P_+   \F_{\mu^\lambda_-}\oplus P_-  \F_{\mu^\lambda_+}\nonumber\\
 &\psi_\lambda(y):= (P_+ y(\mu^\lambda_-), P_- y(\mu^\lambda_+)),
\end{align}
for a $c$-point
\begin{align}\label{Eq:IsoLamc}
&\psi_c:\tilde{Y}_c \xrightarrow{\sim} P_+  \F_{\mu^c_-}\oplus P_-\F_{\mu^c_+}\nonumber\\ 
 &\psi_c(y):= (P_+ y(\mu^c_-), P_- y(\mu^c_+)).
\end{align}
The fact that these are isomorphisms follows from  the inequality
\begin{align}\label{important}
|P_+ y(\mu^\lambda_+)| +|E_m P_- y(\mu^\lambda_+)|+|P_- y(\mu^\lambda_-)|   \leq C t^{-\frac{1}{3}}|\psi_\lambda(y)|,
\end{align}
proved as Corollary \ref{last}  in Section \ref{ylambda}.

For $\lambda\in\Lambda^-$, the above inequality holds with $+$ and $-$ subscripts interchanged.  For $y\in \tilde{Y}_\lambda$ with $\lambda\in\Lambda^0$ and for $y\in \tilde{Y}_c$, the above inequality holds with $E_m=0.$

For a $p\pm$-point, there is only one $\mu^\pm$ point on its interval, and we have the following analogs of \eqref{important}:
\begin{align}\label{pimportant}
&|P_- y(\mu^p_-)|  +|P_+ y(p-)|   \leq C t^{-\frac{1}{3}}(|P_+ y(\mu^p_-)|  +|P_- y(p-)|); \end{align} 
\begin{align}|P_+ y(\mu^p_+)|  +|P_- y(p+)|   \leq C t^{-\frac{1}{3}}(|P_- y(\mu^p_+)|  +|P_+ y(p+)|).
\end{align}  Thus, we have linear isomorphisms
\begin{align}
&\psi_{p-}:\tilde{Y}_{p-} \xrightarrow{\sim}P_+ \F^*_{\mu^p_-}\oplus 
 P_-\F^*_{p_-} ,\\
&\psi_{p+}:\tilde{Y}_{p+} \xrightarrow{\sim} 
 P_+\F^*_{p_+} 
\oplus P_- \F^*_{\mu^p_+}.
\end{align}

 {
Now, given these isomorphisms we can extend each bundle by extending each summand on the right-hand-side to a bundle over $\Mbar.$ E.g. $P_+\F_\mu=P_+\S\otimes \FE_\mu\otimes \Fe^*_\mu$ near infinity is isomorphic to $\mathcal{H}\otimes \FE_\mu\otimes \mathcal{H}^{-\mathfrak{l}_\mu}$, thus we define $\overline{P_+\F_\mu}\to\Mbar$ on $\mathcal{M}_R^c$ using the pullback of $\mathcal{H}^{1-\mathfrak{l}_\mu}\otimes \FE_\mu.$
}
  
Having defined the $\bar{Y}$ and $\bar{X}$ bundles, we now study the extension of $\Delta$ to  $\bar{\Delta}:\bar {X}\oplus \bar{ N}\oplus 
\overline{\mathcal{W}}\to \bar{Y}^\vee.$  Equivalently, we consider the extension of the maps $F_\lambda$ and $F_p$.   For concreteness, consider $\lambda\in\Lambda^+$.  Then for any $y_L\in P_+ \F_{\mu^\lambda_-}$ and $y_R\in P_-(1-E_m)\F_{\mu^\lambda_+}$, $\psi_\lambda^{-1}((y_L,y_R))$ defines a unique $y\in Y_\lambda$, satisfying $P_+y(\mu^\lambda_-)=y_L$ and $P_-(1-E_m)y(\mu^\lambda_+)=y_R$.  {The inequality \eqref{important}} ensures that  $y(\mu^\lambda_-)$   and $y(\mu^\lambda_+)$ are $O(t^{-\frac13})$ close to $y_L$ and $y_R$ in the respective fiber norms at $\mu^\lambda_\pm:$
\begin{align}
y(\mu^\lambda_-) &= y_L + O(t^{-\frac13} (|y_L|+|y_R|)),\\
y(\mu^\lambda_+) &= y_R + O(t^{-\frac13} (|y_L|+|y_R|)).
\end{align}
Thus the pairing \eqref{Eq:Flam}
\begin{align}
F_\lambda(x)(y)&=\langle x(\mu^\lambda_+), y(\mu^\lambda_+)\rangle - \langle x(\mu^\lambda_-),y(\mu^\lambda_-)\rangle \\
&=\langle x(\mu^\lambda_+), y_R\rangle - \langle x(\mu^\lambda_-),y_L\rangle \nonumber\\
&
+ O(t^{-\frac13} (|y_L|+|y_R|)(|x(\mu^\lambda_+) |+|x(\mu^\lambda_-)|)),
\end{align}
is also $t^{-\frac{1}{3}}$ close to
$\langle x(\mu^\lambda_+), y_R\rangle - \langle x(\mu^\lambda_-),y_L\rangle$.
Since, $P_+ y_L=y_L$ and $P_- y_R=y_R$, we conclude that $F_\lambda$  extends continuously to  $t=\infty$  where it vanishes on $P_- X_{\mu^\lambda_-}$ and on $P_+ X_{\mu^\lambda_+}$.  Therefore, $F_\lambda$ extends continuously to a map at infinity  with the form (for $\lambda\not\in\Lambda_0$)   
\begin{align}
	F_\lambda: P_+ \F_{\mu^\lambda_-}\oplus P_- \F_{\mu^\lambda_+}&\to 
	\tilde{Y}_\lambda^\vee,
\end{align}
and for $\lambda\in\Lambda_0$ by 
\begin{align}
	F_\lambda: P_+ \F_{\mu^\lambda_-}\oplus P_- \F_{\mu^\lambda_+}\oplus \W_\lambda&\to 
	\tilde{Y}_\lambda^\vee.
	\end{align}
	A similar discussion applies to other $\tilde{Y}$ spaces.  For $\tilde{Y}_p$ spaces, recall the domain of $\Delta$  is augmented by  the space $N_p=N_p^+\oplus N_p^-$,  defined in Eq.~\eqref{Eq:Ndef} with 
$N_p^+=\FE_{p_{\sigma}+}\otimes\Fe^*_{p_{\sigma}-}$ and 
$N_p^-=\FE_{p_{\sigma}-}\otimes\Fe^*_{p_{\sigma}+}$.  The maps $F_{p\pm}$ defined in 
\eqref{Eq:Fpm} and 
\eqref{Eq:Fpp} involve maps $\mathcal{A}$ and $\mathcal{B}$ satisfying for some $c_1,c_2>0$, 
\begin{align}
c_1t^\frac12 |n_+|\leq |\mathcal{A}(n_+,0)|&\leq c_2t^\frac12 |n_+|,& 
|\mathcal{A}(0,n_-)|&=O(|n_-|),\\ 
c_1t^\frac12|n_-|\leq |\mathcal{B}(0,n_-)|&\leq c_2t^\frac12|n_-|, & |\mathcal{B}(n_+,0)|&=O(|n_+|),
\end{align}
with $\mathcal{A}(n_+,0)\in P_+\F_{p+}$ and $\mathcal{B}(0,n_-)\in P_-\F_{p-}.$ Hence $F_p$ must be modified to extend continuously to $\infty$. The purpose of $F_p$ is to impose the boundary condition: $\exists \nu\in N_p$ such that 
$x(p_-) = \mathcal{B}(\nu)$ and $x(p_+) = \mathcal{A}(\nu).$  This is equivalent to the condition : $\exists \nu\in N_p$ such that 
$x(p_-) = t^{-\frac{1}{2}}\mathcal{B}(\nu)$ and $x(p_+) = t^{-\frac{1}{2}}\mathcal{A}(\nu).$ Hence, for $t$ large, we set 
$F_{p\pm}:X\oplus N\rightarrow Y_{p\pm}^\vee$  by 
\begin{align}\label{Eq:Fpmt}
F_{p-}(x,\nu)(y_p)&=
\phi(t)\ \langle \mathcal{B}\nu,y_p\rangle_{p-}-\langle x,y_p\rangle_{p-},\\
\label{Eq:Fppt}
F_{p+}(x,\nu)(y_p)&= 
 \langle x,y_p\rangle_{p+}-\phi(t)\langle \mathcal{A}\nu,y_p\rangle_{p+},
\end{align}  
where $\phi$ is continuous and satisfies $\phi(t)=1$ for $t$ small and $\phi(t)=t^{-\frac{1}{2}}$ for $t$ large. 
This ensures the continuous extension of $\bar{\Delta}$ at each $p$-point.

\subsubsection{Splitting at Infinity}
From the preceding discussion, we see that up to $O(t^{-\frac{1}{3}})$ error terms,  the map $\Delta$ decomposes into a sum  of maps which are   block-diagonal in the sense that (i) they are local, mapping $\F_\mu$ to $\F_\mu$, and (ii) they commute with $P_\pm$.   

\begin{align}
\Delta=\sum_\lambda\Delta_\lambda+\sum_p(\Delta_{p-}+\Delta_{p+})+O(t^{-\frac13}).
\end{align}
Thus, at infinity the bundle $\Ebar$ splits into a direct sum of bundles, as the exact sequence \eqref{Eq:Exact1} splits into the direct sum of exact sequences:
\begin{align}
	&0\rightarrow\Ebar_\lambda\rightarrow P_+\overline{\F}_{\mu^\lambda_-}\oplus P_-\overline{\F}_{\mu^\lambda_+}\xrightarrow{\Delta_\lambda} \overline{Y}_\lambda^\vee\rightarrow 0,\  \text{for } \lambda\in \Lambda_+\cup \Lambda_-\\
	&0\rightarrow\Ebar_\lambda\rightarrow P_+\overline{\F}_{\mu^\lambda_-}\oplus P_-\overline{\F}_{\mu^\lambda_+}\oplus W_\lambda\xrightarrow{\Delta_\lambda} \overline{Y}_\lambda^\vee\rightarrow 0,\ \text{for } \lambda\in\Lambda_{0}\\
	&0\rightarrow\Ebar_{p}\rightarrow P_- \overline{\F}_{\mu^p_+}\oplus P_+ \overline{\F}_{\mu^p_-}\oplus N_p \xrightarrow{\Delta_p} (Y_{p+}\oplus Y_{p-})^\vee\rightarrow 0.
	\end{align}
 
Recalling  the ranks of respective bundles,
\begin{align*}
&\dim P_\mp \F_{\mu^a_\pm}=R(\mu^a_\pm),&
&\dim Y_\lambda=R(\lambda+)+R(\lambda-)-1, \text{for } \lambda\not\in\Lambda_0\\
&\dim N_p^\pm=R(p\pm),&&\dim Y_{p\pm}=2R(p\pm),\\
&\dim W_\lambda=1,& &\dim Y_\lambda=2 R(\lambda), \quad 
 \text{for } \lambda\in\Lambda_0,
\end{align*}
we conclude that over $\mathcal{M}_R^c$, for sufficiently large $R$, $\Ebar_{p\pm}=0$, and each $\mathcal{E}_{\lambda}$ is a line bundle.  
So, the bundle $\Ebar$ over the sphere at infinity is a Whitney sum of line bundles $\Ebar=\oplus_\lambda\Ebar_\lambda,$ as expected.

We now compute the Chern classes of the resulting line bundles at infinity: 
\begin{multline*}
\ch_1(\Ebar_\lambda|_{S^2_\infty}) = \ch_1(\overline{P_+\F}_{\mu^\lambda_-}|_{S^2_\infty}) + \ch_1(\overline{P_-\F}_{\mu^\lambda_+}|_{S^2_\infty})-\ch_1(\overline{Y}^\vee_{\lambda}|_{S^2_\infty})\\
=R(\lambda-)(1-\mathfrak{l}_\lambda)+R(\lambda+)(-1-\mathfrak{l}_\lambda)-(R(\lambda-)+R(\lambda+)-1)(-\mathfrak{l}_\lambda)\\
  =-\Delta R_\lambda-\mathfrak{l}_\lambda=-\hat{m}_\lambda.
\end{multline*}
(E.g., $\ch\overline{P_+\F}_\mu|_{S^2_\infty}= \ch \mathcal{H}\otimes \FE_\mu \otimes \Fe_\mu^*=R(\mu)e^{c_1(\mathcal{H})} e^{-c_1(\ebar_\mu)}$.)  
Thus,   $\Ebar_\lambda=\kappa^*(\mathcal{H}^{-\hat{m}_\lambda}),$ and $\Ebar$ decomposes near $\bar{C}_\infty$ as $\oplus_\lambda\kappa^*(\mathcal{H}^{-\hat{m}_\lambda}),$
which is compatible with the splitting dictated by the asymptotic connection \eqref{Eq:SplitA}.

\subsection{The Distinguished Interval}
\label{sec:distinguished}
 The elementary interval $[a,b]$ containing the point $s=0$ requires special attention.  Here $a$ is a $\lambda$-, $p$-, or $c$-point with the largest value of $s$ attained by such points, and $b$ is a $\lambda$-, $p$-, or $c$-point with the smallest value of $s.$  We call $[a,b]$ the {\em distinguished interval} and denote the $\mu$-point within it by $\hat{\mu}.$

 The bundle $\Fe^*_s$ is presented at infinity as  $\kappa^*\mathcal{H}^{-\mathfrak{l}_s}$. For all other (not distinguished) elementary intervals, $\mathfrak{l}_s$ is constant, and therefore  the extended bundle $\bar{\Fe}^*_s$ is the same for all points in that interval.  For the distinguished interval, however,   $\mathfrak{l}_s=k$ for $s<0$, and  $\mathfrak{l}_s=0$ for $s\geq 0$, with discontinuity at $s=0.$  As a result, on the left $\bar{\Fe}^*_s = \bar{\Fe}^*_\ell=\kappa^*\mathcal{H}^{-k}$ and on the right $\bar{\Fe}^*_s=\bar{\Fe}^*_0$ is trivial.  This calls for special care in defining $\overline{X}_{\hat{\mu}}.$

For all other intervals, we identified $X_\mu$ with $\mathcal{F}_\mu$ in the interior and, after splitting near infinity $\F_\mu=P_-\F_\mu\oplus P_+\F_\mu$, we used the identifications 
$P_+\mathcal{F}_\mu=\FE_\mu\otimes\kappa^*(\mathcal{H}^{1-\mathfrak{l}_\mu})$ and 
$P_-\mathcal{F}_\mu=\FE_\mu\otimes\kappa^*(\mathcal{H}^{-1-\mathfrak{l}_\mu})$   to extend it to  $\Mbar.$ The above splitting and extension are compatible with the pairings $F$, which pairs the $P_+$ factor of $X_\mu$ with the $\tilde{Y}$ space on its right and the $P_-$ factor with the $\tilde{Y}$ space on its left. This ensures that the map $\Delta$ extends to infinity.

For the distinguished interval, we identify $X_{\hat{\mu}}$  with the fiber $\mathcal{F}_{\hat{\mu}}\simeq\S\otimes\FE_{\hat{\mu}}\otimes e_{0}$, which is a trivial bundle over the interior.  Once we split it near infinity into  $P_-\F_{\hat{\mu}}\oplus  P_+\F_{\hat{\mu}}$, however, we must identify  
the  summand $P_+\F_{\hat{\mu}}$ as $P_+ \S\otimes\FE_{\hat{\mu}}\otimes \Fe^*_0=\FE_{\hat{\mu}}\otimes\kappa^* \mathcal{H}$ (to pair with the space $\tilde{Y}_b$) and the  summand $P_-\F_{\hat{\mu}}$ as $P_- \S\otimes\FE_{\hat{\mu}}\otimes \Fe^*_\ell = \FE_{\hat{\mu}}\otimes\kappa^* \mathcal{H}^{-1-k}$ (to pair with the space $\tilde{Y}_b$). 
Thus, the bundle $\overline{X}_{\hat{\mu}}\to\Mbar$ is defined to be the trivial bundle $\S\otimes\FE_{\hat{\mu}}\otimes \Fe^*_0$ in the interior $ \mathcal{M}$  and 
$\FE_{\hat{\mu}}\otimes \overline{\mathcal{H}^{-1-k}\oplus  \mathcal{H} }$ over the neighborhood of infinity $ \mathcal{M}_R^c$.

Now that we understand how $\overline{X}_{\hat{\mu}}$ differs from all other $\overline{X}_\mu$ bundles, we can compute its Chern character.

\subsubsection{Compactifying the Trivial Bundle \texorpdfstring{$S$}{S}}
\label{sec:trivial}

Consider the following compactification of the trivial rank two bundle $S\to\mathcal{M}$. 
Using the decomposition of  $S$ on $\mathcal{M}\setminus\mathcal{M}_R$  into two eigen-line-bundles of $\t$, with    $P_\pm S|_{\mathcal{M}\setminus\mathcal{M}_R}=\kappa^*(\mathcal{H}^{\pm 1}),$  and define $\bar{S}\to\Mbar$ by
\begin{align}
\bar{S}|_\mathcal{M}&=S,&
\bar{S}|_{\Mbar\setminus\mathcal{M}_R}&=\kappa^*(\mathcal{H})\oplus\kappa^*(\mathcal{H}^{-1}).
\end{align}

Of course,  $\bar{S}\to\Mbar$ is  trivial, so $\ch(\bar{S})=2+0+0$. Nevertheless, we present it at infinity as a sum of these subbundles 
and compute its Chern classes now in order to set our notation.   

The bundle  $\bar{S}$ trivializes over $\mathcal{M}$, and we write $\bar{S}|_\mathcal{M} =V_{in}\times\mathcal{M}.$  Let $d$ denote the flat connection induced by the  above trivialization  on $\mathcal{M}$. In $\mathcal{M}\setminus \mathcal{M}_R,$ the connection 
$D:= \Pi_+d\Pi_+ + \Pi_-d\Pi_-=d+\Pi_+[d,\Pi_+] + \Pi_-[d,\Pi_-]$ 
preserves the splitting and extends to $\bar{S}|_{\mathcal{M}_R^c}$. 
With the bump function $\Psi$ introduced  in Section \ref{Sec:AbIn}, we set
$$\nabla:= (1-\Psi)d+\Psi D.$$
This is a connection on the whole of $\bar{S}.$ Then we have (using repeatedly $\Pi[d,\Pi]\Pi = 0$, for any smooth projection $\Pi$)
\begin{align}F_{\nabla}&= d\Psi\wedge (\Pi_+[d,\Pi_+] +  \Pi_-[d,\Pi_-])   +(2\Psi -\Psi^2 )[d,\Pi_+] \wedge [d,\Pi_+]   .
\end{align}
Since $F_\nabla$ manifestly has trace zero,  $c_1(\bar S) = 0$. 
Moreover, 
\begin{align}F_{\nabla}\wedge F_{\nabla}&= 2d\Psi\wedge (\Pi_+[d,\Pi_+] +  \Pi_-[d,\Pi_-]) \wedge (2\Psi -\Psi^2 )[d,\Pi_+]^{ \wedge 2} 
\end{align}
also has trace zero since it maps the image of $\Pi_+$ to the image of $\Pi_-$ and vice versa. 
  Thus,  $c_2(\bar{S})=0$ also, and we conclude that  $\ch(\bar{S})=\ch_0+\ch_1+\ch_2=2+ 0+ 0$ 
  and the bundle $\bar{S}\to\Mbar$ is trivial.

  \subsubsection{Chern Characters of \texorpdfstring{$X_{\hat{\mu}}$}{Xmu}}

Ignoring for now the trivial factor $\FE_{\hat{\mu}}$ of $\overline{X}_{\hat{\mu}}$, we identify
$\S\otimes\Fe^*_{\hat{\mu} }$ near infinity as  
$ (P_-S\otimes\Fe^*_\ell\oplus P_+S\otimes\Fe^*_0)|_{\mathcal{M}\setminus \mathcal{M}_R}
=\kappa^*\mathcal{H}^{-1}\otimes\Fe^*_\ell\oplus \kappa^*\mathcal{H}\otimes\Fe^*_0\equiv  \kappa^*\mathcal{H}^{-1-k}\oplus \kappa^*\mathcal{H}^{1}.$  Thus, we need to study the bundle $\overline{\kappa^*\mathcal{H}^{-1-k}\oplus\kappa^*\mathcal{H}^{1}}$, which in the interior is the trivial bundle $\S$ and over $\mathcal{M}_R^c$ is $\kappa^*\mathcal{H}^{-1-k}\oplus\kappa^*\mathcal{H}^{1}.$

We recall that we can identify a neighborhood $\mathcal{M}_R^c$ of $\overline{C}_\infty$ with the total space of a bundle: $\mathrm{Tot}\, \mathcal{H}^{-k}=\mathcal{M}_R^c.$  
The line bundle $\kappa^*\mathcal{H}^{-k}$ over $\mathcal{M}\setminus \mathcal{M}_R$ is the pullback of that bundle  to itself.  It is trivial and (in that trivialization) inherits the natural nontrivial connection $d+\ii\hat{\eta}$ from its base $\mathcal{M}_R^c$.   
We recall that  $d\hat{\eta}$ is $-k$ times the pullback of the volume form on $S^2$ and $\hat{\eta}(\partial/\partial_\tau)=1.$ 

The connection on $\kappa^*\mathcal{H}^{-1}$, induced from the trivial ambient bundle $S$ is $P_- dP_-.$ After tensoring with  $\kappa^*\mathcal{H}^{-k}$,  
the connection of $\mathcal{H}^{-1-k}$ is $P_- (d+\hat{\eta})P_-.$ And using the partition of unity, as in Subsection~\ref{sec:trivial}, to combine   this connection over $\mathcal{M}_R^c$ with the trivial connection over $\mathcal{M}$  we obtain the following connection on
$\overline{\kappa^*\mathcal{H}^{-1-k}\oplus\kappa^*\mathcal{H}^{1}}$ 
\begin{multline}
\nabla = (1-\Psi)d + \Psi(P_+dP_+  +P_-dP_-+ P_-\ii\hat\eta) \\
=d+\Psi(P_+[d,P_+]  +P_-[d,P_-]+ P_-\ii\hat\eta) ,
\end{multline}
with curvature
\begin{multline}
F_\nabla =  d\Psi\wedge (P_+[d,P_+]  +P_-[d,P_-]+ P_-\ii\hat\eta) \\  
+  (2\Psi-\Psi^2) [d,P_+]\wedge [d,P_+]    
 + \Psi   P_-\ii d\hat\eta + (\Psi^2-\Psi) [d,P_+]\wedge  \ii\hat\eta    .
\end{multline}
Then 
\begin{align}\tr F_\nabla &=\ii  d(\Psi\wedge \hat\eta ),
\end{align}
and  
\begin{align}
\tr F_\nabla\wedge F_\nabla  =  & \ii d\left((2\Psi^3-4\Psi^2)    \hat\eta \wedge  \tr  P_+ [d,P_+]   \wedge    [d,P_+]\right)\nonumber\\
&   -   d(\Psi^2  \hat\eta  \wedge    d\hat\eta ) . 
\end{align}
Hence 
\begin{align}
	\int_{\mathcal{M}}tr F_\nabla\wedge F_\nabla &= - \int_{S_\infty}(\hat\eta  \wedge    d\hat\eta + 2\ii\hat\eta \wedge  tr  P_+ [d,P_+]   \wedge    [d,P_+])\nonumber\\
&= -4\pi^2 (k+2). 
\end{align}
 
Hence we deduce  
\begin{align}
\ch(\overline{\mathcal{H}^{-1-k}\oplus\mathcal{H}^{1}} ) 
=2\oplus ch_1(\bar{\Fe}^*_l)\oplus \left(-1-\frac{k}{2}\right) \Omega,
\end{align}
where $\Omega$ is the generator of $H^4(\Mbar,\mathbb{Z})$; i.e. $\Omega(\Mbar)=1.$

We conclude that 
$
\ch_1(\bar{X}_{\hat{\mu}})[\bar{C}_\sigma]=-R_0, 
\ch_2(\bar{X}_{\hat{\mu}})[\Mbar]=- R_0\left(1+\frac{k}{2}\right), 
$ 
which we write as
\begin{align}
\ch(\bar{X}_{\hat{\mu}})= R_0\ch(\Fe^*_0)+R_0\ch(\Fe^*_\ell)- R_0\Omega.
\end{align}

\subsection{Evaluating the Chern Characters of \texorpdfstring{$\Ebar$}{Ebar}.}
\label{sec:chern_characters}

The identifications 
$\bar{N}_p=\FE_{p+}\otimes\bar{\Fe}^*_{p-}\oplus \FE_{p-}\otimes\bar{\Fe}^*_{p+}, \bar{\W}_\lambda=W_{\lambda}\otimes\bar{\Fe}_{\lambda}^*$, 
and $\bar{X}_\mu=\bar{\F_\mu},$ for $\mu\neq\hat{\mu}$  (of spaces of sections over the bow with particular fibers over the bow) are global.  The identification of the $\bar{Y}_\lambda$ spaces is more subtle.  For example, the isomorphisms \eqref{Eq:IsoLamP} and \eqref{Eq:IsoLamM} are valid only for sufficiently large $t.$ Note, that since \eqref{Eq:IsoLam0} and \eqref{Eq:IsoLamc} are valid in the limit of $\mu_{\pm}^\lambda\to\lambda$ and $\mu^c_\pm\to c$, we have
$\bar{Y}_\lambda=\bar{\F}_\lambda$ for $\lambda \in \Lambda_0,$  $\bar{Y}_c =\bar{\F}_c.$ Thus, we only need to find a global description of $\bar{Y}_\lambda$ for $\lambda \not\in \Lambda_0.$

First we modify the isomorphisms \eqref{Eq:IsoLamP} and \eqref{Eq:IsoLamM} near infinity.  By Corollary \ref{blasts},  
$\bar{Y}_\lambda \simeq Q_0\bar{\F}_\lambda\oplus V^{m-1}\otimes\Fe^*_\lambda.$ 
This  agrees with the global isomorphism $Q_0 Y_\lambda\simeq\F_\lambda\oplus V^{m-1}\otimes \Fe_\lambda^*$ that follows from Props.~\ref{impthm} and \ref{impthmc}.

We now compute the Chern character $\Ebar$ by applying the splitting principle to its defining short exact sequence \eqref{Eq:SplitBar}. 
\begin{align*}
	\ch \bar{X}_\mu&= \ch \bar{\S}\otimes\FE_\mu\otimes\bar{\Fe}^*_\mu=2 R_\mu \ch \bar{\Fe}^*_\mu, \quad \text{for \ } \mu\neq\hat{\mu},\\
	\ch \bar{X}_{\hat{\mu}}&= R_0 \ch \bar{\Fe}^*_0 + R_\ell \ch \bar{\Fe}^*_\ell - R_0 \Omega, \\
	\ch \bar{N}_p&= \ch \FE_{p+}\otimes\bar{\Fe}^*_{p-}\oplus\FE_{p-}\otimes\bar{\Fe}_{p+}^*=R_{p+} \ch\bar{\Fe}^*_{p-}+R_{p-}\otimes\ch\bar{\Fe}_{p+}^*,\\
	\ch \bar{\W}_\lambda&= \ch W_\lambda\otimes\bar{\Fe}^*_{\lambda}
,\end{align*}
for $\lambda \not\in \Lambda_0$
\begin{align*}
	\ch \bar{Y}_{\lambda}&= \ch \bar{\S}\otimes\FE^{cont}\otimes\bar{\Fe}_\lambda^* +\ch V^{m-1}\otimes\bar{\Fe}_\lambda^* = (R_{\lambda-}+R_{\lambda+}-1) \ch \bar{\Fe}_\lambda^*
,\end{align*}
and for $\lambda \in \Lambda_0$
\begin{align*}
	\ch \bar{Y}_\lambda - \ch\bar{\W}_\lambda&= \ch \bar{\S}\otimes\FE_\lambda\otimes\bar{\Fe}^*_\lambda - \ch W_\lambda\otimes\bar{\Fe}^*_{\lambda} = (2 R_\lambda-1) \ch \bar{\Fe}^*_\lambda 
.\end{align*}
Assembling these we have
\begin{align*}
\ch \Ebar &= \ch \bar{X} +\ch \bar{N}  - (\ch  \bar{Y}^{\vee}-\ch\bar{\W}) \\ 
	&=\sum_\mu 2 R_\mu \ch \bar{\Fe}^*_{\mu} +
	\sum_{p}(R_{p+}\ch\bar{\Fe}^*_{p-}+R_{p-}\ch\bar{\Fe}^*_{p+})
- \sum_c 2 R_c\ch\bar{\Fe}^*_c \\
&-\sum_{\lambda}(R_{\lambda+}+R_{\lambda-}-1)\ch\bar{\Fe}^*_\lambda 
-\sum_{p} (2R_{p+}\ch \bar{\Fe}^*_{p+} +2R_{p-}\ch\bar{\Fe}^*_{p-})\\
&- R_0\Omega
,\end{align*}
and, after some rearrangement 
\begin{multline}\label{Eq:Che}
\ch \Ebar =\sum_{\lambda} \ch\bar{\Fe}^*_\lambda 
- \sum_{p}(R_{p+}-R_{p-})(\ch\bar{\Fe}^*_{p+}-\ch\bar{\Fe}^*_{p-}) 
-R_0\Omega \\
+ \sum_{\mu}2R_\mu \ch\bar{\Fe}^*_\mu - \sum_\lambda(R_{\lambda+}+R_{\lambda-}) \ch \bar{\Fe}^*_\lambda
- \sum_{c}2R_c \ch\bar{\Fe}^*_c  
- \sum_{p\pm}R_{p\pm} \ch\bar{\Fe}^*_{p\pm}
.\end{multline}
Since all bundles $\bar{\Fe}_s$  are isomorphic within any given interval 
$I_\sigma=[p_\sigma+,p_{\sigma+1}-]$ (and within the intervals $[0,p_1-]$ and $[p_k+,\ell]$), traversing this interval left to right, one concludes that the consecutive terms of the last line of \eqref{Eq:Che} one encounters cancel each other.  Thus, the last line of \eqref{Eq:Che}  in fact vanishes.  

The Chern data of the line bundles $\bar{\Fe}^*_s$ follow from Lemma~\ref{lem:charFe}:
\begin{align}
	\ch_1^{\bar{\Fe}_s^*}[\Cbar_\sigma]&=-c_1({\bar{\Fe}_s})[\Cbar_\sigma]=
	\begin{cases} -1,& \text{for\ } s > p_\sigma\\ 
	0,& \text{for\ } s<p_\sigma \end{cases},\\
\ch_2^{\bar{\Fe}_s^*}[\Mbar]&=\frac12 c_1^2(\bar{\Fe}_s)[\Mbar]=\sum_{p_\sigma <s}\frac12 (-1) = -\frac{\mathfrak{l}_s}{2}.
\end{align}

Evaluating $\ch_1^{\Ebar}[\bar{C}_\sigma]$ using \eqref{Eq:Che} in the first term (the $\lambda$ sum) only terms with $\lambda>p_\sigma$ contribute, each contributing $-1$, and in the second term (the $p$ sum)  only $p=p_\sigma$ term does not vanish, contributing $\Delta R_{p_\sigma}$.  The last line vanishes. 

Similarly, evaluating $\ch_2^{\Ebar}[\Mbar]$ with \eqref{Eq:Che}, same $\lambda$-terms  contribute, each contributing $-\frac12 \mathfrak{l}_\lambda$, with each $p$-term contributing $-\frac12 \Delta R_p$, the last  line cancels.  All together:
\begin{align}
ch_1^{\Ebar}[\bar{C}_\sigma]&=-r_{p_\sigma}+\Delta R_{p_\sigma},\\
ch_2^{\Ebar}[\Mbar]&=\sum_\lambda(-\frac12 \mathfrak{l}_\lambda) 
- \sum_p\Delta R_p(-\frac12) -  R_0 
=-\frac12\sum_\lambda\hat{m}_\lambda-  R_0,
\end{align}
where we used the fact that the total rank change around the bow is zero: $\sum_p \Delta R_p+\sum_\lambda \Delta R_\lambda=0.$

Consequently, using \eqref{Eq:ch1bar} and \eqref{Eq:ch2bar} we have
\begin{align}\label{Eq:Ch1}
\frac{\ii}{2\pi}\int_{C_\sigma}\mathrm{tr} F&=
\Delta R_{p_\sigma}-r_{p_\sigma} + 
\sum_{\lambda} \frac{\lambda}{\ell},\\
\label{Eq:Ch2}
\frac12\left(\frac{\ii}{2\pi}\right)^2\int_{\mathcal{M}}\mathrm{tr} F\wedge F&=
-\frac12\sum_\lambda\hat{m}_\lambda-  R_0 
+ 
\sum_\lambda\frac{\lambda}{\ell}\hat{m}_\lambda 
-\frac{k}{2}\sum_\lambda  \left(\frac{\lambda}{\ell}\right)^2.
\end{align}

\section{\texorpdfstring{$Y_\lambda$}{Ylambda} Estimates}\label{ylambda}
In this section we prove a quantitative one dimensional analog of the stable manifold theorem  for elements of $Y$ stated in Proposition \ref{prop1}. We recall that $P_\pm$ denotes orthogonal projection onto the $\pm t$  eigenspaces of $i\e_jt^j.$
Let $[\mu^\lambda_-,\mu^\lambda_+]\ni\lambda$ be an interval containing no other $\lambda$- or $p$-points and let $m=|\Delta_\lambda|$ be the dimension of the space of terminating components at $\lambda.$ 
Any section $y\in Y_\lambda,$ satisfies the equation 
$\frac{dy}{ds} = i{\mathcal T}y.$
For $\lambda\in \Lambda^+$ (respectively $\lambda\in \Lambda^-$), we further decompose $i{\mathcal T}$ on the higher rank side $(\lambda,\mu^\lambda_+]$ (respectively on $[\mu^\lambda_-,\lambda)$) as 
\begin{align}i{\mathcal T} =  \frac{(m+1)Q_+-(m-1) Q_-}{2(s-\lambda)}+t(P_--P_+)+B,
\end{align}
where $B$ is bounded, with bound independent of $s,t$. Here, as on page~\pageref{zerodirac},  $Q_-$ denotes the orthogonal projection onto $V^{m +1}$, and $Q_+$ is orthogonal projection onto  $V^{m -1}$. (We recall that these subspaces, defined in Subsection \ref{BowIndex}, are the $su(2)$ irreducible summands of $V^2\otimes V^{m}$). We also set $Q_0:=I-Q_+-Q_-$ be the projection on the space of continuous components, $Q_1=I-Q_0$, and let $E_m$ denote projection onto the lowest weight vector of $\ii t^j(\e_j\otimes 1_\FE + 1_\S\otimes \bm{\rho}_j)$ in $V^{m +1}$.

\begin{proposition}\label{prop1}Let 
$y\in Y_\lambda,$ with $\lambda\in \Lambda^+ $. Assume $\lambda-\mu^\lambda_-\geq \mu^\lambda_+-\lambda\geq \frac{\ln(t)}{t}.$ Then there exists $C>0$ so that for $t$ large and $(\mu^\lambda_+-\lambda)$ small, 
\begin{multline}\label{propeq1}
|P_+y(\mu^\lambda_+)|+|E_my(\mu^\lambda_+)|+|P_-y(\mu^\lambda_-)|\\ 
\leq C(\mu^\lambda_+-\lambda)(|P_-(I-E_m)y(\mu^\lambda_+)|+|P_+y(\mu^\lambda_-)|).
\end{multline}
For $y\in Y_\lambda,$ with $\lambda\in \Lambda^- $, assume $\mu^\lambda_+-\lambda\geq\lambda-\mu^\lambda_-\geq  t^{-\frac{1}{2}}.$ Then there exists $C>0$ so that for $t$ large, 
\begin{multline}|P_-y(\mu^\lambda_-)|+|E_my(\mu^\lambda_-)|+|P_+y(\mu^\lambda_+)|\\ 
\leq C(\lambda-\mu^\lambda_-)(|P_+(I-E_m)y(\mu^\lambda_-)|+|P_-y(\mu^\lambda_+)|).
\end{multline}
For $y\in Y_\lambda$, $\lambda\in \Lambda^0$ or $y\in Y_c$, the analog of the above inequalities holds with $E_m$ set to $0$. For $y\in Y_{p\mp}$, the analogous inequality holds with $E_m$ set to zero, and $\mu^\lambda_\pm$ replaced by $p\mp$ (and $\mu^\lambda_\mp$ replaced by $\mu^p_\mp$).
\end{proposition}

\begin{corollary}\label{last}
For $t$ large, $\mu^\lambda_+-\mu^\lambda_-$ small, and $\mu^\lambda_+$, $\mu^\lambda_-$, constrained as in Proposition \ref{prop1}, the bundle 
$Y_\lambda$, $\lambda\in \Lambda^+$ is isomorphic to $(P_+S\otimes {\mathcal{E}}\otimes e)_{\mu^\lambda_-}\oplus (I-E_m)(P_-S\otimes {\mathcal{E}}\otimes e)_{\mu^\lambda_+} .$ For $\lambda\in \Lambda^-$, $Y_\lambda$  is isomorphic to $ (I-E_m)(P_+S\otimes {\mathcal{E}}\otimes e)_{\mu^\lambda_-}\oplus (P_-S\otimes {\mathcal{E}}\otimes e)_{\mu^\lambda_+} .$ Analogous isomorphisms hold for $Y_\lambda$ with $\lambda\in \Lambda^0$, $Y_c$, and $Y_p$. 
\end{corollary} 
\begin{proof}Let $\lambda\in \Lambda^+$. By \eqref{propeq1}, the natural evaluation map from 
$Y_\lambda$ to   $(P_+S\otimes {\mathcal{E}}\otimes e)_{\mu^\lambda_-}\oplus (I-E_m)(P_-S\otimes {\mathcal{E}}\otimes e)_{\mu^\lambda_+}$ is injective. It is surjective by a dimension count. The other cases follow similarly. 
\end{proof}

\begin{corollary}\label{blasts}
For $t$ large, and $\mu^\lambda_+$, $\mu^\lambda_-$, constrained as in Proposition \ref{prop1}, the bundle 
$Y_\lambda$, $\lambda\in \Lambda^+$ is isomorphic to $(P_+S\otimes {\mathcal{E}}\otimes e)_{\mu^\lambda_-}\oplus Q_+( S\otimes {\mathcal{E}}\otimes e)_{\mu^\lambda_+} .$ For $\lambda\in \Lambda^-$, $Y_\lambda$  is isomorphic to $ Q_+(S\otimes {\mathcal{E}}\otimes e)_{\mu^\lambda_-}\oplus (P_-S\otimes {\mathcal{E}}\otimes e)_{\mu^\lambda_+} .$
\end{corollary} 
\begin{proof}Let
 $\lambda\in \Lambda^+$. The vector spaces $Q_+( S\otimes {\mathcal{E}}\otimes e)_{\mu^\lambda_+}$ and $(I-E_m)(P_- S\otimes {\mathcal{E}}\otimes e)_{\mu^\lambda_+}$ both have dimension $m-1$. Hence the projection $Q_+: (I-E_m)(P_- S\otimes {\mathcal{E}}\otimes e)_{\mu^\lambda_+}\to Q_+(  S\otimes {\mathcal{E}}\otimes e)_{\mu^\lambda_+}$ is an isomorphism if injective.   The kernel of the map is contained  in the image of  $P_-$ intersected with the image of $Q_-$, but that space  is precisely the kernel of $(I-E_m)$. 
\end{proof}

The rest of this section is 
 devoted to a proof of Proposition \ref{prop1}. We will treat the case $\lambda\in \Lambda^+$. The proof of the other cases is similar. We will also restrict attention to  $m \geq 2$, as the $m =1$ case is significantly simpler. 

\subsection{Integral estimates}
 Given a $C^1$ section $\psi$ of $\mathrm{End}\,\F$ and a section $y$ of $\F$ satisfying $Dy:=\frac{d}{ds}y-i\mathcal{T}y = 0$  on an interval $[a,b]$, we have
\begin{align}\label{this}
&\|[D,\psi]y\|^2_{L^2([a,b])}=\|D(\psi y)\|^2_{L^2([a,b])}\nonumber\\
&
=\langle  \psi y, D^*D(\psi y) \rangle_{L^2([a,b])} +\langle \psi y, D\psi y\rangle|_a^b\nonumber\\
&=\langle \psi y, \nabla^*\nabla(\psi y) \rangle_{L^2([a,b])} +\langle \psi y, D\psi y\rangle|_a^b\nonumber\\
&=\|\tilde \nabla(\psi y)\|^2_{L^2([a,b])}-\langle \psi y, \frac{d}{ds}(\psi y)\rangle|_a^b+\langle \psi y, D\psi y\rangle|_a^b\nonumber\\
&=\|\tilde \nabla(\psi y)\|^2_{L^2([a,b])} - \langle \psi y, i\mathcal{T}\psi y\rangle|_a^b,
\end{align}
where we have set $\tilde \nabla z: =(\frac{dz}{ds}, i(T_1-t_1)z,i(T_2-t_2)z,i(T_3-t_3)z) .$ 

We continue to assume $y\in Y_\lambda$, $\lambda\in \Lambda^+$. 
Then, for $[a,b]=[\lambda,\mu^\lambda_+]$, rewrite \eqref{this} as 
\begin{align}\label{azalea}&
\int_{\lambda}^{\mu^\lambda_+}(|\frac{d (\psi y)}{ds}|^2+\sum_k|(T^k-it^k) \psi y|^2-|[D,\psi] y|^2)ds \nonumber\\
&
=  \langle  \psi y ,i\mathcal{T}\psi y\rangle(\mu^\lambda_+) - \langle   Q_0\psi y ,i\mathcal{T}Q_0\psi y\rangle(\lambda).   \end{align}
By the fundamental theorem of calculus, for $s>\lambda$, 
\begin{align}\label{FTC}\frac{d}{ds}\langle  \psi y ,i\mathcal{T}\psi y\rangle = |\frac{d (\psi y)}{ds}|^2+\sum_k|(T^k-it^k) \psi y|^2-|[D,\psi] y|^2.
\end{align}
Hence $ \langle   \psi y ,i\mathcal{T}\psi y\rangle(s)$ is increasing  if $|\frac{d (\psi y)}{ds}|^2+\sum_k|(T^k-it^k) \psi y|^2-|[D,\psi] y|^2\geq 0$. 
In particular, 
in any interval $I$ on which $|\frac{d   y }{ds}|^2+\sum_k|(T^k-it^k)   y|^2 \geq \frac{1}{2}\alpha_I^2 t^2|y|^2,$ 
we have 
 \begin{align}\label{gnu1}\frac{d}{ds}[\frac{1}{2}|y|^2+\frac{1}{\alpha_I t} \langle   y ,i\mathcal{T} y\rangle]\geq 
t\alpha_I[ \frac{1}{2}|y|^2+\frac{1}{\alpha_I t} \langle   y , i\mathcal{T}y\rangle], 
\end{align} 
and 
\begin{align}\label{gnu2}\frac{d}{ds}[\frac{1}{2}|y|^2-\frac{1}{\alpha_I t} \langle   y ,i\mathcal{T} y\rangle]\leq 
-t\alpha_I[ \frac{1}{2}|y|^2-\frac{1}{\alpha_I t} \langle   y ,i\mathcal{T} y\rangle]. 
\end{align} 
Let 
\begin{align*}
	\kappa_{\p I}&:= \max_{x\in \p I}|{\mathcal T}(x)|,& &\text{ and }& 
	\kappa_I&:= \max_{x\in  I}|{\mathcal T}(x)|
.\end{align*}
Then for $x\in \p I$, we have 
\begin{multline}t(|P_-y|^2(x)-|P_+y|^2(x))+\kappa_{\p I}|y|^2(x)\\
\leq \langle y,i\mathcal{T}y\rangle \leq t(|P_-y|^2(x)-|P_+y|^2(x))+\kappa_{\p I}|y|^2(x).
\end{multline} 
With this notation, integrating \eqref{gnu1} and \eqref{gnu2} over $I=[a,b]$ yields 
\begin{align}\label{gnu12}&    |P_-y|^2(a)
+e^{-\alpha_I t(b-a)} \frac{ 2-\alpha_I -\frac{\kappa_{\p I}}{t}}{2+\alpha_I -\frac{\kappa_{\p I}}{t} }|P_+y|^2(b)\nonumber\\
&\leq  \frac{ 2-\alpha_I +\frac{\kappa_{\p I}}{t}}{2+\alpha_I -\frac{\kappa_{\p I}}{t} }|P_+y|^2(a) 
 + e^{-\alpha_I t(b-a)}  \frac{ 2+\alpha_I +\frac{\kappa_{\p I}}{t}}{2+\alpha_I -\frac{\kappa_{\p I}}{t} }|P_-y|^2(b) , 
\end{align} 
and 
\begin{align}\label{gnu13}& |P_+y|^2(b)+e^{-\alpha_I t(b-a)}\frac{ 2-\alpha_I -\frac{\kappa_{\p I}}{t}}{2+\alpha_I -\frac{\kappa_{\p I}}{t} }|P_-y|^2  (a)
 \nonumber\\
&\leq  \frac{ 2-\alpha_I +\frac{\kappa_{\p I}}{t}}{2+\alpha_I -\frac{\kappa_{\p I}}{t} }|P_-y|^2  (b)+e^{-\alpha_R t(b-a)}\frac{ 2+\alpha_I +\frac{\kappa_{\p I}}{t}}{2+\alpha_I -\frac{\kappa_{\p I}}{t} }|P_+y|^2(a). 
\end{align} 
When $\kappa_I<<t$, we sharpen  \eqref{gnu12} and \eqref{gnu13} by choosing $\alpha_I:= 2(1-\frac{\kappa_I}{t})$ (using $|\frac{d y}{ds}|^2\geq (t-\kappa)^2|y|^2$).  We obtain

\begin{align}  \label{b12} |P_-y|^2(a)
 \leq  \frac{e^{ -\alpha_I t(b-a)}}{(1-\frac{\kappa_I}{t})}|P_-y|^2(b)+\frac{\kappa_I}{ t (1-\frac{\kappa_I}{t})}  |P_+y|^2(a) ,
\end{align}
and 
\begin{align} \label{b13}  |P_+y|^2(b)  
\leq   \frac{e^{-\alpha_I t(b-a)}}{(1-\frac{\kappa_I}{t}) }  |P_+y|^2(a) +\frac{\frac{\kappa_I}{t}}{(1-\frac{\kappa_I}{t})}  |P_-y|^2(b) .
\end{align}
Consider now 3 consecutive intervals: $I_1 = [\mu^\lambda_-,\lambda)$, $I_2=(\lambda, \frac{\mu^\lambda_++\lambda}{2}]$,  and $I_3 = [\frac{\mu^\lambda_++\lambda}{2},\mu^\lambda_+].$ Summing the 6 inequalities obtained by specializing \eqref{b12} and \eqref{b13}  to both $I_1$ and $I_3$, and 
\eqref{gnu12} and \eqref{gnu13} applied to $I_2$ yields for $t (\mu^\lambda_+-\lambda)$ and $t (\lambda-\mu^\lambda_-)$ large

\begin{multline}  \label{tout} 
|P_-y|^2(\mu^\lambda_-)  +|P_+y|^2(\mu^\lambda_+)\\
 \leq  ( \frac{\kappa_{I_1}}{ t (1-\frac{\kappa_{I_1}}{t})}  +  \frac{e^{-\alpha_{I_1} t(\lambda-\mu_L)}}{(1-\frac{\kappa_{I_1}}{t}) } ) |P_+y|^2(\mu^\lambda_-)   \\
 + (\frac{\frac{\kappa_{I_3}}{t}}{(1-\frac{\kappa_{I_3}}{t})}+\frac{e^{ -\alpha_{I_3} t( \frac{\mu_R-\lambda}{2})}}{(1-\frac{\kappa_{I_3}}{t})})|P_-y|^2(\mu^\lambda_+) .
\end{multline}

 We are now left to estimate $|E_my(\mu^\lambda_+)|$. Integrate the differential inequality  
\begin{align}\frac{d}{ds}(e^{-t(s-\mu^\lambda_+)}(s-\lambda)^{\frac{m-1}{2}}|E_my|(s))\leq e^{-t(s-\mu^\lambda_+)}(s-\lambda)^{\frac{m-1}{2}}|By|(s) 
\end{align}
to obtain 
\begin{align}\label{emu1} |E_my|(\mu^\lambda_+)&\leq \int_\lambda^{\mu^\lambda_+}e^{-t(u-\mu^\lambda_+)}(\frac{s-\lambda}{\mu^\lambda_+-\lambda })^{\frac{m-1}{2}}|By|(u)du \nonumber\\
&\leq \|B\|_{L^\infty}\sqrt{\mu^\lambda_+-\lambda}\sqrt{\int_\lambda^{\mu^\lambda_+}e^{-2t(u-\mu^\lambda_+)}(\frac{s-\lambda}{\mu^\lambda_+-\lambda })^{m-1}|y|^2(u)du} .
\end{align}
This leads us to choose the weight function  
$$\psi:= (s-\lambda)^{\frac{m-1 }{2}}e^{t(\mu^{\lambda}_+-s)}Q_1 + (s-\lambda+\frac{m-1}{ t})^{\frac{m-1}{2}}e^{t(\mu^{\lambda}_+-s)}Q_0 .$$ 
Then 
\begin{multline}\label{burr}
\sum_j| (\frac{iT_j}{2(s-\lambda)}+t_j) \psi y|^2-|[D,\psi] y|^2\\
= 
\sum_a [(\frac{w_a}{2(s-\lambda)}-t)^2+ \frac{(m-1)^2+2(m-1)-w_a^2}{4(s-\lambda)^2}]| \pi_aQ_1\psi y|^2\\
 - (\frac{m-1}{2(s-\lambda)}-t )^2|Q_1\psi y|^2 + t^2|Q_0\psi y|^2 \\ 
 -  (\frac{m-1}{2(s-\lambda+\frac{m-1}{t})}-t
 )^2|Q_0\psi y|^2 +O(t)|\psi y|^2
\\
= 
\sum_a [(\frac{t(m-1-w_a)}{(s-\lambda)}+ \frac{ 2(m-1) }{4(s-\lambda)^2}]| \pi_aQ_1\psi y|^2\\
  +   (2t -\frac{m-1}{2(s-\lambda)+\frac{m-1}{t}})(\frac{m-1}{2(s-\lambda)+\frac{m-1}{t}}) |Q_0\psi y|^2 +O(t)|\psi y|^2\\
\geq C\frac{t}{|\mu^\lambda_+-\lambda|}|\psi y|^2,
\end{multline}
where $\pi_a$ is projection onto the $w_a$ weight space, $w_a$ a weight. For the last inequality, we have made the additional assumption that $|\mu^\lambda_+-\lambda|>\frac{1}{t}.$ 

Combining \eqref{azalea} and \eqref{burr} gives us for   $|\mu^\lambda_{+}-\lambda|>>\frac{1}{t},$ and $C_1,C_2$ independent of $t$ and $\mu^\lambda_+$, for $t$ large, 
\begin{align}\label{podoph0}&\int_\lambda^{\mu^\lambda_{+}} e^{2t(\mu^\lambda_{+}-u)}(\frac{u-\lambda}{\mu^\lambda_{+}-\lambda})^{ m-1}|y|^2ds \leq C_1|\mu^\lambda_+-\lambda||P_-y(\mu^\lambda_{+})|^2 \nonumber\\
&+ C_2e^{2t(\mu^\lambda_{+}-\lambda)}t^{-m }|\mu^\lambda_+-\lambda|^{2-m}(|P_+y|^2(\lambda)-|P_-y|^2(\lambda)).
\end{align} 
To control the exponentially large term in \eqref{podoph0}, we apply \eqref{b13} for the interval $[\mu^\lambda_-,\lambda)$ to obtain
\begin{align} \label{b132}  |P_+y|^2(\lambda)  - |P_-y|^2(\lambda)
\leq   \frac{e^{-2(1-\frac{\kappa_{I_1}}{t}) t(\lambda-\mu^\lambda_-)}}{(1-\frac{\kappa_I}{t}) }  |P_+y|^2(\mu^\lambda_-). 
\end{align}
Inserting this estimates into \eqref{emu1} yields for some $C_3>0$, 
\begin{multline}\label{emu2} 
|E_my|(\mu^\lambda_+)\\
\leq    C_3((\mu^\lambda_+-\lambda)|P_-y(\mu^\lambda_{+})|   +  e^{t(\mu^\lambda_{+}+\mu^\lambda_--2\lambda)}t^{-\frac{m}{2} }|\mu^\lambda_+-\lambda|^{\frac{3-m}{2}}   |P_+y|(\mu^\lambda_-) ) .
\end{multline}
Proposition \ref{prop1} now follows readily from the estimates \eqref{tout} and \eqref{emu2}. 
\newline\hspace*{\fill}{\em q.e.d.}

\section{Conclusion}
By analyzing the Up transform in detail, we proved that the resulting connection is  anti-self-dual and has $L^2$ curvature; i.e. it is, indeed, an instanton.  We also computed the rank of the resulting bundle (Theorem~\ref{rankth}), the asymptotic form of the resulting connection (Theorem~\ref{MainTh}), and its Chern character values (Eqs.~\eqref{Eq:Ch1} and \eqref{Eq:Ch2})  in terms of the original bow representation.   In particular, the monopole charges $\hat{m}_j$ appearing in the asymptotic expression \eqref{goal} equal Chern numbers $R_{\lambda+}-R_{\lambda-}+\mathfrak{l}_{\lambda_j}$ of the  line bundles associated to the corresponding $\lambda$-points.  The basis of homology two cycles in $\mathrm{TN}_k^\nu$ is  $\{[C_{\sigma+1}]-[C_\sigma]\}_{\sigma=1}^{k-1},$ and the values of the first Chern classes on $C_\sigma$ is equal to the sum of $\sum_j\frac{\lambda_j}{\ell}$ and the index $R_{p_\sigma+} - R_{p_\sigma-} - r_{p_\sigma}$ of the corresponding $p$-point.  This is in complete agreement with the brane picture used in \cite{Witten:2009xu}.  

To conclude, our result now allows us to evaluate the index of the instanton  Dirac operator $D_A^-$.  For any instanton this index was computed in Theorem D of \cite{First}:
\begin{align}
\mathrm{ind}_{L^2}\, D_A^-=&\sum_j \left(
 (\left\{ \frac{\lambda_j}{\ell} \right\} -\frac12)  
 (m_j-k \lfloor\frac{\lambda_j}{\ell}\rfloor) 
 - \frac{k}{2}\left\{\frac{\lambda_j}{\ell}\right\}^2 
 \right)\\
 &+\frac{1}{8\pi^2}\int\mathrm{tr}\, F\wedge F.\nonumber
\end{align}
This expression is not evidently integer.  With our established value \eqref{Eq:Ch2} of the second Chern character, it becomes
\begin{align}
\mathrm{ind}_{L^2} D_A^-&= R_0 + 
\sum_j \left(
\frac{k}{2} \lfloor \frac{\lambda_j}{\ell}\rfloor(\lfloor \frac{\lambda_j}{\ell}\rfloor-1)
- \lfloor \frac{\lambda_j}{\ell}\rfloor(\hat{m}_j - k\lfloor \frac{\lambda_j}{\ell}\rfloor) \right). 
\end{align}
For the bow, all $\lambda_j\in[0,\ell)\Rightarrow \lfloor \frac{\lambda_j}{\ell}\rfloor=0. $  Hence, the connection resulting from the Up transform satisfies $\mathrm{ind}_{L^2} D_A^-= R_0.$  The index of the instanton Dirac operator therefore equals  the rank of the bow representation at $s=0.$

\begin{appendices}
\section{Anti-self-duality}\label{ASDproof}
The proof of the anti-self-duality of the curvature of the index bundle differs only slightly from the comparable computation in flat space \cite{Corrigan:1978ce,Nahm:1983sv}, as we now show. 
We recall that $\Pi_{t,b}$ denotes the orthogonal $L^2$ projection onto the kernel of $\D^\dagger_{t,b}$.
Then $\nabla_X = \Pi_{t,b}X^H\Pi_{t,b}.$ In local coordinates set 
$d^H:=e(dy^\alpha)(\frac{\p}{\p y^\alpha})^H ,$ where $e(\omega)$ denotes left exterior multiplication by $\omega$. Then 
\begin{align}F_{A}
&=  \Pi_{t,b}d^H \Pi_{t,b} d^H\Pi_{t,b} =  \Pi_{t,b}d^H (I-\D_{t,b}\mathfrak{G}_{t,b}\D_{t,b}^\dagger) d^H\Pi_{t,b}\nonumber\\
&=  \Pi_{t,b}(d^H)^2\Pi_{t,b}+\Pi_{t,b}[d^H ,\D_{t,b}]\mathfrak{G}_{t,b}[d^H,\D_{t,b}^\dagger ]\Pi_{t,b}.
\end{align}
Since the curvature of the line bundle $\mathcal{L}_v$ is anti-self-dual for each $v$, $(d^H)^2$ is anti-self-dual and therefore $\Pi_{t,b}(d^H)^2\Pi_{t,b}$ is anti-self-dual. In the coordinates $(\tau,t_1,t_2,t_3)$ we have 
\begin{align}[d^H,\D_{t,b}] = -\frac{\ii d\tau}{V} +\ii dt_j\otimes \e_j\text{ and }[d^H,\D_{t,b}^\dagger] =  \frac{\ii d\tau}{V} +\ii dt_j\otimes \e_j.
\end{align}
Hence $[d^H,\D_{t,b}]\wedge [d^H,\D_{t,b}^\dagger] $ is anti-self-dual. Because $\mathfrak{G}_{t,b}$ commutes with quaternions, we have 
$$\Pi_{t,b}[d^H ,\D_{t,b}]\mathfrak{G}_{t,b}[d^H,\D_{t,b}^\dagger ]\Pi_{t,b} = \Pi_{t,b}\mathfrak{G}_{t,b}[d^H ,\D_{t,b}]\wedge [d^H,\D_{t,b}^\dagger ]\Pi_{t,b}$$ is anti-self-dual. Thus $F_A$ is also anti-self-dual.

\section{\texorpdfstring{$\mathbb{R}^4$}{R4} in the Taub-NUT Coordinates}\label{App:R4}
\subsection{The Euclidean Metric}
We identify $\mathbb{R}^4\cong\mathbb{C}^2\cong \S$ and arrange its coordinates in a spinor $b$. For $\mathbb{R}^3=\mathrm{Im}\,\mathbb{H}$ we also arrange its coordinates into  $\tslash:=\ii\, \t.$ 
Given $bb^\dagger=t+\tslash$, $b^\dagger b=2t$ and
\begin{multline*}
2 d\vec{t}\,^2=\tr d\tslash d\tslash=\tr d(bb^\dagger-t)d(bb^\dagger-t)\\
=(b^\dagger db)^2+(db^\dagger b)^2+2 \tr db b^\dagger b db^\dagger-2 d t^2\\
=(b^\dagger db)^2+(db^\dagger b)^2+4t db^\dagger db-\frac12(db^\dagger b+b^\dagger db)^2\\
=4t db^\dagger db+\frac12(b^\dagger db-db^\dagger b)^2.
\end{multline*}
And the Euclidean metric on $\S\ni b$ is
\begin{align}
db^\dagger db=\frac{d\vec{t}\,^2}{2t}-2t\left(\frac{b^\dagger db-db^\dagger b}{4t}\right)^2
=\frac{d\vec{t}\,^2}{2t}+2t\hat{\eta}^2,
\end{align}
with the one-form $\hat{\eta}:=\ii \frac{b^\dagger db-db^\dagger b}{4t}$ with
\begin{align}
d\hat{\eta}=\ii\frac{db^\dagger\wedge db}{2t}-\ii\frac{db^\dagger b b^\dagger db}{4t^2}
=\ii\frac{db^\dagger(t-\tslash)db}{4 t^2}.
\end{align}
Considering that 
$
*_3 dt=-\frac{\ii}{4}\tr \frac{\tslash}{t} d\tslash d\tslash=-\frac{\ii}{2}db^\dagger(t-\tslash)db,
$
we have
\begin{align}
d\hat{\eta}=-\frac{*_3 dt}{2t^2}=*_3 d\frac{1}{2t}.
\end{align}

\subsection{The Intersection Number}
The identification $\S\cong\mathbb{C}^2\ni b=\left(\begin{smallmatrix}\bar{u}\\v\end{smallmatrix}\right)$ uses complex structure $\e_3=-\ii\sigma_3=\left(\begin{smallmatrix} -\ii&0\\0&\ii \end{smallmatrix}\right).$ The two coordinate planes, the $u$-plane $P_1:=\{v=0\}$ and the $v$-plane $P_2:=\{u=0\}$ have the usual orientation specified by the complex structure, in which, away from the origin, the frames $(u,\ii u)$ and $(v,\ii v)$ consisting of the radial vector and its counterclockwise normal, are positively oriented.  With these orientations, the intersection number of $P_1$ and $P_2$ is $+1.$ 

Under the map to $\mathbb{R}^3=\mathbb{C}\times\mathbb{R}$
\begin{align}
\pi: (u,v)\mapsto (t_1+\ii t_2, t_3)=(uv, \frac{u\bar{u}-v\bar{v}}{2}),
\end{align}
the fiber (away from the origin) is a circle.  The rotation along this circle is $b\mapsto e^{\ii\epsilon}b$, which in terms of $u$ and $v$ is $(u,v)\mapsto(e^{-\ii \epsilon}u, e^{\ii\epsilon}v).$  The Taub-NUT coordinates consist of $t_1,t_2,t_3,$ and a local coordinate $\tau$ along this circle fiber.  Under the above circle rotation $\tau\mapsto\tau+\epsilon.$

The image of the $u$-plane $\pi(P_1)=\{t_1=t_2=0, t_3\ge 0\}$ is the upward pointing ray, while the image of the $v$-plane $\pi(P_1)=\{t_1=t_2=0, t_3\le 0\}$ is the downward pointing ray.  From the Taub-NUT point of view $P_1$ and $P_2$ are two cigars originating at the same Taub-NUT center at the origin of $\mathbb{R}^3.$  The cigar orientation, however, is oriented by the frame consisting of the radial vector along the ray and  $\partial/\partial\tau$.  According to the circle action above, this corresponds to orienting $P_1$ with $(u,\ii u)$ frame and $P_2$ with $(v,\ii v)$ frame. Thus, in this orientation, the intersection number of the two cigars is $-1.$  

Note, that in the HHM compactification, cigars corresponding to the rays originating at the same NUT are homotopic.  Thus, the self-intersection of any such cigar is $-1.$

\end{appendices}

\bibliographystyle{alphaurl}
\bibliography{IDbib1}

\newcommand{\noop}[1]{}
\begin{thebibliography}{AHDM78}

\bibitem[AHDM78]{Atiyah:1978ri}
Michael~F. Atiyah, Nigel~J. Hitchin, Vladimir~G. Drinfel'd, and Yuri~I. Manin.
\newblock Construction of {I}nstantons.
\newblock {\em \href{http://dx.doi.org/10.1016/0375-9601(78)90141-X}{Phys.
  Lett. A}}, 65(3):185--187, 1978.

\bibitem[Ati57]{Atiyah57}
Michael~F. Atiyah.
\newblock Complex {A}nalytic {C}onnections in {F}ibre {B}undles.
\newblock {\em \href{https://doi.org/10.2307/1992969}{Trans. Amer. Math.
  Soc.}}, 85:181--207, 1957.

\bibitem[Biq97]{Biquard97}
Olivier Biquard.
\newblock Fibr\'{e}s de {H}iggs et {C}onnexions {I}nt\'{e}grables: le {C}as
  {L}ogarithmique (diviseur lisse).
\newblock {\em \href{https://doi.org/10.1016/S0012-9593(97)89915-6}{Ann. Sci.
  \'{E}cole Norm. Sup. (4)}}, 30(1):41--96, 1997.

\bibitem[BvB89]{Braam:1988qk}
Peter~J. Braam and Pierre van Baal.
\newblock Nahm's {T}ransformation for {I}nstantons.
\newblock {\em \href{http://projecteuclid.org/euclid.cmp/1104178397}{Comm.
  Math. Phys.}}, 122(2):267--280, 1989.

\bibitem[CFTG78]{Corrigan:1978ce}
Edward~F. {Corrigan}, David~B. {Fairlie}, Stephen {Templeton}, and Peter
  {Goddard}.
\newblock A {G}reen {F}unction for the {G}eneral {S}elf-{D}ual {G}auge {F}ield.
\newblock {\em \href{http://dx.doi.org/10.1016/0550-3213(78)90311-5}{Nuclear
  Physics B}}, 140:31--44, 1978.

\bibitem[CG84]{Corrigan:1983sv}
Edward~F. Corrigan and Peter Goddard.
\newblock Construction of {I}nstanton and {M}onopole {S}olutions and
  {R}eciprocity.
\newblock {\em \href{http://dx.doi.org/10.1016/0003-4916(84)90145-3}{Ann.
  Physics}}, 154(1):253--279, 1984.

\bibitem[CH19]{Cherkis:2017pop}
Sergey~A. Cherkis and Jacques Hurtubise.
\newblock Monads for instantons and bows.
\newblock {\em \href{https://doi.org/10.4310/ATMP.2019.v23.n1.a4}{Adv. Theor.
  Math. Phys.}}, 23(1):167--251, 2019.
\newblock \href {http://arxiv.org/abs/1709.00145} {\path{arXiv:1709.00145}}.

\bibitem[Che09]{Cherkis:2008ip}
Sergey~A. Cherkis.
\newblock Moduli {S}paces of {I}nstantons on the {T}aub-{NUT} {S}pace.
\newblock {\em \href{http://dx.doi.org/10.1007/s00220-009-0863-8}{Comm. Math.
  Phys.}}, 290(2):719--736, 2009.
\newblock \href {http://arxiv.org/abs/0805.1245} {\path{arXiv:0805.1245}}.

\bibitem[Che10]{Cherkis:2009jm}
Sergey~A. Cherkis.
\newblock Instantons on the {T}aub-{NUT} space.
\newblock {\em \href{https://projecteuclid.org/euclid.atmp/1288619154}{Adv.
  Theor. Math. Phys.}}, 14(2):609--641, 2010.
\newblock \href {http://arxiv.org/abs/0902.4724} {\path{arXiv:0902.4724}}.

\bibitem[Che11]{Cherkis:2010bn}
Sergey~A. Cherkis.
\newblock Instantons on {G}ravitons.
\newblock {\em \href{http://dx.doi.org/10.1007/s00220-011-1293-y}{Comm. Math.
  Phys.}}, 306(2):449--483, 2011.
\newblock \href {http://arxiv.org/abs/1007.0044} {\path{arXiv:1007.0044}}.

\bibitem[Che14]{Cherkis:2014vfa}
Sergey~A. Cherkis.
\newblock {Phases of Five-dimensional Theories, Monopole Walls, and Melting
  Crystals}.
\newblock {\em \href{http://dx.doi.org/10.1007/JHEP06(2014)027}{JHEP}}, 06:027,
  2014.
\newblock \href {http://arxiv.org/abs/1402.7117} {\path{arXiv:1402.7117}}.

\bibitem[CLHSip]{Third}
Sergey~A. Cherkis, Andr\'es Larra\'in-Hubach, and Mark Stern.
\newblock {I}nstantons on {M}ulti-{T}aub-{NUT} {S}paces {III}: {D}own
  {T}ransform.
\newblock {\em \noop{3000}}, in preparation\noop{3001ip}.

\bibitem[CLHS16]{First}
Sergey~A. Cherkis, Andr\'es Larra\'in-Hubach, and Mark Stern.
\newblock {Instantons on {M}ulti-{T}aub-{NUT} {S}paces I: {A}symptotic {F}orm
  and {I}ndex {T}heorem}.
\newblock {To appear in \em J. Diff. Geom.}
\newblock \href {http://arxiv.org/abs/1608.00018} {\path{arXiv:1608.00018}}.

\bibitem[Cro15]{Cross:2015hla}
Rebekah Cross.
\newblock {Asymptotic Dynamics of Monopole Walls}.
\newblock {\em \href{http://dx.doi.org/10.1103/PhysRevD.92.045029}{Phys.
  Rev.}}, D92(4):045029, 2015.
\newblock \href {http://arxiv.org/abs/1506.07606} {\path{arXiv:1506.07606}}.

\bibitem[CW12]{Cherkis:2012qs}
Sergey~A. Cherkis and Richard~S. Ward.
\newblock {Moduli of Monopole Walls and Amoebas}.
\newblock {\em \href{http://dx.doi.org/10.1007/JHEP05(2012)090}{JHEP}}, 05:090,
  2012.
\newblock \href {http://arxiv.org/abs/1202.1294} {\path{arXiv:1202.1294}}.

\bibitem[DK90]{DK1}
Simon~K. Donaldson and Peter~B. Kronheimer.
\newblock {\em The {G}eometry of {F}our-{M}anifolds}.
\newblock The Clarendon Press, Oxford University Press, New York, 1990.

\bibitem[EV86]{EV86}
H.~{Esnault} and E.~{Viehweg}.
\newblock {Logarithmic De Rham Complexes and Vanishing Theorems}.
\newblock {\em
  \href{http://adsabs.harvard.edu/abs/1986InMat..86..161E}{Inventiones
  Mathematicae}}, 86:161, 1986.
\newblock \href {https://doi.org/10.1007/BF01391499}
  {\path{doi:10.1007/BF01391499}}.

\bibitem[GRG97]{Gibbons:1996nt}
Gary~W. Gibbons, Paulina Rychenkova, and Ryushi Goto.
\newblock {HyperK\"ahler quotient construction of BPS monopole moduli spaces}.
\newblock {\em \href{https://doi.org/10.1007/s002200050121}{Commun. Math.
  Phys.}}, 186:585--599, 1997.
\newblock \href {http://arxiv.org/abs/hep-th/9608085}
  {\path{arXiv:hep-th/9608085}}.

\bibitem[HHM04]{HHM}
Tam{\'a}s Hausel, Eugenie Hunsicker, and Rafe Mazzeo.
\newblock Hodge {C}ohomology of {G}ravitational {I}nstantons.
\newblock {\em \href{http://dx.doi.org/10.1215/S0012-7094-04-12233-X}{Duke
  Math. J.}}, 122(3):485--548, 2004.
\newblock \href {http://arxiv.org/abs/math/0207169}
  {\path{arXiv:math/0207169}}.

\bibitem[Hit83]{Hitchin:1983ay}
Nigel~J. Hitchin.
\newblock On the {C}onstruction of {M}onopoles.
\newblock {\em \href{http://projecteuclid.org/euclid.cmp/1103922679}{Comm.
  Math. Phys.}}, 89(2):145--190, 1983.

\bibitem[HKLR87]{Hitchin:1986ea}
Nigel~J. Hitchin, Anders Karlhede, Ulf Lindstr{{\"o}}m, and Martin.
  Ro{\v{c}}ek.
\newblock Hyper-{K}{\"a}hler {M}etrics and {S}upersymmetry.
\newblock {\em \href{http://projecteuclid.org/euclid.cmp/1104116624}{Comm.
  Math. Phys.}}, 108(4):535--589, 1987.

\bibitem[KN90]{KN}
Peter~B. Kronheimer and Hiraku Nakajima.
\newblock Yang-{M}ills {I}nstantons on {ALE} {G}ravitational {I}nstantons.
\newblock {\em \href{http://dx.doi.org/10.1007/BF01444534}{Math. Ann.}},
  288(2):263--307, 1990.

\bibitem[Nah80]{Nahm:1979yw}
Werner Nahm.
\newblock {A Simple {F}ormalism for the {BPS} {M}onopole}.
\newblock {\em \href{http://dx.doi.org/10.1016/0370-2693(80)90961-2}{Phys.
  Lett.}}, B90:413, 1980.

\bibitem[Nah82]{NahmADHM}
Werner Nahm.
\newblock The {C}onstruction of all {S}elf-{D}ual {M}ultimonopoles by the
  {ADHM} {M}ethod.
\newblock In {\em Monopoles in {Q}uantum {F}ield {T}heory ({T}rieste, 1981)},
  pages 87--94. World Sci. Publishing, Singapore, 1982.

\bibitem[Nah84]{Nahm:1983sv}
Werner Nahm.
\newblock Self-{D}ual {M}onopoles and {C}alorons.
\newblock In {\em Group Theoretical Methods in Physics ({T}rieste, 1983)},
  volume 201 of {\em \href{http://dx.doi.org/10.1007/BFb0016145}{Lecture Notes
  in Phys.}}, pages 189--200. Springer, Berlin, 1984.

\bibitem[Nak93]{Nakajima:1990zx}
Hiraku Nakajima.
\newblock Monopoles and {N}ahm's {E}quations.
\newblock In {\em Einstein Metrics and {Y}ang-{M}ills Connections ({S}anda,
  1990)}, volume 145 of {\em Lecture Notes in Pure and Appl. Math.}, pages
  193--211. Dekker, New York, 1993.

\bibitem[Sza07]{Szabo07}
Szil\'ard Szab\'o.
\newblock Nahm {T}ransform for {I}ntegrable {C}onnections on the {R}iemann
  {S}phere.
\newblock {\em M\'em. Soc. Math. Fr. (N.S.)}, (110):ii+114 pp. (2008), 2007.
\newblock \href {http://arxiv.org/abs/arXiv:math/0511471}
  {\path{arXiv:arXiv:math/0511471}}.

\bibitem[Sza17]{Szabo17}
Szil\'{a}rd Szab\'{o}.
\newblock Nahm transformation for parabolic {H}iggs bundles on the projective
  line---case of non-semisimple residues.
\newblock {\em \href{https://doi.org/10.1016/j.geomphys.2016.10.018}{J. Geom.
  Phys.}}, 122:80--91, 2017.

\bibitem[Wit09]{Witten:2009xu}
Edward Witten.
\newblock Branes, {I}nstantons, and {T}aub-{NUT} {S}paces.
\newblock {\em \href{http://dx.doi.org/10.1088/1126-6708/2009/06/067}{J. High
  Energy Phys.}}, 6:067, 55, 2009.
\newblock \href {http://arxiv.org/abs/0902.0948} {\path{arXiv:0902.0948}}.

\end{thebibliography}

\end{document}